\documentclass[a4paper,12pt,twoside]{article}
\usepackage[francais]{babel}
\usepackage{amssymb, amsmath, eucal}
\usepackage[all]{xy}
\usepackage{mathrsfs}
\usepackage[dvips]{graphics}
\begin{document}
\textwidth15.5cm
\textheight22.5cm
\voffset=-13mm
\newtheorem{The}{Theorem}[section]
\newtheorem{Lem}[The]{Lemma}
\newtheorem{Prop}[The]{Proposition}
\newtheorem{Cor}[The]{Corollary}
\newtheorem{Rem}[The]{Remark}
\newtheorem{Obs}[The]{Observation}
\newtheorem{SConj}[The]{Standard Conjecture}
\newtheorem{Titre}[The]{\!\!\!\! }
\newtheorem{Conj}[The]{Conjecture}
\newtheorem{Question}[The]{Question}
\newtheorem{Prob}[The]{Problem}
\newtheorem{Def}[The]{Definition}
\newtheorem{Not}[The]{Notation}
\newtheorem{Claim}[The]{Claim}
\newtheorem{Conc}[The]{Conclusion}
\newtheorem{Ex}[The]{Example}
\newtheorem{Fact}[The]{Fact}
\newcommand{\C}{\mathbb{C}}
\newcommand{\R}{\mathbb{R}}
\newcommand{\N}{\mathbb{N}}
\newcommand{\Z}{\mathbb{Z}}
\newcommand{\Q}{\mathbb{Q}}
\newcommand{\Proj}{\mathbb{P}}

\begin{center}

{\Large\bf Holomorphic Deformations of Balanced Calabi-Yau $\partial\bar\partial$-Manifolds}

\end{center}

\begin{center}

 {\large Dan Popovici}

\end{center}

\vspace{1ex}

\noindent {\small {\bf Abstract.} Given a compact complex $n$-fold $X$ satisfying the $\partial\bar\partial$-lemma and supposed to have a trivial canonical bundle $K_X$ and to admit a balanced (=semi-K\"ahler) Hermitian metric $\omega$, we introduce the concept of deformations of $X$ that are {\bf co-polarised} by the balanced class $[\omega^{n-1}]\in H^{n-1,\,n-1}(X,\,\C)\subset H^{2n-2}(X,\,\C)$ and show that the resulting theory of balanced co-polarised deformations is a natural extension of the classical theory of K\"ahler polarised deformations in the context of Calabi-Yau or even holomorphic symplectic compact complex manifolds. The concept of Weil-Petersson metric still makes sense in this strictly more general, possibly non-K\"ahler context, while the Local Torelli Theorem still holds.}

\vspace{3ex}

\section{Introduction}\label{section:introduction}

  Let $X$ be a compact complex manifold ($\mbox{dim}_{\C}X=n$). Recall that a Hermitian metric $\omega$ on $X$ (identified throughout with the corresponding $C^{\infty}$ positive-definite $(1,\,1)$-form $\omega$) is said to be {\it balanced} (see [Gau77b] for the actual notion called {\it semi-K\"ahler} there, [Mic82] for the actual term) if $$d(\omega^{n-1})=0,$$  \noindent while $X$ is said to be a {\it balanced manifold} if it carries such a metric. In dimension $n\geq 3$, the balanced condition on $X$, besides being weaker than the K\"ahler one, is even weaker than the Fujiki {\it class} ${\cal C}$ condition thanks to a theorem of Alessandrini and Bassanelli [AB93]. (Recall that a compact complex manifold $X$ is said to be of {\it class} ${\cal C}$ if it is bimeromorphically equivalent to a compact K\"ahler manifold, i.e. if there exists a holomorphic bimeromorphic map $\mu:\widetilde{X}\rightarrow X$, called modification, from a compact K\"ahler manifold $\widetilde{X}$.) 

 On the other hand, all {\it class} ${\cal C}$ compact complex manifolds are known to satisfy the {\it $\partial\bar\partial$-lemma} (cf. [DGMS75]) in the following sense\!\!: \\

\noindent {\it for every pure-type $d$-closed form on $X$, the properties of $d$-exactness, $\partial$-exactness, $\bar\partial$-exactness and $\partial\bar\partial$-exactness are equivalent}. \\

\noindent There exist compact complex manifolds satisfying the $\partial\bar\partial$-lemma that are not of {\it class} ${\cal C}$ (see e.g. Observation \ref{Obs:bal-ddbar-nonC}), while the $\partial\bar\partial$-lemma implies the Hodge decomposition and the Hodge symmetry on $X$, i.e. it defines canonical isomorphisms (the latter one by conjugation)\!\!: $$H^k_{DR}(X,\, \C)\simeq\bigoplus\limits_{p+q=k}H^{p,\,q}(X,\,\C) \hspace{2ex}\mbox{and}\hspace{2ex} H^{q,\,p}(X,\,\C)\simeq\overline{H^{p,\,q}(X,\,\C)},$$ \noindent relating the De Rham cohomology groups $H^k_{DR}(X,\, \C)$ ($k=0, 1, \dots , 2n$) to the Dolbeault cohomology groups $H^{p,\,q}(X,\,\C)$ ($p,q=0, 1, \dots , n$). These manifolds, sometimes referred to as {\it cohomologically K\"ahler}, will play a key role in this work and thus deserve a name in their own right.

\begin{Def}\label{Def:d-dbar} A compact complex manifold $X$ will be said to be a {\bf $\partial\bar\partial$-manifold} if the $\partial\bar\partial$-lemma holds on $X$.

 If, furthermore, the canonical bundle $K_X$ of $X$ is trivial, $X$ will be called a {\bf Calabi-Yau $\partial\bar\partial$-manifold}.

\end{Def}

While there are plenty of examples of compact balanced manifolds that do not satisfy the $\partial\bar\partial$-lemma (e.g. the Iwasawa manifold), the answer to the question of whether non-balanced $\partial\bar\partial$-manifolds exist does not seem to be known and constitutes in our opinion a problem worth investigating. Part of the difficulty stems from the fact that neither twistor spaces (which are always balanced by [Gau91]), nor nilmanifolds (which are never $\partial\bar\partial$ unless they are K\"ahler), nor any other familiar class of compact complex non-K\"ahler manifolds (e.g. the Calabi-Eckmann manifolds are never either balanced or $\partial\bar\partial$) can produce an example. 

 Our main object of study in this paper will be the class of {\it balanced} $\partial\bar\partial$-{\it manifolds}. The principal reason behind our interest in them stems from their remarkable stability properties under both modifications and small deformations. Indeed, if $\mu\,\,:\,\,\widetilde{X}\rightarrow X$ is a modification of compact complex manifolds, Corollary 5.7 in [AB96] (see also [AB95]) states that $$X \hspace{1ex}\mbox{\it is balanced if and only if}\hspace{1ex} \widetilde{X} \hspace{1ex}\mbox{\it is balanced},$$ \noindent while by [DGMS75] we know that $X$ is a $\partial\bar\partial$-manifold whenever $\widetilde{X}$ is one. On the other hand, although the balanced property of $X$ is not open under holomorphic deformations by another result of Alessandrini and Bassanelli (cf. [AB90]), the $\partial\bar\partial$-{\it property} and the simultaneous occurence of the {\it balanced} and $\partial\bar\partial$-{\it properties} are both deformation open by two results of Wu (cf. [Wu06], a partial survey of which can be found in [Pop11, 4.3]). Specifically, if $(X_t)_{t\in\Delta}$ is a holomorphic family of compact complex manifolds, C.-C. Wu proved in [Wu06] the following two theorems\!\!: 

\vspace{1ex}

$(i)$\, {\it if $X_0$ is a $\partial\bar\partial$-manifold, then $X_t$ is a $\partial\bar\partial$-manifold for all $t\in\Delta$ sufficiently close to $0\in\Delta$;} 

\vspace{1ex}

$(ii)$\, {\it if $X_0$ is a balanced $\partial\bar\partial$-manifold, then $X_t$ is a balanced $\partial\bar\partial$-manifold for all $t\in\Delta$ sufficiently close to $0\in\Delta$.}

\vspace{1ex}

 This paper will hopefully provide further evidence to substantiate the view that non-K\"ahler, {\it balanced Calabi-Yau $\partial\bar\partial$-manifolds} form a class that is well worth studying. We point out examples of such manifolds in $\S.$\ref{subsection:examples}.

 Our first observation (cf. section \ref{section:BTT}) will be that the K\"ahler assumption can be weakened to the $\partial\bar\partial$ assumption in the Bogomolov-Tian-Todorov theorem. This fact, hinted at in the introduction to [Tia87], is probably known, but we take this opportunity to point out how the $\partial\bar\partial$-lemma can be solely relied upon in the proofs given in [Tia87] and [Tod89].

\begin{The}(Bogomolov-Tian-Todorov for Calabi-Yau $\partial\bar\partial$-manifolds) \label{The:BTT} Let $X$ be a compact complex manifold satisfying the $\partial\bar\partial$-lemma and whose canonical bundle $K_X$ is trivial. Then the Kuranishi family of $X$ is unobstructed.

\end{The}

 Here, as usual, unobstructedness means that the base space of the Kuranishi family is isomorphic to an open subset of $H^{0, \, 1}(X, \, T^{1, \, 0}X)$. \\

 Given a compact {\it balanced Calabi-Yau $\partial\bar\partial$-manifold} $X$ of complex dimension $n$, by a {\it balanced class} $[\omega^{n-1}]\in H^{n-1,\,n-1}(X,\,\C)\subset H^{2n-2}(X,\,\C)$ we shall mean the Dolbeault cohomology class of type $(n-1,\, n-1)$ (or the De Rham\footnote{Dolbeault cohomology classes will be denoted by $[\hspace{2ex}]$, De Rham cohomology classes will be denoted by $\{\,\,\,\}$.} cohomology class of degree $2n-2$ that is the image of the former under the above canonical inclusion which holds thanks to the $\partial\bar\partial$ assumption -- see e.g. Lemma \ref{Lem:d-closed-rep}) of the $(n-1)^{st}$ power of a balanced metric $\omega$ on $X$. If $\Delta\subset H^{0, \, 1}(X, \, T^{1, \, 0}X)$ denotes the open subset that is the base space of the Kuranishi family $(X_t)_{t\in\Delta}$ of $X$ (with $X_0=X$), we define in section \ref{section:co-polarisations} the notion of local deformations $X_t$ of $X$ that are {\bf co-polarised by the balanced class $[\omega^{n-1}]$} by requiring that the De Rham cohomology class $\{\omega^{n-1}\}\in H^{2n-2}(X,\,\C)$ be of type $(n-1,\, n-1)$ for the complex structure $J_t$ of $X_t$. Since all nearby fibres $X_t$ are again $\partial\bar\partial$-manifolds by Wu's first theorem [Wu06], the De Rham cohomology space $H^{2n-2}(X,\,\C)$ (which is independent of the complex structure $J_t$) admits for all $t\in\Delta$ a Hodge decomposition \begin{eqnarray}\label{eqn:Hodge-decomp2n-2}H^{2n-2}(X,\,\C)=H^{n,\,n-2}(X_t,\,\C)\oplus H^{n-1,\,n-1}(X_t,\,\C)\oplus H^{n-2,\,n}(X_t,\,\C)\end{eqnarray} \noindent depending on the complex structure $J_t$ and satisfying the Hodge symmetry $H^{n,\,n-2}(X_t,\,\C)\simeq\overline{H^{n-2,\,n}(X_t,\,\C)}$ (after possibly shrinking $\Delta$ about $0$). The original balanced class $\{\omega^{n-1}\}\in H^{2n-2}(X,\,\R)$ is said to be of type $(n-1,\, n-1)$ for the complex structure $J_t$ of $X_t$ if its components of types $(n,\,n-2)$ and $(n-2,\,n)$ in the Hodge decomposition (\ref{eqn:Hodge-decomp2n-2}) vanish. The balanced class being real, this is equivalent to either of these components vanishing. The condition is still equivalent to the De Rham class $\{\omega^{n-1}\}\in H^{2n-2}(X,\,\C)$ being representable by a form of $J_t$-pure type $(n-1,\, n-1)$. We denote by $$\Delta_{[\omega^{n-1}]}\subset\Delta \hspace{3ex}  \mbox{and by}\hspace{3ex} \pi:{\cal X}_{[\omega^{n-1}]}\rightarrow\Delta_{[\omega^{n-1}]}$$ \noindent the open subset of local deformations of $X$ co-polarised by the balanced class $[\omega^{n-1}]$, resp. the local universal family of co-polarised deformations of $X$. Since all sufficiently nearby fibres $X_t$ are again balanced $\partial\bar\partial$-manifolds by Wu's second theorem [Wu06], Observation \ref{Obs:balanced-def-class_appendix} shows that the co-polarising De Rham cohomology class $\{\omega^{n-1}\}\in H^{2n-2}(X,\,\C)$ can still be represented by a form $\omega_t^{n-1}$ with $\omega_t$ a $J_t$-balanced metric for every $t\in\Delta_{[\omega^{n-1}]}$ (after possibly shrinking $\Delta$ about $0$). 
  
  We go on to show in $\S.$\ref{subsection:comparisonK} that in the special case where $X$ is K\"ahler and $\omega$ is a K\"ahler metric on $X$, the local deformations of $X$ that are co-polarised by the balanced class $[\omega^{n-1}]\in H^{n-1,\, n-1}(X,\,\C)$ are precisely those that are polarised by the K\"ahler class $[\omega]\in H^{1,\, 1}(X,\,\C)$. Thus the theory of balanced co-polarisations is a natural extension to the balanced case of the classical theory of K\"ahler polarisations.

 The tangent space to $\Delta_{[\omega^{n-1}]}$ at $0$ is isomorphic under the Kodaira-Spencer map to a subspace $H^{0,\,1}(X,\,T^{1,\, 0}X)_{[\omega^{n-1}]}\subset H^{0,\,1}(X,\,T^{1,\, 0}X)$ which, in turn, is isomorphic under the canonical isomorphism defined by the Calabi-Yau property of $X$ to a subspace $H^{n-1,\, 1}_{prim}(X,\,\C)\subset H^{n-1,\, 1}(X,\,\C)$\!\!: $$T_0\Delta_{[\omega^{n-1}]}\simeq H^{0,\,1}(X,\,T^{1,\, 0}X)_{[\omega^{n-1}]}\simeq H^{n-1,\, 1}_{prim}(X,\,\C).$$ \noindent It is well known that in the case of deformations polarised by a K\"ahler class $[\omega]$, $H^{n-1,\, 1}_{prim}(X,\,\C)$ is the space of {\it primitive} (for the K\"ahler class $[\omega]$) Dolbeault classes of type $(n-1,\, 1)$. In the balanced case, since $\omega$ need not be closed, the standard definition of {\it primitive} $(n-1,\, 1)$-classes is no longer meaningful, but we use the analogy with the K\"ahler case to define $H^{n-1,\, 1}_{prim}(X,\,\C)$ in an ad hoc way (cf. Definition \ref{Def:primitive-classes_n-11}). 

 In section \ref{section:period-WP}, we give two applications of this construction. The Calabi-Yau condition being deformation open when the Hodge number $h^{n,\,0}(t):=\mbox{dim}_{\C}H^{n,\,0}(X_t,\,\C)$ does not jump in a neighbourhood of $t=0$ (and this is indeed the case here since we assume $X_0$ to be a $\partial\bar\partial$-manifold -- the non-jumping actually holds more generally if we merely assume the Fr\"olicher spectral sequence of $X_0$ to degenerate at $E_1$), $K_{X_t}$ is still trivial for all $t\in\Delta$ sufficiently close to $0$, so when the complex structure $J_t$ of $X$ varies, we get complex lines $$H^{n,\, 0}(X_t,\,\C)\subset H^n(X,\,\C),  \hspace{3ex} t\in\Delta,$$ \noindent varying in a fixed complex vector space in a holomorphic way with $t\in\Delta$. (The above inclusion follows from the $\partial\bar\partial$-lemma holding on $X_t$ for all $t$ close to $0$ by Wu's first theorem [Wu06].) We show (cf. Theorem \ref{The:localTorelli}) that the resulting (holomorphic) period map $$\Delta\ni t\stackrel{{\cal P}}{\longmapsto} H^{n,\, 0}(X_t,\,\C)\in\Proj H^n(X,\,\C)$$ \noindent in which every complex line $H^{n,\, 0}(X_t,\,\C)$ has been identified with a point in the complex projective space $\Proj H^n(X,\,\C)$ is locally an immersion. This means that the {\bf Local Torelli Theorem} still holds in this context. 

 We also propose two variants $\omega^{(1)}_{WP}, \omega^{(2)}_{WP}$ (cf. Definition \ref{Def:WP1}) of a notion of {\bf Weil-Petersson metric} on $\Delta_{[\omega^{n-1}]}$ associated with a $C^{\infty}$ family of {\it balanced metrics} $(\omega_t)_{t\in\Delta_{[\omega^{n-1}]}}$ on the fibres $X_t$ such that each $\omega_t^{n-1}$ lies in the co-polarising balanced class $\{\omega^{n-1}\}\in H^{2n-2}(X,\,\C)$. The metrics $\omega^{(1)}_{WP}$ and $\omega^{(2)}_{WP}$ coincide if $\mbox{Ric}(\omega_t)=0$ for all $t\in\Delta_{[\omega^{n-1}]}$ (cf. Observation \ref{Obs:WP-metrics-coincide}). Although in the case of K\"ahler polarised deformations, the Weil-Petersson metric $\omega_{WP}$ coincides with the pullback ${\cal P}^{\star}\omega_{FS}$ of the Fubini-Study metric of $\Proj H^n(X,\,\C)$ under the period map (hence $\omega_{WP}$ is K\"ahler) by Tian's Theorem 2 in [Tia87], this need not be the case in our context of balanced co-polarised deformations. However, we can compare these two metrics (cf. Theorem \ref{The:final-metric-formulae} and Corollary \ref{Cor:metric-comparison}), show that $$\omega^{(2)}_{WP}\geq{\cal P}^{\star}\omega_{FS}>0 \hspace{3ex} \mbox {on}\hspace{1ex}\Delta_{[\omega^{n-1}]}$$ \noindent and make the difference $\omega^{(2)}_{WP}-{\cal P}^{\star}\omega_{FS}$ explicit. The obstruction to the identity $\omega^{(2)}_{WP}={\cal P}^{\star}\omega_{FS}$ holding appears now clearly to be the possible non-existence of a form that is both {\it primitive} and {\it $d$-closed} in an arbitrary Dolbeault cohomology class of type $(n-1,\, 1)$ which is assumed to be {\it primitive} in the ad hoc balanced sense of section \ref{section:co-polarisations} (cf. $\S.$\ref{subsection:prim-n-11}). 

 It is natural to ask whether the Weil-Petersson metric $\omega^{(2)}_{WP}$ is K\"ahler. We ignore the answer at this stage, although we cannot see why this should be the case for balanced, non-K\"ahler fibres $X_t$.

 Another natural question is whether there is a canonical choice of balanced metrics $\omega_t$ on the fibres $X_t$ such that each $\omega_t^{n-1}$ lies in the co-polarising balanced class $\{\omega^{n-1}\}\in H^{2n-2}(X,\,\C)$ and such that $\mbox{Ric}(\omega_t)=0$ for all $t$. This would induce a {\it canonical} (i.e. depending only on the co-polarising balanced class $\{\omega^{n-1}\}$) Weil-Petersson metric $\omega_{WP}$ on $\Delta_{[\omega^{n-1}]}$. The answer would follow from the answer to another tantalising question: \begin{Question}\footnote{Since the first version of this paper was posted on the arXiv, an analogue of this problem for Gauduchon metrics in the context of the Aeppli cohomology has been introduced in [Pop15], also considered in [TW18] and completely solved in [STW17], but the balanced case put forward here currently seems out of reach.} Does there exist a balanced analogue of Yau's theorem on the Calabi conjecture? In other words, is it true that every balanced class contains the $(n-1)^{st}$ power of a balanced metric for which the volume form has been prescribed?

\end{Question}

In the last section \ref{section:hol-symp} we briefly discuss the case of co-polarised deformations of holomorphic symplectic manifolds.

\section{Preliminaries}\label{section:preliminaries}

\subsection{Examples of non-K\"ahler, balanced Calabi-Yau $\partial\bar\partial$-manifolds}\label{subsection:examples}

 We now pause to point out a few examples in support of the theory that will be developed in the next sections. 

\vspace{1ex}

$(1)$\,  We have seen that all {\it class} ${\cal C}$ manifolds are both balanced and $\partial\bar\partial$. However, the implications are strict and even the simultaneous occurence of the balanced and $\partial\bar\partial$ conditions does not ensure the {\it class} ${\cal C}$ property.

\begin{Obs}\label{Obs:bal-ddbar-nonC} There exist compact balanced $\partial\bar\partial$-manifolds that are not of class ${\cal C}$. In other words, the class of compact {\bf balanced $\partial\bar\partial$-manifolds} strictly contains Fujiki's class ${\cal C}$.

\end{Obs}

\noindent {\it Proof.} To show that the {\it class} ${\cal C}$ property of compact complex manifolds is not deformation open, [Cam91a] and [LP92] exhibit holomorphic families of twistor spaces $(X_t)_{t\in\Delta}$ in which the central fibre $X_0$ is Moishezon (hence is also a $\partial\bar\partial$-manifold), while, for every $t\in\Delta\setminus\{0\}$ sufficiently close to $0$, the fibre $X_t$ has vanishing algebraic dimension (hence is non-Moishezon, hence is not of {\it class} ${\cal C}$ since, by another result of Campana [Cam91b], the Moishezon and class ${\cal C}$ properties of twistor spaces are equivalent). However, all the fibres $X_t$ are {\it balanced} since, by a result of Gauduchon [Gau91], all twistor spaces are balanced. Moreover, when $t$ is close to $0$, $X_t$ is a $\partial\bar\partial$-manifold by Wu's openness theorem for deformations of $\partial\bar\partial$-manifolds. Thus any of the fibres $X_t$ with $t\neq 0$ but $t$ close to $0$ provides an example as stated.  \hfill $\Box$

\vspace{1ex}

 Notice that the above examples are not Calabi-Yau manifolds since the restriction of the canonical bundle of any twistor space to any twistor line is isomorphic to ${\cal O}_{\Proj^1}(-4)$, hence it cannot be trivial. 

\vspace{1ex}

$(2)$\, On the other hand, examples of compact {\bf non-K\"ahler, class ${\cal C}$}, {\bf holomorphic symplectic manifolds} were constructed by Yoshioka in [Yos01, section 4.4]. In particular, Yoshioka's manifolds are compact, non-K\"ahler, balanced Calabi-Yau $\partial\bar\partial$-manifolds. Thus they fall into the category of manifolds that will be investigated in this paper. While Yoshioka's manifolds are of {\it class} ${\cal C}$, it is natural to wonder whether compact, {\it non-class} ${\cal C}$, balanced Calabi-Yau $\partial\bar\partial$-manifolds (i.e. manifolds as in Observation \ref{Obs:bal-ddbar-nonC} having, in addition, a trivial canonical bundle) exist. They actually do (see $(3)$ and $(4)$ below) and it is tempting to expect that such an example could be constructed by taking our cue from e.g. [FLY12]\!: by starting off with a compact K\"ahler Calabi-Yau manifold $X$, contracting $X$ under a crepant map to some (possibly non-K\"ahler, but necessarily {\it class} ${\cal C}$) manifold $Y$ and then slightly deforming $Y$ to some (possibly {\it non-class} ${\cal C}$, but necessarily balanced Calabi-Yau $\partial\bar\partial$) manifold $Y_t$. The stability properties of balanced $\partial\bar\partial$-manifolds under both contractions and small deformations could thus be taken full advantage of if an explicit example obtained in this way (in which $Y$ is smooth) could be written down.

\vspace{1ex}

$(3)$\, We will now point out a class of examples, recently constructed in the literature, of compact {\bf balanced Calabi-Yau $\partial\bar\partial$-manifolds} that are {\bf not of class ${\cal C}$}. In [FOU14, Theorem 5.2] (see also [Kas13]), a (compact) solvmanifold $M$ of real dimension $6$ and a holomorphic family of complex structures $(J_a)_{a\in\Delta}$ on $M$ are constructed (where $\Delta:=\{a\in\C\,;\,|a|<1\}$) such that $X_a:=(M,\,J_a)$ is a balanced Calabi-Yau $\partial\bar\partial$-manifold for every $a\in\Delta\setminus\{0\}$. Furthermore, it can be easily checked that $X_a=(M,\,J_a)$ is not of class ${\cal C}$ for any $a\in\Delta$ by either of the next two arguments. 

\vspace{1ex}

$(a)$\, A direct calculation shows the existence of a $C^{\infty}$ positive definite $(1,\,1)$-form $\omega$ on $X_a$ such that $i\partial\bar\partial\omega\geq 0$. Then, by Theorem 2.3 in [Chi14], if $X_a$ were of class ${\cal C}$, it would have to be K\"ahler. However, a direct calculation shows that no K\"ahler metrics exist on any $X_a$. This argument has kindly been communicated to the author by L. Ugarte. 

\vspace{1ex}

$(b)$\, Since the fundamental group is a bimeromorphic invariant of compact complex manifolds, if $X_a$ were of class ${\cal C}$, its fundamental group would also occur as the fundamental group of a compact K\"ahler manifold. However, this is impossible as follows from [Cam04] (where it is proved that the Albanese morphism $\alpha_X \,\,:\,\,X\to\mbox{Alb}(X)$ of any Calabi-Yau class ${\cal C}$ manifold $X$ is {\it surjective}) combined with [Cam95] (where $\pi_1(X)$ is studied when $\alpha_X$ is surjective). This argument has kindly been communicated to the author by F. Campana.

\vspace{1ex}

$(4)$\, Let us finally mention that in the very recent preprint [Fri17], a large class of compact, {\it non-class} ${\cal C}$, balanced Calabi-Yau $\partial\bar\partial$-manifolds obtained via a construction of Clemens's (that was subsequently used by many authors, including Friedman himself and in [FLY12]) was produced. The original compact $3$-fold $X$ with trivial $K_X$ is only assumed to be $\partial\bar\partial$ (so not necessarily K\"ahler), to have $h^{0,\,1}=h^{0,\,2}=0$ (hence also $h^{1,\,0}=h^{2,\,0}=0$ thanks to the Hodge symmetry that holds on any $\partial\bar\partial$-manifold) and to have disjoint smooth rational $(-1,\,-1)$-curves $C_1,\dots , C_r$ whose classes $[C_i]\in H^4(X,\,\C)$ satisfy a linear dependence relation but generate $H^4(X,\,\C)$. Friedman shows in [Fri17] that the singular compact $3$-fold obtained from any such $X$ by contracting the $C_i$'s has smooth small deformations that are $\partial\bar\partial$-manifolds, but are not of {\it class} ${\cal C}$. They are not even deformation equivalent to any {\it class} ${\cal C}$ manifold. Actually, all the small deformations lying in an open dense subset of the moduli space are shown in [Fri17] to be $\partial\bar\partial$-manifolds, while all small smoothings are conjectured to be. They are balanced by [FLY12].

\subsection{The  balanced Ricci-flat Bochner principle}\label{section:Bochner-balRflat}\label{subsection:Bochner-principle}

 We collect here some known facts that will come in handy later on.

 \begin{Prop}\label{Prop:Bprinc-bal} Let $(X,\,\omega)$ be a compact complex Hermitian manifold with $\mbox{dim}_{\C}X=n$. Let $D=D'+\bar\partial$ be the Chern connection of $K_X$ equipped with the metric induced by $\omega$ and let $|\,\cdot\,|_{\omega}$ (resp. $||\,\cdot\,||$) be the pointwise norm (resp. the $L^2$-norm) of $K_X$-valued forms w.r.t. this metric.

  If $\omega$ is {\bf balanced} and $\mbox{Ric}\,(\omega)=0$, every $u\in C^{\infty}_{n,\, 0}(X,\,\C)\simeq C^{\infty}(X,\,K_X)$ satisfies $$||\bar\partial u||^2=||D'u||^2.$$ \noindent In particular, every holomorphic $n$-form $u$ on $X$ (i.e. $u\in H^0(X,\,K_X)$) is {\bf parallel} (i.e. $Du=0$) and satisfies \begin{eqnarray}\label{eqn:uubar-parallel}|u|_{\omega}^2=C \hspace{3ex} (\mbox{hence also} \hspace{2ex} i^{n^2}\,u\wedge\bar{u} = C\,\omega^n)  \hspace{2ex}\mbox{on}\hspace{1ex} X\end{eqnarray} \noindent for some constant $C\geq 0$.

\end{Prop}

 It turns out that the balanced assumption is unnecessary in the last statement of Proposition \ref{Prop:Bprinc-bal}. Indeed, every holomorphic section of a flat line bundle is parallel.

\begin{Obs}\label{Obs:flat-parallel} Let $(L,\, h)\rightarrow X$ be a Hermitian holomorphic line bundle over a compact complex manifold such that the curvature form $i\Theta_h(L)$ vanishes identically on $X$. Then any global holomorphic section $\sigma\in H^0(X,\, L)$ satisfies $D\sigma=0$, where $D$ is the Chern connection of $(L,\, h)$.

\end{Obs}

Both the statement and the proof of this observation have been kindly pointed out to the author by J.-P. Demailly. These facts are actually well known, cf. e.g. [KW70] and [Gau77b], so the proofs can be omitted.

\section{The Bogomolov-Tian-Todorov theorem for Calabi-Yau $\partial\bar\partial$-manifolds}\label{section:BTT}
 
 In this section we prove Theorem \ref{The:BTT} by adapting to the $\partial\bar\partial$ context the proofs given in [Tia87] and [Tod89] for the K\"ahler context. Two preliminary facts are needed. The first is a simple but very useful consequence of the $\partial\bar\partial$-lemma (also observed in [DGMS75]) that will play a key role in this paper. 

\begin{Lem}\label{Lem:d-closed-rep} Let $X$ be a compact $\partial\bar\partial$-manifold. Then any Dolbeault cohomology class of any type $(p, \, q)$ on $X$ can be represented by a $d$-closed form.

\end{Lem}

 In particular, this lemma defines a canonical injection $H^{p,\, q}(X,\,\C)\subset H^{p+q}(X,\,\C)$ (for any $p,q$) by mapping any Dolbeault class $[u]\in H^{p,\, q}(X,\,\C)$ to the De Rham class $\{u\}\in H^{p+q}(X,\,\C)$ of any of its $d$-closed representatives $u$. Fresh applications of the $\partial\bar\partial$-lemma easily show that this map is independent of the choice of the $d$-closed representative of the Dolbeault class and that it is injective.

\vspace{2ex}

\noindent {\it Proof of Lemma \ref{Lem:d-closed-rep}.} Let $\alpha$ be an arbitrary $\bar\partial$-closed $(p, \, q)$-form on $X$. It represents a Dolbeault cohomology class $[\alpha]\in H^{p, \, q}(X, \, \C)$. We will show that there exists a $d$-closed $(p, \, q)$-form $\beta$ on $X$ that is Dolbeault cohomologous to $\alpha$. In other words, we are looking for a $(p, \, q-1)$-form $v$ on $X$ such that

$$\beta:=\alpha + \bar\partial v$$

\noindent is $d$-closed. Since $\alpha$ is $\bar\partial$-closed, the condition $d(\alpha + \bar\partial v)=0$ amounts to  

\begin{eqnarray}\label{eqn:del-prooflemma}\partial\bar\partial v = -\partial\alpha.\end{eqnarray}

\noindent Now, $\partial\alpha$ is both $\partial$-closed (even $\partial$-exact) and $\bar\partial$-closed (since $\alpha$ is $\bar\partial$-closed and $\partial$ anti-commutes with $\bar\partial$), hence $\partial\alpha$ is $d$-closed. Being a $\partial$-exact pure-type form that is $d$-closed, $\partial\alpha$ must be $\partial\bar\partial$-exact by the $\partial\bar\partial$-lemma. Hence a $(p, \, q-1)$-form $v$ satisfying (\ref{eqn:del-prooflemma}) exists. \hfill $\Box$

\vspace{2ex}

 The second preliminary fact is peculiar to manifolds with a trivial canonical bundle. Fix an arbitrary Hermitian metric $\omega$ on a given compact complex manifold $X$ ($n=\mbox{dim}_{\C}X$). Thus $\omega$ is a Hermitian metric on the holomorphic vector bundle $T^{1, \, 0}X$ of vector fields of type $(1, \, 0)$ of $X$. Let $D$ denote the corresponding Chern connection of $(T^{1, \, 0}X, \, \omega)$ and let $D=D' + D''$ be its splitting into components of type $(1, \, 0)$ and respectively $(0, \, 1)$, where $D''=\bar\partial$ is the $\bar\partial$ operator defining the complex structure of $X$. We denote, as usual, by $D^{''\star}=\bar\partial^{\star}$ the formal adjoint of $D''=\bar\partial$ w.r.t. $\omega$ and by 

$$\Delta''=\bar\partial\bar\partial^{\star} + \bar\partial^{\star}\bar\partial: C^{\infty}_{p, \, q}(X, \, T^{1, \, 0}X)\rightarrow C^{\infty}_{p, \, q}(X, \, T^{1, \, 0}X)$$   

\noindent the corresponding anti-holomorphic Laplacian on $T^{1, \, 0}X$-valued $(p, \, q)$-forms. Since $\Delta''$ is elliptic and $\bar\partial^2=0$, the standard three-space orthogonal decomposition holds\!\!:

\begin{eqnarray}\label{eqn:3space-decomp-tangent}C^{\infty}_{0, \, 1}(X, \, T^{1, \, 0}X) = \ker\Delta''\oplus\mbox{Im}\,\bar\partial\oplus\mbox{Im}\,\bar\partial^{\star}.\end{eqnarray}

\noindent On the other hand, a similar orthogonal decomposition holds for the space of scalar-valued $(n-1, \, 1)$-forms on $X$ with respect to $\bar\partial$ acting on scalar forms, its formal adjoint $\bar\partial^{\star}$ and the induced anti-holomorphic Laplacian $\Delta''$. 

 Suppose now that the canonical bundle $K_X$ of $X$ is trivial. Fix a holomorphic $n$-form $u$ with no zeroes on $X$ (i.e. $u$ is a non-vanishing holomorphic section of the trivial line bundle $K_X$). Thus $u$ can be identified with the class $[u]\in H^{n, \, 0}(X,\,\C)\simeq\C$, so it is unique up to a nonzero constant factor. It is then clear that, for every $q=0, \dots , n$, $u$ defines an isomorphism (that may well be called the {\bf Calabi-Yau isomorphism})\!\!:

\begin{eqnarray}\label{eqn:u-isom}T_u\,\,:\,\,C^{\infty}_{0, \, q}(X, \, T^{1, \, 0}X)\stackrel{\cdot\lrcorner u}{\longrightarrow}C^{\infty}_{n-1,\, q}(X, \,\C)\end{eqnarray}

\noindent mapping any $\theta\in C^{\infty}_{0, \, q}(X, \, T^{1, \, 0}X)$ to $T_u(\theta)\!:=\theta\lrcorner u$, where the operation denoted by $\cdot\lrcorner$ combines the contraction of $u$ by the $(1, \, 0)$-vector field component of $\theta$ with the exterior multiplication by the $(0,\,q)$-form component. 

We now record the following well-known fact for future reference.

\begin{Lem}\label{Lem:u-parallel} Let $X$ be a compact complex manifold ($\mbox{dim}_{\C}X=n$) endowed with a Hermitian metric $\omega$ such that $\mbox{Ric}(\omega)=0$. If $K_X$ is trivial and if $u\in C^{\infty}_{n,\,0}(X,\,\C)$ such that $\bar\partial u=0$, $u$ has no zeroes and

\begin{eqnarray}\label{eqn:u-normalisation}i^{n^2}\int\limits_Xu\wedge\bar{u}=\int\limits_XdV_{\omega},\hspace{3ex} (\mbox{where}\,\,\,dV_{\omega}\!:=\frac{\omega^n}{n!})\end{eqnarray}

\noindent then the Calabi-Yau isomorphism $T_u: C^{\infty}_{0,\,1}(X,\,T^{1,\,0}X)\rightarrow C^{\infty}_{n-1,\,1}(X,\,\C)$ (see (\ref{eqn:u-isom}) with $q=1$) is an isometry w.r.t. the {\bf pointwise} (hence also the $L^2$) scalar products induced by $\omega$ on the vector bundles involved.  

\end{Lem}

\noindent {\it Proof.} Fix an arbitrary point $x_0\in X$ and choose local holomorphic coordinates $z_1, \dots , z_n$ about $x_0$ such that 

$$\omega(x_0)=i\,\sum\limits_{j=1}^n\lambda_j\,dz_j\wedge d\bar{z}_j  \hspace{3ex} \mbox{and} \hspace{3ex} u(x_0)= f\,dz_1\wedge\dots\wedge dz_n.$$

\noindent A simple calculation shows that for any $\theta, \eta\in C^{\infty}_{0,\,1}(X,\,T^{1,\,0}X)$, the pointwise scalar products at $x_0$ are related by 

$$\langle\theta,\,\eta\rangle = \frac{\lambda_1\dots\lambda_n}{|f|^2}\,\langle\theta\lrcorner u,\,\eta\lrcorner u\rangle.$$

\noindent Thus having $\langle\theta,\,\eta\rangle = \langle\theta\lrcorner u,\,\eta\lrcorner u\rangle$ at $x_0$ is equivalent to having $|f|^2=\lambda_1\dots\lambda_n$. On the other hand, the identity $i^{n^2}\,u\wedge\bar{u} = |u|_\omega^2\,\omega^n$ implies that 

$$|f|^2=(n!)\,|u|_{\omega}^2\,(\lambda_1\dots\lambda_n).$$ 

\noindent Thus $T_u$ is an isometry w.r.t. the pointwise scalar products induced by $\omega$ if and only if 

\begin{eqnarray}\label{eqn:isometry-cond}|u|_{\omega}^2=\frac{1}{n!}  \hspace{3ex} \mbox{at every point of}\,\,\,X.\end{eqnarray} 

\noindent Since we know from (\ref{eqn:uubar-parallel}) of Proposition \ref{Prop:Bprinc-bal} and from Observation \ref{Obs:flat-parallel} that $|u|_{\omega}^2$ is constant on $X$, we see from the identity $i^{n^2}\,u\wedge\bar{u} = |u|_\omega^2\,\omega^n$ that the normalisation (\ref{eqn:u-normalisation}) of $u$ is equivalent to (\ref{eqn:isometry-cond}), i.e. to $T_u$ being an isometry w.r.t. the pointwise scalar products induced by $\omega$ on the vector bundles involved.  \hfill $\Box$

\vspace{2ex}

We shall now compare in the case $q=1$ the image under the operation $\cdot\lrcorner u$ of the three-space decomposition (\ref{eqn:3space-decomp-tangent}) of $C^{\infty}_{0, \, 1}(X, \, T^{1, \, 0}X)$ with the analogous three-space decomposition of $C^{\infty}_{n-1,\, 1}(X, \,\C)$.

\begin{Lem}\label{Lem:3space-decomp-contr} Let $X$ be a compact complex manifold ($n=\mbox{dim}_{\C}X$) such that $K_X$ is trivial. Then, for $q=1$, the isomorphism $T_u$ of (\ref{eqn:u-isom}) satisfies\!\!:

\begin{eqnarray}\label{eqn:Tuker}T_u(\ker\bar\partial)=\ker\bar\partial \hspace{3ex}\mbox{and}\hspace{3ex} T_u(\mbox{Im}\,\bar\partial)=\mbox{Im}\,\bar\partial.\end{eqnarray}

\noindent Hence $T_u$ induces an isomorphism in cohomology

\begin{eqnarray}\label{eqn:u-isom-classes}T_{[u]}\,\,:\,\,H^{0, \, 1}(X, \, T^{1, \, 0}X)\stackrel{\cdot\lrcorner [u]}{\longrightarrow}H^{n-1,\, 1}(X, \,\C)\end{eqnarray}

\noindent defined by $T_{[u]}([\theta])=[\theta\lrcorner u]$ for all $[\theta]\in H^{0, \, 1}(X, \, T^{1, \, 0}X)$.

\noindent If $\omega$ is any Hermitian metric on $X$ such that $\mbox{Ric}(\omega)=0$, $T_u$ also satisfies\!\!:

\begin{eqnarray}\label{eqn:Tu-dbarstar}T_u(\mbox{Im}\,\bar\partial^{\star})=\mbox{Im}\,\bar\partial^{\star}\hspace{3ex}\mbox{and}\hspace{3ex}T_u(\ker\Delta'')=\ker\Delta''.\end{eqnarray}

\end{Lem}

\noindent {\it Proof.} It relies on the easily checked formulae\!\!:

\begin{eqnarray}\label{eqn:Leibniz-dbar-u}\bar\partial(\theta\lrcorner u)=(\bar\partial\theta)\lrcorner u + \theta\lrcorner(\bar\partial u)=(\bar\partial\theta)\lrcorner u, \hspace{1ex} \bar\partial(\xi\lrcorner u)=(\bar\partial\xi)\lrcorner u - \xi\lrcorner(\bar\partial u)=(\bar\partial\xi)\lrcorner u\end{eqnarray}

\noindent for all $\theta\in C^{\infty}_{0, \, 1}(X, \, T^{1, \, 0}X)$ and all $\xi\in C^{\infty}(X, \, T^{1, \, 0}X)$. Note, however, that the analogous identities for $\partial$ fail. These formulae imply the inclusions\!\!:

$$T_u(\ker\bar\partial)\subset\ker\bar\partial \hspace{3ex} \mbox{and}  \hspace{3ex} T_u(\mbox{Im}\,\bar\partial)\subset\mbox{Im}\,\bar\partial.$$

 To prove the reverse inclusion of the former equality in (\ref{eqn:Tuker}), suppose that $\theta\lrcorner u\in\ker\bar\partial$ for some $\theta\in C^{\infty}_{0, \, 1}(X, \, T^{1, \, 0}X)$. By (\ref{eqn:Leibniz-dbar-u}), this means that $(\bar\partial\theta)\lrcorner u=0$, which is equivalent to $\bar\partial\theta=0$ since the map $T_u$ of (\ref{eqn:u-isom}) is an isomorphism. (Here $q=2$.) 

 To prove the reverse inclusion of the latter equality in (\ref{eqn:Tuker}), let $\theta\in C^{\infty}_{0, \, 1}(X, \, T^{1, \, 0}X)$ such that $\theta\lrcorner u=\bar\partial v$ for some $(n-1, \, 0)$-form $v$. With respect to local holomorphic coordinates $z_1, \dots , z_n$ on some open subset $U\subset X$, let

$$\theta=\sum\limits_{j, \, k}\theta^j_k\,d\bar{z}_k\otimes\frac{\partial}{\partial z_j} \hspace{3ex} \mbox{and}  \hspace{3ex} u=f\,dz_1\wedge\dots\wedge dz_n,$$

\noindent where $f$ is a holomorphic function with no zeroes on $U$. Then

$$\theta\lrcorner u = \sum\limits_{j, \, k}(-1)^{j-1}f\theta^j_k\,d\bar{z}_k\wedge dz_1\wedge\dots\wedge\widehat{dz_j}\wedge\dots\wedge dz_n.$$

\noindent Letting $v=\sum\limits_j v_j\, dz_1\wedge\dots\wedge\widehat{dz_j}\wedge\dots\wedge dz_n$, the condition $\theta\lrcorner u=\bar\partial v$ reads

$$\sum\limits_{j, \, k}(-1)^{j-1}f\theta^j_k\,d\bar{z}_k\wedge dz_1\wedge\dots\wedge\widehat{dz_j}\wedge\dots\wedge dz_n = \sum\limits_{j, \, k}\frac{\partial v_j}{\partial\bar{z}_k}\,d\bar{z}_k\wedge dz_1\wedge\dots\wedge\widehat{dz_j}\wedge\dots\wedge dz_n,$$

\noindent which is equivalent to $\theta^j_k=\frac{\partial}{\partial\bar{z}_k}((-1)^{j-1}\,\frac{v_j}{f})$ for all $j, k$ since $f$ is holomorphic without zeroes. Setting $\xi_j:=(-1)^{j-1}\,\frac{v_j}{f}$ for all $j$, we get the local representative of a global vector field

$$\xi:=\sum\limits_j\xi_j\,\frac{\partial}{\partial z_j}\in C^{\infty}(X, \, T^{1, \, 0}X)$$

\noindent satisfying $\theta = \bar\partial\xi$ on $X$. Hence $\theta\in\mbox{Im}\,\bar\partial.$ We have thus proved that $\mbox{Im}\,\bar\partial\subset T_u(\mbox{Im}\,\bar\partial)$, hence the latter identity in (\ref{eqn:Tuker}).

 Thus the identities (\ref{eqn:Tuker}) are proved. Then so is (\ref{eqn:u-isom-classes}), an obvious consequence of (\ref{eqn:Tuker}).

 To get (\ref{eqn:Tu-dbarstar}), recall that we know from (\ref{eqn:uubar-parallel}) of Proposition \ref{Prop:Bprinc-bal} and from Observation \ref{Obs:flat-parallel} that $|u|_{\omega}^2$ is constant on $X$ whenever $\mbox{Ric}(\omega)=0$. Then the proof of Lemma \ref{Lem:u-parallel} shows that $\langle\theta,\,\eta\rangle=\mbox{Const}\cdot\langle\theta\lrcorner u,\,\eta\lrcorner u\rangle$ for all $\theta, \eta\in C^{\infty}_{0,\,1}(X,\,T^{1,\,0}X)$, hence $\theta\perp\eta$ if and only if $\theta\lrcorner u\perp\eta\lrcorner u$. (The notation is the obvious one.) This fact suffices to deduce (\ref{eqn:Tu-dbarstar}) from the pairwise orthogonality of $\ker\Delta''$, $\mbox{Im}\,\bar\partial$ and $\mbox{Im}\,\bar\partial^{\star}$ in the three-space decompositions of $C^{\infty}_{0,\,1}(X,\,T^{1,\,0}X)$ and $C^{\infty}_{n-1,\,1}(X,\,\C)$ and from the identities (\ref{eqn:Tuker}). The proof is complete.   \hfill $\Box$

\vspace{2ex}

Now recall that the isomorphisms $T_u$ of (\ref{eqn:u-isom}) and their inverses allow one to define a Lie bracket on $\oplus_qC^{\infty}_{n-1, \, q}(X, \, \C)$ by setting (cf. [Tia87, p. 631])\!\!:

\begin{eqnarray}\label{eqn:Liebracket}[\zeta_1,\, \zeta_2]:=T_u\bigg[T_u^{-1}\zeta_1, \, T_u^{-1}\zeta_2\bigg]\in C^{\infty}_{n-1, \, q_1+q_2}(X, \, T^{1, \, 0}X)\end{eqnarray}

\noindent for any forms $\zeta_1\in C^{\infty}_{n-1, \, q_1}(X, \, \C)$ and $\zeta_2\in C^{\infty}_{n-1, \, q_2}(X, \, \C)$, where the operation $[\,\,,\,\,]$ on the right-hand side of (\ref{eqn:Liebracket}) combines the Lie bracket of the $T^{1, \, 0}X$-parts of $T_u^{-1}\zeta_1\in C^{\infty}_{0, \, q_1}(X, \, T^{1, \, 0}X)$ and $T_u^{-1}\zeta_2\in C^{\infty}_{0, \, q_2}(X, \, T^{1, \, 0}X)$ with the wedge product of their $(0, \, q_1)-$ and respectively $(0, \, q_2)$-form parts. The definition (\ref{eqn:Liebracket}) can be reformulated as\!\!:

\begin{eqnarray}\label{eqn:Lie-bracket-comm}[\Phi_1,\,\Phi_2]\lrcorner u = [\Phi_1\lrcorner u,\,\Phi_2\lrcorner u] \hspace{2ex} \mbox{for all}\hspace{1ex}\Phi_1,\Phi_2\in C^{\infty}_{0,\,1}(X,\,T^{1,\,0}X).\end{eqnarray}

The main technical ingredient in the proofs of [Tia87] and [Tod89] was the following general observation, the so-called {\it Tian-Todorov lemma}.

\begin{Lem}\label{Lem:basic-trick}(cf. Lemma 3.1. in [Tia87], Lemma 1.2.4. in [Tod89]) Let $X$ be a compact complex manifold ($n=\mbox{dim}_{\C}X$) such that $K_X$ is trivial. Then, for every forms $\zeta_1, \zeta_2\in C^{\infty}_{n-1, \, 1}(X, \, \C)$ such that $\partial\zeta_1=\partial\zeta_2=0$, we have

$$[\zeta_1, \, \zeta_2]\in\mbox{Im}\,\partial.$$

\noindent More precisely, the identity $[\theta_1\lrcorner u,\, \theta_2\lrcorner u] = \partial(\theta_1\lrcorner(\theta_2\lrcorner u))$ holds for $\theta_1, \theta_2\in C^{\infty}_{0,\, 1}(X,\, T^{1,\, 0}X)$ whenever $\partial(\theta_1\lrcorner u) = \partial(\theta_2\lrcorner u)=0$.

\end{Lem}

 We can now briefly review the main arguments in the proofs of [Tia87] and [Tod89] by pointing out that they are still valid when the K\"ahler assumption is weakened to the $\partial\bar\partial$ assumption.

\vspace{2ex}

\noindent {\it Proof of Theorem \ref{The:BTT}.} Let $[\eta]\in H^{0, \, 1}(X, \, T^{1, \, 0}X)$ be an arbitrary nonzero class. Pick any $d$-closed representative $w_1$ of the class $[\eta]\lrcorner [u]\in H^{n-1, \, 1}(X, \, \C)$. Such a $d$-closed representative exists by Lemma \ref{Lem:d-closed-rep} thanks to the $\partial\bar\partial$ assumption on $X$. This is virtually the only modification of the proof compared to the K\"ahler case where the $\Delta''$-harmonic representative of the class $[\eta]\lrcorner [u]$ was chosen. Since $\Delta'=\Delta''$ in the K\"ahler case, $\Delta''$-harmonic forms are also $\partial$-closed, hence $d$-closed, but this no longer holds in the non-K\"ahler case. 

 Since $T_u$ is an isomorphism, there is a unique $\Phi_1\in C^{\infty}_{0, \, 1}(X, \, T^{1, \, 0}X)$ such that $\Phi_1\lrcorner u=w_1$. Now $\bar\partial w_1=0$, so the former equality in (\ref{eqn:Tuker}) implies that $\bar\partial \Phi_1=0$. Moreover, since $[\Phi_1\lrcorner u]=[w_1]$, (\ref{eqn:u-isom-classes}) implies that $[\Phi_1]=[\eta]\in H^{0, \, 1}(X, \, T^{1, \, 0}X)$ and this is the original class we started off with. However, $\Phi_1$ need not be the $\Delta''$-harmonic representative of the class $[\eta]$ in the non-Kaehler case (in contrast to the K\"ahler case of [Tia87] and [Tod89]). Meanwhile, by the choice of $w_1$, we have 

$$\partial(\Phi_1\lrcorner u)=0,$$

\noindent so Lemma \ref{Lem:basic-trick} applied to $\zeta_1=\zeta_2=\Phi_1\lrcorner u$ yields $[\Phi_1\lrcorner u, \, \Phi_1\lrcorner u]\in\mbox{Im}\,\partial$. On the other hand, $[\Phi_1\lrcorner u, \, \Phi_1\lrcorner u]\in\ker\bar\partial$ as can be easily checked and is well-known (see e.g. [Tod89, Lemma 1.2.5]). By the $\partial\bar\partial$-lemma applied to the $(n-1, \, 2)$-form $1/2\,[\Phi_1\lrcorner u, \, \Phi_1\lrcorner u]$, there exists $\psi_2\in C^{\infty}_{n-2, \, 1}(X,\, \C)$ such that

$$\bar\partial\partial\psi_2=\frac{1}{2}\,[\Phi_1\lrcorner u, \, \Phi_1\lrcorner u].$$

\noindent We can choose $\psi_2$ of minimal $L^2$-norm with this property (i.e. $\psi_2\in\mbox{Im}(\partial\bar\partial)^{\star}$, see e.g. the explanation after Definition \ref{def:min-d-closed} in terms of the Aeppli cohomology). Put $w_2:=\partial\psi_2\in C^{\infty}_{n-1, \, 1}(X, \, \C)$. Since $T_u$ is an isomorphism, there is a unique $\Phi_2\in C^{\infty}_{0, \, 1}(X, \, T^{1, \, 0}X)$ such that $\Phi_2\lrcorner u=w_2$. Implicitly, $\partial(\Phi_2\lrcorner u)=0$. Moreover, using (\ref{eqn:Leibniz-dbar-u}), we get

$$(\bar\partial\Phi_2)\lrcorner u=\bar\partial(\Phi_2\lrcorner u)=\frac{1}{2}\,[\Phi_1\lrcorner u, \, \Phi_1\lrcorner u] = \frac{1}{2}\, [\Phi_1, \, \Phi_1]\lrcorner u,$$

\noindent where the last identity follows from (\ref{eqn:Lie-bracket-comm}). Hence

$$(\mbox{Eq.}\,1) \hspace{6ex} \bar\partial\Phi_2=\frac{1}{2}\, [\Phi_1, \, \Phi_1].$$

 We can now continue inductively. Suppose we have constructed $\Phi_1, \dots , \Phi_{N-1}\in C^{\infty}_{0, \, 1}(X, \, T^{1, \, 0}X)$ such that

$$\partial(\Phi_k\lrcorner u)=0 \hspace{3ex} \mbox{and} \hspace{3ex} \bar\partial(\Phi_k\lrcorner u)=\frac{1}{2}\,\sum\limits_{l=1}^{k-1}[\Phi_l\lrcorner u, \, \Phi_{k-l}\lrcorner u], \hspace{3ex} 1\leq k\leq N-1.$$

\noindent By formulae (\ref{eqn:Leibniz-dbar-u}), (\ref{eqn:Lie-bracket-comm}) and since $T_u$ is an isomorphism, the latter identity above is equivalent to

$$(\mbox{Eq.}\,(k-1)) \hspace{6ex} \bar\partial\Phi_k=\frac{1}{2}\,\sum\limits_{l=1}^{k-1}[\Phi_l, \, \Phi_{k-l}], \hspace{3ex} 1\leq k\leq N-1.$$

\noindent Then it is easily seen and well known (cf. [Tod89, Lemma 1.2.5]) that 

$$\frac{1}{2}\,\sum\limits_{l=1}^{N-1}[\Phi_l\lrcorner u, \, \Phi_{N-l}\lrcorner u]\in\ker\bar\partial.$$

\noindent On the other hand, since $\Phi_1\lrcorner u, \dots , \Phi_{N-1}\lrcorner u\in\ker\partial$, Lemma \ref{Lem:basic-trick} gives

$$[\Phi_l\lrcorner u, \, \Phi_{N-l}\lrcorner u]\in\mbox{Im}\, \partial \hspace{3ex} \mbox{for all}\hspace{3ex} l=1, \dots, N-1.$$

\noindent Thanks to the last two relations, the $\partial\bar\partial$-lemma implies the existence of a form $\psi_N\in C^{\infty}_{n-2, \, 1}(X, \, \C)$ such that

$$\bar\partial\partial\psi_N=\frac{1}{2}\,\sum\limits_{l=1}^{N-1}[\Phi_l\lrcorner u, \, \Phi_{N-l}\lrcorner u].$$

\noindent We can choose $\psi_N$ of minimal $L^2$-norm with this property (i.e. $\psi_N\in\mbox{Im}(\partial\bar\partial)^{\star}$). Letting $w_N:=\partial\psi_N\in C^{\infty}_{n-1, \, 1}$, there exists a unique $\Phi_N\in C^{\infty}_{0, \, 1}(X, \, T^{1, \, 0}X)$ such that $\Phi_N\lrcorner u=w_N$. Implicitly

$$\partial(\Phi_N\lrcorner u)=0.$$

\noindent We also have $\bar\partial(\Phi_N\lrcorner u)=\frac{1}{2}\,\sum\limits_{l=1}^{N-1}[\Phi_l\lrcorner u, \, \Phi_{N-l}\lrcorner u]$ by construction. By formulae (\ref{eqn:Leibniz-dbar-u}), (\ref{eqn:Lie-bracket-comm}) and since $T_u$ is an isomorphism, this amounts to

$$(\mbox{Eq.}\,(N-1))\hspace{6ex} \bar\partial\Phi_N=\frac{1}{2}\,\sum\limits_{l=1}^{N-1}[\Phi_l, \, \Phi_{N-l}].$$

\noindent We have thus shown inductively that the equation $(\mbox{Eq.}\,k)$ is
solvable for every $k\in\N^{\star}$. It is well-known (cf. [Kur62]) that in this case the series $\Phi(t):=\Phi_1\,t + \Phi_2\,t^2 + \dots + \Phi_N\,t^N + \dots$ converges in a H\"older norm for all $t\in\C$ such that $|t|<\varepsilon\ll 1$. This produces a form $\Phi(t)\in C^{\infty}_{0,\,1}(X,\,T^{1,\,0}X)$ which defines a complex structure $\bar\partial_t$ on $X$ that identifies with $\bar\partial-\Phi(t)$ and is the deformation of the original complex structure $\bar\partial$ of $X$ in the direction of the originally given $[\eta]\in H^{0,\,1}(X,\, T^{1,\,0}X)$. The proof is complete.  \hfill $\Box$

\hspace{3ex}

 We end this section by noticing that the full force of the $\partial\bar\partial$ assumption is not needed in Theorem \ref{The:BTT}, but only a special case thereof, since only two applications in very particular situations have been made of it.

 First, we needed any Dolbeault cohomology class $[\alpha]\in H^{n-1,\,1}(X,\,\C)$ (denoted by $[\eta]\lrcorner [u]$ in the proof) to be representable by a $d$-closed form. The proof of Lemma \ref{Lem:d-closed-rep} shows this to be equivalent to requiring that any $\partial$-exact $(n,\,1)$-form $\partial\alpha$ for which $\bar\partial\alpha = 0$ be $\partial\bar\partial$-exact. This is equivalent to requiring the following linear map (which is always well defined)

\begin{equation}\label{eqn:A_1-def}A_1\,\,:\,\,H^{n-1,\,1}_{\bar\partial}(X,\,\C)\longrightarrow H^{n,\,1}_{BC}(X,\,\C), \hspace{3ex} [\alpha]_{\bar\partial}\mapsto [\partial\alpha]_{BC}\end{equation}

\noindent to vanish identically, where the subscript $BC$ indicates a Bott-Chern cohomology group. By duality, the vanishing of $A_1$ is equivalent to the vanishing of its dual map

\begin{equation}\label{eqn:A_1star-def}A_1^{\star}\,\,:\,\,H^{0,\,n-1}_A(X,\,\C)\longrightarrow H^{1,\,n-1}_{\bar\partial}(X,\,\C), \hspace{3ex} [u]_A\mapsto [\partial u]_{\bar\partial},\end{equation}

\noindent where the subscript $A$ indicates an Aeppli cohomology group.

 The other special case of the $\partial\bar\partial$ lemma needed in the proof of Theorem \ref{The:BTT} was the requirement that any $\partial$-exact and $d$-closed $(n-1,\,2)$-form $\beta$ (denoted by $[\Phi_1\lrcorner u, \, \Phi_1\lrcorner u]$ in the proof) be $\partial\bar\partial$-exact. This is equivalent to requiring the following linear map (which is always well defined) 

\begin{equation}\label{eqn:B-def}B\,\,:\,\,H^{n-1,\,2}_{BC}(X,\,\C)\longrightarrow H^{n-1,\,2}_{\partial}(X,\,\C), \hspace{3ex} [\beta]_{BC}\mapsto [\beta]_{\partial}\end{equation}

\noindent to be injective. From the exact sequence

$$H^{n-2,\,2}_A(X,\,\C)\stackrel{A_2}{\longrightarrow} H^{n-1,\,2}_{BC}(X,\,\C)\stackrel{B}{\longrightarrow} H^{n-1,\,2}_{\partial}(X,\,\C),$$

\noindent we infer that $B$ being injective is equivalent to the linear map $A_2$ vanishing identically, where

\begin{equation}\label{eqn:A_2-def}A_2\,\,:\,\,H^{n-2,\,2}_A(X,\,\C)\longrightarrow H^{n-1,\,2}_{BC}(X,\,\C), \hspace{3ex} [v]_A\mapsto [\partial v]_{BC}.\end{equation}

 This discussion can be summed up as follows.

\begin{Obs}\label{Obs:BTTgeneral} Let $X$ be a compact complex manifold with $\mbox{dim}_{\C}X=n$ whose canonical bundle $K_X$ is trivial such that the linear maps $A_1$ and $A_2$ defined in (\ref{eqn:A_1-def}) and (\ref{eqn:A_2-def}) vanish identically. Then the Kuranishi family of $X$ is unobstructed.

\end{Obs}

\section{Co-polarised deformations of balanced Calabi-Yau $\partial\bar\partial$-manifolds}\label{section:co-polarisations}

\subsection{Definitions}\label{subsection:definitions}

 Let $(X, \, \omega)$ be a compact {\it balanced} Calabi-Yau $\partial\bar\partial$- manifold ($n=\mbox{dim}_{\C}X$). Denote by $\pi:{\cal X}\rightarrow\Delta$ the Kuranishi family of $X$. Thus $\pi$ is a proper holomorphic submersion from a complex manifold ${\cal X}$, while the fibres $X_t$ with $t\in\Delta\setminus\{0\}$ can be seen as deformations of the given manifold $X_0=X$. The base space $\Delta$ is smooth and can be viewed as an open subset of $H^{0, \, 1}(X, \, T^{1, \, 0}X)$ (or as a ball containing the origin in $\C^N$, where $N=\mbox{dim}_{\C}H^{0, \, 1}(X, \, T^{1, \, 0}X)$) by Theorem \ref{The:BTT}. Hence the tangent space at $0$ is

$$T_0\Delta\simeq H^{0, \, 1}(X, \, T^{1, \, 0}X).$$

\noindent By Wu's result [Wu06, Theorem 5.13, p. 56], small deformations of balanced $\partial\bar\partial$-manifolds are again balanced $\partial\bar\partial$-manifolds. Hence, in our case, $X_t$ is a balanced Calabi-Yau $\partial\bar\partial$-manifold for all $t\in\Delta$ sufficiently close to $0$.

 Recall that in the special case where the class $[\omega]\in H^{1, \, 1}(X, \, \C)\subset H^2(X, \, \C)$ is K\"ahler (and is furthermore often required to be integral, but we deal with arbitrary, possibly non-rational classes here), it is standard to define the deformations of $X_0=X$ {\it polarised} by $[\omega]$ as those nearby fibres $X_t$ on which the De Rham class $\{\omega\}\in H^2(X, \, \C)$ is still a K\"ahler class (hence, in particular, of type $(1, \, 1)$) for the complex structure $J_t$ of $X_t$. In the more general balanced case treated here, $\omega$ need not define a class, but $\omega^{n-1}$ does. Taking our cue from the standard K\"ahler case, we propose the following dual notion in the balanced context.

\begin{Def}\label{Def:co-polarisation} Having fixed a balanced class 

$$[\omega^{n-1}]\in H^{n-1, \, n-1}(X, \, \C)\subset H^{2n-2}(X, \, \C),$$ 

\noindent we say that a fibre $X_t$ is {\bf co-polarised} by $[\omega^{n-1}]$ if the De Rham class 

$$\{\omega^{n-1}\}\in H^{2n-2}(X, \, \C)$$ 

\noindent is of type $(n-1, \, n-1)$ for the complex structure $J_t$ of $X_t$.

 The restricted family $\pi:{\cal X}_{[\omega^{n-1}]}\rightarrow\Delta_{[\omega^{n-1}]}$ will be called the universal family of deformations of $X$ that are co-polarised by the balanced class $[\omega^{n-1}]$, where $\Delta_{[\omega^{n-1}]}$ is the set of $t\in\Delta$ such that $X_t$ is co-polarised by $[\omega^{n-1}]$ and $
{\cal X}_{[\omega^{n-1}]}=\pi^{-1}(\Delta_{[\omega^{n-1}]})\subset{\cal X}$.

\end{Def}

 After possibly shrinking $\Delta_{[\omega^{n-1}]}$ about $0$, we can assume that $[\omega^{n-1}]\in H^{n-1, \, n-1}(X_t, \, \C)$ is a balanced class for the complex structure $J_t$ of the fibre $X_t$ for every $t\in\Delta_{[\omega^{n-1}]}$ (cf. Observation \ref{Obs:balanced-def-class_appendix}). 

 Note that in the special case where $\omega$ is K\"ahler on $X_0=X$, the $(2n-2)$-class $\{\omega^{n-1}\}$ is a balanced class for $J_t$ whenever the $2$-class $\{\omega\}$ is a K\"ahler class for $J_t$. We shall see further down that the converse also holds, meaning that in the special K\"ahler case the notion of co-polarised deformations of $X$ coincides with that of polarised deformations. Recall that when $\omega$ is K\"ahler, the deformations of $X$ polarised by $[\omega]$ are parametrised by the following subspace of $H^{0, \, 1}(X, \, T^{1, \, 0}X)$\!\!:

\begin{eqnarray}\label{eqn:pol-space}H^{0, \, 1}(X, \, T^{1, \, 0}X)_{[\omega]}:=\bigg\{[\theta]\in H^{0, \, 1}(X, \, T^{1, \, 0}X)\,\,;\,\, [\theta\lrcorner\omega]=0\in H^{0, \, 2}(X, \, \C)\bigg\}\end{eqnarray}

\noindent which is isomorphic under the restriction of $T_{[u]}$ (cf. (\ref{eqn:u-isom-classes})) to the space of {\it primitive} Dolbeault classes of type $(n-1, \, 1)$\!\!:

\begin{eqnarray}\label{eqn:pol-space-isom}H^{0, \, 1}(X, \, T^{1, \, 0}X)_{[\omega]}\stackrel{T_{[u]}}{\longrightarrow}H^{n-1, \, 1}_{prim}(X, \, \C).\end{eqnarray}

 We shall now see that the co-polarised deformations of $X$ are parametrised by an analogous subspace.

\begin{Lem}\label{Lem:co-pol-base} For a given balanced class $[\omega^{n-1}]\in H^{n-1, \, n-1}(X, \, \C)$, consider the following vector subspace of $H^{0, \, 1}(X, \, T^{1, \, 0}X)$\!\!:

\begin{eqnarray}\label{eqn:co-pol-space} H^{0, \, 1}(X, \, T^{1, \, 0}X)_{[\omega^{n-1}]}:=\bigg\{[\theta]\in H^{0, \, 1}(X, \, T^{1, \, 0}X)\,\,;\,\, [\theta\lrcorner\omega^{n-1}]=0\in H^{n-2, \, n}(X, \, \C)\bigg\}.\end{eqnarray}

 \noindent Then\!\!: 

 \noindent $(a)$\, the space $H^{0, \, 1}(X, \, T^{1, \, 0}X)_{[\omega^{n-1}]}$ is well-defined (i.e. the class $[\theta\lrcorner\omega^{n-1}]\in H^{n-2, \, n}(X, \, \C)$ is independent of the choice of representative $\theta$ in the class $[\theta]\in H^{0, \, 1}(X, \, T^{1, \, 0}X)$ and of the choice of representative $\omega^{n-1}$ in the class $[\omega^{n-1}]\in H^{n-1, \, n-1}(X, \, \C)$). We can therefore denote

\begin{eqnarray}\label{eqn:not-contr-classes}[\theta]\lrcorner [\omega^{n-1}]:=[\theta\lrcorner\omega^{n-1}].\end{eqnarray}

\noindent $(b)$\, the open subset $\Delta\subset H^{0, \, 1}(X, \, T^{1, \, 0}X)$ satisfies

$$\Delta_{[\omega^{n-1}]}=\Delta\cap H^{0, \, 1}(X, \, T^{1, \, 0}X)_{[\omega^{n-1}]}.$$

\noindent Implicitly, $T_0\Delta_{[\omega^{n-1}]}\simeq H^{0, \, 1}(X, \, T^{1, \, 0}X)_{[\omega^{n-1}]}.$

\end{Lem}

\noindent {\it Proof.} $(a)$ follows from Lemma \ref{Lem:dbar-contr-omega} below. Indeed, if $\theta + \bar\partial\xi$ is another representative of the class $[\theta]$ for some vector field $\xi\in C^{\infty}(X,\,T^{1, \, 0}X)$, then 

$$(\theta + \bar\partial\xi)\lrcorner\omega^{n-1} = \theta\lrcorner\omega^{n-1} + \bar\partial(\xi\lrcorner\omega^{n-1})$$

\noindent since $\omega$ is balanced. Hence $[(\theta + \bar\partial\xi)\lrcorner\omega^{n-1}] = [\theta\lrcorner\omega^{n-1}]$. Similarly, if $\omega^{n-1}+\bar\partial\lambda$ is another representative of the Dolbeault class $[\omega^{n-1}]$ for some $(n-1,\,n-2)$-form $\lambda$, then 

$$\theta\lrcorner(\omega^{n-1} + \bar\partial\lambda) = \theta\lrcorner\omega^{n-1} + \bar\partial(\theta\lrcorner\lambda)$$

\noindent since $\bar\partial\theta=0$. Hence $[\theta\lrcorner(\omega^{n-1} + \bar\partial\lambda)]=[\theta\lrcorner\omega^{n-1}]$.

$(b)$\, Since $X_t$ is a $\partial\bar\partial$-manifold for every $t$ close to $0$, it admits a Hodge decomposition which in degree $2n-2$ spells

$$H^{2n-2}(X, \, \C)=H^{n, \, n-2}(X_t, \, \C)\oplus H^{n-1, \, n-1}(X_t, \, \C)\oplus H^{n-2, \, n}(X_t, \, \C),$$ 

\noindent with $H^{n-2, \, n}(X_t, \, \C)\simeq\overline{H^{n, \, n-2}(X_t, \, \C)}$. In our case, the real De Rham class $\{\omega^{n-1}\}\in H^{2n-2}(X, \, \R)$ splits accordingly as

$$\{\omega^{n-1}\} = \{\omega^{n-1}\}_t^{n, \, n-2} + \{\omega^{n-1}\}_t^{n-1, \, n-1} + \{\omega^{n-1}\}_t^{n-2, \, n},$$

\noindent with $\{\omega^{n-1}\}_t^{n-2, \, n}=\overline{\{\omega^{n-1}\}_t^{n, \, n-2}}$ and $\{\omega^{n-1}\}_t^{n-1, \, n-1}$ real. Thus the definition of $\Delta_{[\omega^{n-1}]}$ translates to

$$\Delta_{[\omega^{n-1}]}=\bigg\{t\in\Delta\,\,;\,\,\{\omega^{n-1}\}_t^{n-2, \, n}=0\in H^{n-2, \, n}(X_t, \, \C)\bigg\}.$$

\noindent Moreover, $\{\omega^{n-1}\}$ is of type $(n-1, \, n-1)$ for $J_0$, so $\{\omega^{n-1}\}_0^{n-2, \, n}=0$ and $\{\omega^{n-1}\}_0^{n, \, n-2}=0$. Let $t_1, \dots , t_N$ be local holomorphic coordinates about $0$ in $\Delta$. So $t=(t_1, \dots , t_N)\in\Delta$ identifies with $[\theta]$ varying in an open subset of $H^{0, \, 1}(X, \, T^{1, \, 0}X)$. Let $[\theta]\in H^{0, \, 1}(X, \, T^{1, \, 0}X)$ be the image of $\frac{\partial}{\partial t_i}_{|t_i=0}$ under the Kodaira-Spencer map $\rho : T_0\Delta\stackrel{\simeq}{\longrightarrow} H^{0, \, 1}(X, \, T^{1, \, 0}X)$. Then, under the Gauss-Manin connection on the Hodge bundle $\Delta\ni t\mapsto H^{2n-2}(X_t, \, \C)$, the derivative of the class $[\omega^{n-1}]_t^{n-2, \, n}\in H^{n-2, \, n}(X_t, \, \C)$ in the direction of $t_i$ at $t_i=0$ is the class $[\theta\lrcorner\omega^{n-1}]\in H^{n-2, \, n}(X, \, \C)$. \hfill $\Box$

\vspace{2ex} Here is the lemma that has been used in the proof of $(a)$ above.

\begin{Lem}\label{Lem:dbar-contr-omega} Let $X$ be a compact complex manifold ($\mbox{dim}_{\C}X=n$) equipped with an arbitrary Hermitian metric $\omega$. Then: \\

\noindent $(i)$\, $\bar\partial(\xi\lrcorner\omega^{n-1})=(\bar\partial\xi)\lrcorner\omega^{n-1} - \xi\lrcorner\bar\partial\omega^{n-1},  \hspace{3ex}\mbox{for every}\hspace{1ex} \xi\in C^{\infty}(X,\,T^{1, \, 0}X).$ \\

\noindent Hence, if $\omega$ is balanced, we have $\bar\partial(\xi\lrcorner\omega^{n-1})=(\bar\partial\xi)\lrcorner\omega^{n-1}$.\\

\noindent $(ii)$\, $\bar\partial(\theta\lrcorner\omega)=(\bar\partial\theta)\lrcorner\omega + \theta\lrcorner\bar\partial\omega,  \hspace{3ex}\mbox{for every}\hspace{1ex} \theta\in C^{\infty}_{0, \, 1}(X,\,T^{1, \, 0}X).$ \\

 Analogous identities hold for forms of any type in place of $\omega$ or $\omega^{n-1}$. However, the analogous identities for $\partial$ in place of $\bar\partial$ fail (intuitively because $\partial$ increases the holomorphic degree of forms, while the contraction by a vector field of type $(1, \, 0)$ decreases the same holomorphic degree). 

\end{Lem}

\noindent {\it Proof.} See Appendix (section \ref{section:appendix}).  \hfill $\Box$

\subsection{Comparison to polarisations of the K\"ahler case}\label{subsection:comparisonK}

 We now pause to observe that in the special case of a K\"ahler class $[\omega]\in H^{1, \, 1}(X, \, \C)$, co-polarised deformations of $X$ coincide with polarised deformations. Thus, although the space $H^{0, \, 1}(X, \, T^{1, \, 0}X)_{[\omega]}$ of (\ref{eqn:pol-space}) no longer makes sense for a non-K\"ahler $\omega$, $H^{0, \, 1}(X, \, T^{1, \, 0}X)_{[\omega^{n-1}]}$ defined in (\ref{eqn:co-pol-space}) naturally extends its meaning to the case of a balanced class $[\omega^{n-1}]$.

\begin{Prop}\label{Prop:pol-copol-equality} Let $(X, \, \omega)$ be a compact K\"ahler manifold ($n=\mbox{dim}_{\C}X$) such that $K_X$ is trivial. Then the following identity holds\!\!:

\begin{eqnarray}\label{eqn:pol-copol-equality}H^{0, \, 1}(X, \, T^{1, \, 0}X)_{[\omega]}= H^{0, \, 1}(X, \, T^{1, \, 0}X)_{[\omega^{n-1}]}.\end{eqnarray}

\end{Prop}

\noindent {\it Proof.} We start by noticing that for any Hermitian metric $\omega$ (no assumption is necessary on $\omega$ here) and any $\theta\in C^{\infty}_{0, \, 1}(X, \, T^{1, \, 0}X)$, we have

\begin{eqnarray}\label{eqn:contr-Leinniz}\theta\lrcorner\omega^k=k\,\omega^{k-1}\wedge(\theta\lrcorner\omega)  \hspace{3ex}\mbox{for any}\hspace{2ex} k.\end{eqnarray}

\noindent This follows from the property $\theta\lrcorner(\omega\wedge\omega^{k-1})=(\theta\lrcorner\omega)\wedge\omega^{k-1} + \omega\wedge(\theta\lrcorner\omega^{k-1})$.

 Suppose now that $\omega$ is K\"ahler and let $[\theta]\in H^{0, \, 1}(X, \, T^{1, \, 0}X)_{[\omega]}$, i.e. $\theta\lrcorner\omega$ is $\bar\partial$-exact. Writing $\theta\lrcorner\omega=\bar\partial v$ for some $(0, \, 1)$-form $v$, from (\ref{eqn:contr-Leinniz}) we get\!\!:

$$\theta\lrcorner\omega^{n-1}=(n-1)\,\omega^{n-2}\wedge\bar\partial v=(n-1)\,\bar\partial(\omega^{n-2}\wedge v)$$

\noindent since $\bar\partial\omega^{n-2}=0$ by the K\"ahler assumption on $\omega$. Thus $\theta\lrcorner\omega^{n-1}$ is $\bar\partial$-exact, proving that $[\theta]\in H^{0, \, 1}(X, \, T^{1, \, 0}X)_{[\omega^{n-1}]}$. This proves the inclusion ``$\subset$''.

  Proving the reverse inclusion ``$\supset$'' in (\ref{eqn:pol-copol-equality}) takes more work. Let us consider the Lefschetz operator

\begin{eqnarray}\label{eqn:HL02forms}L_{\omega}^{n-2}:C^{\infty}_{0, \, 2}(X, \, \C)\rightarrow C^{\infty}_{n-2, \, n}(X, \, \C), \hspace{3ex}  \alpha\mapsto\omega^{n-2}\wedge\alpha,\end{eqnarray}

\noindent of multiplication by $\omega^{n-2}$ which is well known to be an isomorphism for any Hermitian (even non-K\"ahler or non-balanced) metric $\omega$ (see e.g. [Voi02, lemma 6.20, p. 146]). We clearly have $\theta\lrcorner\omega^{n-1}=(n-1)\,L_{\omega}^{n-2}(\theta\lrcorner\omega)$ by (\ref{eqn:contr-Leinniz}).

 The next lemma explains how the three-space decomposition (w.r.t. $\omega$) 

$$C^{\infty}_{0, \, 2}(X, \, \C)=\ker\Delta''\oplus\mbox{Im}\,\bar\partial\oplus\mbox{Im}\,\bar\partial^{\star}$$

\noindent transforms under $L_{\omega}^{n-2}$ and compares to the analogous decomposition of $C^{\infty}_{n-2, \, n}(X, \, \C)$. Note that in $C^{\infty}_{n-2, \, n}(X, \, \C)$ the subspace $\mbox{Im}\,\bar\partial^{\star}$ is reduced to zero for bidegree reasons.

\begin{Lem}\label{Lem:3space-decomp-Ln-2} If $\omega$ is a K\"ahler metric on a compact complex manifold $X$ with $n=\mbox{dim}_{\C}X$, then the operator (\ref{eqn:HL02forms}) satisfies

\begin{eqnarray}\label{eqn:3space-decomp-Ln-2}L_{\omega}^{n-2}(\ker\Delta'')=\ker\Delta'' \hspace{3ex}\mbox{and}\hspace{3ex} L_{\omega}^{n-2}(\mbox{Im}\,\bar\partial\oplus\mbox{Im}\,\bar\partial^{\star})=\mbox{Im}\,\bar\partial.\end{eqnarray}

\end{Lem}

 This will follow from two formulae that have an interest of their own.

\begin{Lem}\label{Lem:dbar-star-wedge} If $\omega$ is K\"ahler, then for every $\alpha\in C^{\infty}_{0, \, 2}(X, \, \C)$ we have

\begin{eqnarray}\label{eqn:dbar-star-wedge}\bar\partial^{\star}(\omega^{n-2}\wedge\alpha)=\omega^{n-2}\wedge\bar\partial^{\star}\alpha + (n-2)\,\omega^{n-3}\wedge i\partial\alpha.\end{eqnarray}

\end{Lem}

\noindent {\it Proof.} Using the K\"ahler commutation relation $\bar\partial^{\star}=-i\,[\Lambda, \, \partial],$ we get 

\begin{eqnarray}\label{eqn:dbar-star-wedge-gen}\bar\partial^{\star}(\omega^{n-2}\wedge\alpha)= -i\,\Lambda(\omega^{n-2}\wedge\partial\alpha) + i\,\partial(\Lambda(\omega^{n-2}\wedge\alpha)).\end{eqnarray}

 In the first term on the right-hand side of (\ref{eqn:dbar-star-wedge-gen}), we have

\begin{eqnarray}\label{eqn:dbar-star-wedge-1}\Lambda(\omega^{n-2}\wedge\partial\alpha)= [\Lambda, \, L^{n-2}](\partial\alpha) + \omega^{n-2}\wedge\Lambda(\partial\alpha) = \omega^{n-2}\wedge\Lambda(\partial\alpha).\end{eqnarray}

\noindent The last identity follows from the well-known formula (cf. [Voi02, p. 148])\!\!: 

\begin{eqnarray}\label{eqn:LrOmega}[L^r,\,\Lambda]=r(k-n+r-1)L^{r-1}  \hspace{2ex}\mbox{on}\,\,k\mbox{-forms,} \hspace{1ex} \mbox{for every}\,\,r,\end{eqnarray} 

\noindent which, when applied with $r=n-2$ to the $3$-form $\partial\alpha$, gives $[\Lambda, \, L^{n-2}](\partial\alpha)=0$.

 In the second term on the right-hand side of (\ref{eqn:dbar-star-wedge-gen}), we have

$$\Lambda(\omega^{n-2}\wedge\alpha)=[\Lambda, \, L^{n-2}](\alpha) + \omega^{n-2}\wedge\Lambda(\alpha) = (n-2)\,\omega^{n-3}\wedge\alpha + \omega^{n-2}\wedge\Lambda(\alpha),$$

\noindent where the last identity follows again from (\ref{eqn:LrOmega}) applied with $r=n-2$ to the $2$-form $\alpha$. (Note that $\Lambda\alpha=0$, but we ignore this here.) Taking $\partial$ on either side of the above identity and using the K\"ahler assumption on $\omega$, we get

\begin{eqnarray}\label{eqn:dbar-star-wedge-2}\partial(\Lambda(\omega^{n-2}\wedge\alpha)) = (n-2)\,\omega^{n-3}\wedge\partial\alpha + \omega^{n-2}\wedge\partial\Lambda(\alpha),\end{eqnarray}

 Thus, putting (\ref{eqn:dbar-star-wedge-1}) and (\ref{eqn:dbar-star-wedge-2}) together, we see that (\ref{eqn:dbar-star-wedge-gen}) transforms to

\begin{eqnarray}\nonumber\bar\partial^{\star}(\omega^{n-2}\wedge\alpha) & = & -i\,\omega^{n-2}\wedge\Lambda(\partial\alpha) + (n-2)\,\omega^{n-3}\wedge i\partial\alpha + \omega^{n-2}\wedge i\partial\Lambda(\alpha)\\
\nonumber & = & \omega^{n-2}\wedge i\,[\partial, \,\Lambda](\alpha) + (n-2)\,\omega^{n-3}\wedge i\partial\alpha\\
\nonumber & = &  \omega^{n-2}\wedge\bar\partial^{\star}\alpha + (n-2)\,\omega^{n-3}\wedge i\partial\alpha.\end{eqnarray}

\noindent This is what we had set out to prove. Note that we have used again the K\"ahler commutation relation $i\,[\partial, \,\Lambda] = -i\, [\Lambda, \, \partial] = \bar\partial^{\star}$.   \hfill $\Box$

\vspace{2ex}

 The next formula we need is the following.

\begin{Lem}\label{Lem:Delta''-wedge} If $\omega$ is K\"ahler, then for every $\alpha\in C^{\infty}_{0, \, 2}(X, \, \C)$ we have

\begin{eqnarray}\label{eqn:Delta''-wedge}\Delta_{\omega}''(\omega^{n-2}\wedge\alpha)=\omega^{n-2}\wedge\Delta_{\omega}''\alpha.\end{eqnarray}

\end{Lem}

\noindent {\it Proof.} This is an immediate consequence of the commutation property

$$[L_{\omega}, \, \Delta_{\omega}'']=0, \hspace{3ex}\mbox{hence}\hspace{3ex} [L_{\omega}^k, \, \Delta_{\omega}'']=0 \hspace{2ex}\mbox{for all}\,\,k,$$

\noindent which in turn follows from the K\"ahler identities. Alternatively, we can use Lemma \ref{Lem:dbar-star-wedge} and the K\"ahler identities to give a direct proof as follows. Since $\bar\partial(\omega^{n-2}\wedge\alpha)=0$ for bidegree reasons, $\Delta_{\omega}''(\omega^{n-2}\wedge\alpha)$ reduces to its first term, so using (\ref{eqn:dbar-star-wedge}) we get

\begin{eqnarray}\nonumber\Delta_{\omega}''(\omega^{n-2}\wedge\alpha) & = & \bar\partial\bar\partial^{\star}(\omega^{n-2}\wedge\alpha)=\bar\partial(\omega^{n-2}\wedge\bar\partial^{\star}\alpha + (n-2)\,\omega^{n-3}\wedge i\partial\alpha)\\
 \label{eqn:Delta''-wedge-gen} & = & \omega^{n-2}\wedge\bar\partial\bar\partial^{\star}\alpha + (n-2)\,\omega^{n-3}\wedge i\bar\partial\partial\alpha.\end{eqnarray}

\noindent Now, using the K\"ahler identity $\bar\partial^{\star}=-i\,[\Lambda, \, \partial]$, we get

\begin{eqnarray}\nonumber\omega^{n-2}\wedge\bar\partial^{\star}\bar\partial\alpha & = & -i\,\omega^{n-2}\wedge [\Lambda, \, \partial]\bar\partial\alpha = -i\,\omega^{n-2}\wedge\Lambda(\partial\bar\partial\alpha) + i\,\omega^{n-2}\wedge\partial\Lambda(\bar\partial\alpha)\\
\label{eqn:Delta''-wedge-1} & = & -i\,\omega^{n-2}\wedge\Lambda(\partial\bar\partial\alpha)\end{eqnarray}

\noindent because $\bar\partial\alpha$ is of type $(0, \, 3)$, so $\Lambda(\bar\partial\alpha)=0$ for bidegree reasons. Meanwhile

\begin{eqnarray}\nonumber \omega^{n-2}\wedge\Lambda(\partial\bar\partial\alpha) & = & [L^{n-2}, \, \Lambda](\partial\bar\partial\alpha) + \Lambda(\omega^{n-2}\wedge\partial\bar\partial\alpha) = [L^{n-2}, \, \Lambda](\partial\bar\partial\alpha)\\
\label{eqn:Delta''-wedge-2}  & = & (n-2)\,\omega^{n-3}\wedge\partial\bar\partial\alpha.\end{eqnarray}

\noindent The second identity on the top line above follows from $\omega^{n-2}\wedge\partial\bar\partial\alpha=0$ for bidegree reasons (since $\omega^{n-2}\wedge\partial\bar\partial\alpha$ is of type $(n-1, \, n+1)$, hence vanishes), while the last identity follows from formula (\ref{eqn:LrOmega}) with $r=n-2$ and $k=4$. 

\noindent The combined identities (\ref{eqn:Delta''-wedge-1}) and (\ref{eqn:Delta''-wedge-2}) yield

$$\omega^{n-2}\wedge\bar\partial^{\star}\bar\partial\alpha = -(n-2)\,\omega^{n-3}\wedge i\partial\bar\partial\alpha = (n-2)\,\omega^{n-3}\wedge i\bar\partial\partial\alpha.$$

\noindent This last identity combines with (\ref{eqn:Delta''-wedge-gen}) to prove the claim.  \hfill $\Box$

\vspace{2ex}

 We need yet another observation.

\begin{Lem}\label{Lem:Ln-2isom} For any Hermitian metric $\omega$ on $X$, the normalised Lefschetz operator

$$\frac{1}{(n-2)!}\,L_{\omega}^{n-2}:C^{\infty}_{0,\,2}(X,\,\C)\rightarrow C^{\infty}_{n-2,\,n}(X,\,\C)$$

\noindent is an isometry w.r.t. the $L^2$ scalar product induced by $\omega$ on scalar-valued forms.

\end{Lem}

\noindent {\it Proof.} We will show that for every $l=3, \dots , n$, the following formula holds

\begin{eqnarray}\label{eqn:L-isom}\langle\langle\omega^{n-2}\wedge\alpha,\,\omega^{n-2}\wedge\beta\rangle\rangle = (n-2)!\,\frac{(l-2)!}{(n-l)!}\,\langle\langle\omega^{n-l}\wedge\alpha,\,\omega^{n-l}\wedge\beta\rangle\rangle\end{eqnarray}

\noindent for all forms $\alpha, \beta\in C^{\infty}_{0,\,2}(X,\,\C)$. We have

\begin{eqnarray}\nonumber\langle\langle\omega^{n-2}\wedge\alpha,\,\omega^{n-2}\wedge\beta\rangle\rangle & = & \langle\langle\Lambda(\omega^{n-2}\wedge\alpha),\,\omega^{n-3}\wedge\beta\rangle\rangle\\
\nonumber & = & \langle\langle[\Lambda,\,L^{n-2}]\,\alpha,\,\omega^{n-3}\wedge\beta\rangle\rangle\\
\nonumber   & = & (n-2)\,\langle\langle\omega^{n-3}\wedge\alpha,\,\omega^{n-3}\wedge\beta\rangle\rangle,\end{eqnarray}

\noindent where in going from the first to the second line, we have used the identities $[\Lambda,\,L^{n-2}]\,\alpha = \Lambda(\omega^{n-2}\wedge\alpha) - \omega^{n-2}\wedge\Lambda\,\alpha = \Lambda(\omega^{n-2}\wedge\alpha)$ since $\Lambda\alpha=0$ for bidegree reasons, while in going from the second to the third line we have used formula (\ref{eqn:LrOmega}) with $r=n-2$ and the anti-commutation $[\Lambda,\,L^{n-2}] = - [L^{n-2},\,\Lambda]$. This proves (\ref{eqn:L-isom}) for $l=3$. We can now continue by induction on $l$. Suppose that (\ref{eqn:L-isom}) has been proved for $l$. We have

\begin{eqnarray}\nonumber\langle\langle\omega^{n-l}\wedge\alpha,\,\omega^{n-l}\wedge\beta\rangle\rangle & = & \langle\langle\Lambda(\omega^{n-l}\wedge\alpha),\,\omega^{n-l-1}\wedge\beta\rangle\rangle\\
\nonumber & = & \langle\langle[\Lambda,\,L^{n-l}]\,\alpha,\,\omega^{n-l-1}\wedge\beta\rangle\rangle\\
\nonumber   & = & (n-l)(l-1)\,\langle\langle\omega^{n-l-1}\wedge\alpha,\,\omega^{n-l-1}\wedge\beta\rangle\rangle\end{eqnarray}

\noindent by arguments similar to those above, where formula (\ref{eqn:LrOmega}) has been used with $r=n-l$. We thus obtain (\ref{eqn:L-isom}) with $l+1$ in place of $l$.
 
 It is now clear that (\ref{eqn:L-isom}) for $l=n$ proves the contention.  \hfill $\Box$

\vspace{3ex}

\noindent {\it End of proof of Lemma \ref{Lem:3space-decomp-Ln-2}.} Since the map $L_{\omega}^{n-2}$ of (\ref{eqn:HL02forms}) is an isomorphism, it follows from Lemma \ref{Lem:Delta''-wedge} that $L_{\omega}^{n-2}(\ker\Delta_{\omega}'') = \ker\Delta_{\omega}''$. Since $L_{\omega}^{n-2}$ maps any pair of orthogonal forms in $C^{\infty}_{0, \, 2}(X, \, \C)$ to a pair of orthogonal forms in $C^{\infty}_{n-2, \, n}(X, \, \C)$ by Lemma \ref{Lem:Ln-2isom}, it follows that the orthogonal complement of $\ker\Delta_{\omega}''$ in $C^{\infty}_{0, \, 2}(X, \, \C)$ (i.e. $\mbox{Im}\,\bar\partial\oplus\mbox{Im}\,\bar\partial^{\star}$) is isomorphic under $L_{\omega}^{n-2}$ to the orthogonal complement of $\ker\Delta_{\omega}''$ in $C^{\infty}_{n-2, \, n}(X, \, \C)$ (i.e. $\mbox{Im}\,\bar\partial$). Note that $\mbox{Im}\,\bar\partial^{\star}=0$ in $C^{\infty}_{n-2, \, n}(X, \, \C)$ for type reasons. The proof is complete.  \hfill $\Box$

\vspace{3ex}

\noindent {\it End of proof of Proposition \ref{Prop:pol-copol-equality}.} Recall that we have yet to prove the inclusion ``$\supset$'' in (\ref{eqn:pol-copol-equality}). Let $[\theta]\in H^{0, \, 1}(X, \,T^{1, \, 0}X)_{[\omega^{n-1}]}$. This means that $\theta\lrcorner\omega^{n-1}\in\mbox{Im}\,\bar\partial\subset C^{\infty}_{n-2, \, n}(X, \, \C)$ (cf. (\ref{eqn:co-pol-space})). Since $\theta\lrcorner\omega^{n-1} = (n-1)\,L_{\omega}^{n-2}(\theta\lrcorner\omega)$ (cf. (\ref{eqn:contr-Leinniz})) and $\theta\lrcorner\omega$ is of type $(0, \, 2)$, we get from Lemma \ref{Lem:3space-decomp-Ln-2} that

\begin{eqnarray}\label{eqn:theta-omega-Im-sum}\theta\lrcorner\omega\in\mbox{Im}\,\bar\partial\oplus\mbox{Im}\,\bar\partial^{\star}\subset C^{\infty}_{0, \, 2}(X, \, \C).\end{eqnarray}

\noindent On the other hand, $\bar\partial\theta=0$ (since $\theta$ represents a class $[\theta]\in H^{0, \, 1}(X, \,T^{1, \, 0}X)$) and $\bar\partial\omega=0$ (since $\omega$ is assumed K\"ahler). Hence $(ii)$ of Lemma \ref{Lem:dbar-contr-omega} gives

\begin{eqnarray}\label{eqn:theta-omega-ker-dbar}\bar\partial(\theta\lrcorner\omega) =0, \hspace{3ex}\mbox{i.e.}\hspace{3ex} \theta\lrcorner\omega\in\ker\bar\partial=\ker\Delta_{\omega}''\oplus\mbox{Im}\,\bar\partial\subset C^{\infty}_{0, \, 2}(X, \, \C).\end{eqnarray}

\noindent Since the three subspaces in the decomposition $C^{\infty}_{0, \, 2}(X, \, \C)= \ker\Delta_{\omega}''\oplus\mbox{Im}\,\bar\partial\oplus\mbox{Im}\,\bar\partial^{\star}$ are mutually orthogonal, (\ref{eqn:theta-omega-Im-sum}) and (\ref{eqn:theta-omega-ker-dbar}) imply that $\theta\lrcorner\omega\in\mbox{Im}\,\bar\partial$, i.e. $[\theta]\in H^{0, \, 1}(X, \,T^{1, \, 0}X)_{[\omega]}$ (cf. (\ref{eqn:pol-space})). \hfill $\Box$

\subsection{Primitive $(n-1, \, 1)$-classes on balanced manifolds}\label{subsection:prim-n-11}

 In the case of a K\"ahler class $[\omega]$, primitive Dolbeault cohomology classes of type $(n-1, \, 1)$ (for $[\omega]$) play a pivotal role in the theory of deformations of $X$ that are polarised by $[\omega]$ thanks to the isomorphism (\ref{eqn:pol-space-isom}) induced by the Calabi-Yau isomorphism. However, if $[\omega]$ is replaced by a balanced class $[\omega^{n-1}]$, primitive classes can no longer be defined in the standard way except in the case of $(1, \, 1)$-classes or, more generally, in that of De Rham $2$-classes (since the definition uses then the $(n-1)^{st}$ power of $\omega$ that is closed by the balanced assumption). In particular, defining an $(n-1, \, 1)$-class $[\alpha]$ as {\it primitive} by requiring that $\omega\wedge\alpha$ be $\bar\partial$-exact would be meaningless if $\omega$ is not closed since this definition would depend on the choice of representative $\alpha$ of the class $[\alpha]$. However, since the space $H^{0, \, 1}(X, \,T^{1, \, 0}X)_{[\omega^{n-1}]}$ carries over the meaning of $H^{0, \, 1}(X, \,T^{1, \, 0}X)_{[\omega]}$ to the balanced case, it is natural to make the following ad hoc definition in the balanced case.

\begin{Def}\label{Def:primitive-classes_n-11} Let $X$ be a compact balanced Calabi-Yau $\partial\bar\partial$-manifold ($n:=\mbox{dim}_{\C}X$). Fix a non-vanishing holomorphic $(n, \, 0)$-form $u$ and a balanced class $[\omega^{n-1}]$ on $X$. The space of {\bf primitive} classes of type $(n-1, \, 1)$ (for $[\omega^{n-1}]$) is defined as the image under the Calabi-Yau isomorphism

$$T_{[u]}: H^{0, \, 1}(X, \, T^{1, \, 0}X)\stackrel{\cdot\lrcorner [u]}{\longrightarrow} H^{n-1, \, 1}(X, \, \C)$$ 

\noindent in (\ref{eqn:u-isom-classes}) of the subspace $H^{0, \, 1}(X, \,T^{1, \, 0}X)_{[\omega^{n-1}]}\subset H^{0, \, 1}(X, \,T^{1, \, 0}X)$, i.e.

$$H^{n-1, \, 1}_{prim}(X, \, \C):=T_{[u]}\bigg(H^{0, \, 1}(X, \,T^{1, \, 0}X)_{[\omega^{n-1}]}\bigg)\subset H^{n-1, \, 1}(X, \, \C).$$

\noindent Explicitly, given the definition (\ref{eqn:co-pol-space}) of $H^{0, \, 1}(X, \,T^{1, \, 0}X)_{[\omega^{n-1}]}$, this means\!\!:

\begin{eqnarray}\label{eqn:explicit-prim-classes}[\theta\lrcorner u]\in H^{n-1, \, 1}_{prim}(X, \, \C) \hspace{3ex}\mbox{iff}\hspace{3ex} [\theta\lrcorner\omega^{n-1}]=0\in H^{n-2, \, n}(X, \, \C)\end{eqnarray}

\noindent for any class $[\theta]\in H^{0, \, 1}(X, \, T^{1, \, 0}X)$.

\end{Def}

 It is clear that $H^{n-1, \, 1}_{prim}(X, \, \C)$ does not depend on the choice of $u$ (which is unique up to a constant factor), but depends on the choice of balanced class $[\omega^{n-1}]$. When $\omega$ is K\"ahler, the ad hoc definition of $H^{n-1, \, 1}_{prim}(X, \, \C)$ coincides with the standard definition thanks to the isomorphism (\ref{eqn:pol-space-isom}) and to Proposition \ref{Prop:pol-copol-equality}. 

 Recall that unlike cohomology classes, primitive forms can be defined in the standard way for any Hermitian metric $\omega$\!\!: for any $k\leq n$, a $k$-form $\alpha$ on $X$ is primitive for $\omega$ if $\omega^{n-k+1}\wedge\alpha=0$. This condition is well known to be equivalent to $\Lambda_{\omega}\alpha=0$. No closedness assumption on $\omega$ is needed.

\vspace{3ex}

 In the rest of this subsection we shall investigate the extent to which the ad hoc primitive $(n-1, \, 1)$-classes defined by a balanced class retain the properties of primitive classes standardly defined by a K\"ahler class. We start with the form analogue of (\ref{eqn:explicit-prim-classes}). By the Calabi-Yau isomorphism (\ref{eqn:u-isom}), all $(n-1, \, 1)$-forms are of the shape $\theta\lrcorner u$ for some $\theta\in C^{\infty}_{0, \, 1}(X, \, T^{1, \, 0}X)$.

\begin{Lem}\label{Lem:prim-forms-omega} Let $(X, \, \omega)$ be an arbitrary Hermitian compact complex manifold ($n:=\mbox{dim}_{\C}X$) with $K_X$ trivial. Fix a non-vanishing holomorphic $(n, 0)$-form $u$. Then for any $\theta\in C^{\infty}_{0, \, 1}(X, \, T^{1, \, 0}X)$, the following equivalences hold\!\!:

\begin{eqnarray}\label{eqn:prim-forms-omega}\theta\lrcorner u \,\,\mbox{is primitive for}\,\,\omega \hspace{2ex} \Longleftrightarrow \theta\lrcorner\omega=0 \hspace{2ex} \Longleftrightarrow \theta\lrcorner\omega^{n-1}=0.\end{eqnarray}

\end{Lem}

\noindent {\it Proof.} By the definition of primitiveness, the $(n-1, \, 1)$-form $\theta\lrcorner u$ is primitive for $\omega$ if and only if $\omega\wedge(\theta\lrcorner u)=0$. Meanwhile

$$0=\theta\lrcorner(\omega\wedge u)=(\theta\lrcorner\omega)\wedge u + \omega\wedge(\theta\lrcorner u),$$

\noindent where the first identity holds for type reasons since the form $\omega\wedge u$ is of type $(n+1, \, 1)$, hence vanishes. Thus the vanishing of $\omega\wedge(\theta\lrcorner u)$ is equivalent to the vanishing of $(\theta\lrcorner\omega)\wedge u$ which, in turn, is equivalent to the vanishing of $\theta\lrcorner\omega$ as can be easily checked using the property $u\neq 0$ at every point of $X$. This proves the first equivalence in (\ref{eqn:prim-forms-omega}). The second equivalence follows from

$$\theta\lrcorner\omega^{n-1} = (n-1)\,\omega^{n-2}\wedge(\theta\lrcorner\omega)$$

\noindent (cf. (\ref{eqn:contr-Leinniz})) and from the map (\ref{eqn:HL02forms}) being an isomorphism.  \hfill $\Box$

\vspace{2ex}

 We have seen in Lemma \ref{Lem:d-closed-rep} that every Dolbeault cohomology class on a $\partial\bar\partial$-manifold can be represented by a $d$-closed form (which is, of course, not unique). The question we will now address is the following.

\begin{Question}\label{Question:prim-dclosed} Is it true that on a balanced Calabi-Yau $\partial\bar\partial$-manifold, every {\bf primitive} $(n-1, \, 1)$-class (in the sense of the ad hoc Definition \ref{Def:primitive-classes_n-11}) can be represented by a form that is both {\bf primitive} and $d$-{\bf closed}?

\end{Question}

 Should the answer to this question be affirmative, it would bear significantly on the discussion of Weil-Petersson metrics in $\S.$\ref{subsection:WP-metrics}. It is clear that in the K\"ahler case the answer is affirmative\!\!: the $\Delta''$-harmonic representative of any {\it primitive} (in the standard sense defined by the K\"ahler class = the ad hoc sense in the case of $(n-1,\,1)$-classes) $(p,\,q)$-class is both {\it primitive} and {\it $d$-closed}. We shall now see that the balanced case is far more complicated.

\begin{Lem}\label{Lem:prim-closed-harm} Let $(X, \, \omega)$ be a compact Hermitian manifold ($n:=\mbox{dim}_{\C}X$) and let $v$ be an arbitrary primitive form of type $(n-1, \, 1)$ on $X$. Then the following equivalences hold\!\!:

\begin{eqnarray}\label{eqn:prim-dstar-d}\bar\partial^{\star}v=0 \Longleftrightarrow \partial v=0  \hspace{3ex} \mbox{and} \hspace{3ex} \partial^{\star}v=0 \Longleftrightarrow \bar\partial v=0.\end{eqnarray}

\end{Lem}

\noindent {\it Proof.} It is well-known (cf. e.g. [Dem97, VI, $\S.$ 5.1]) that $\bar\partial^{\star} = -\star\partial\star$ and $\partial^{\star} = -\star\bar\partial\star$, where $\star : \Lambda^{p, \, q}T^{\star}X\longrightarrow\Lambda^{n-q, n-p}T^{\star}X$ is the Hodge star operator associated with $\omega$. On the other hand, the following formula for {\it primitive} forms is also known\!\!:

\begin{eqnarray}\label{eqn:prim-form-star-formula}\star\, v = i^{n^2+2n-2}\, v  \hspace{3ex}\mbox{for all}\,\,v\in C^{\infty}_{n-1, \, 1}(X, \, \C)_{prim}.\end{eqnarray}

\noindent (Recall that for {\it primitive} forms $v$ of arbitrary type $(p, \, q)$, the formula reads 

\begin{eqnarray}\label{eqn:prim-form-star-formula-gen}\star\, v = (-1)^{k(k+1)/2}\, i^{p-q}\, \frac{\omega^{n-p-q}\wedge v}{(n-p-q)!}, \hspace{2ex} \mbox{where}\,\, k:=p+q,\end{eqnarray}

\noindent see e.g. [Voi02, Proposition 6.29, p. 150].) Since $\star$ is an isomorphism, we see that the identity $\bar\partial^{\star}v=0$ is equivalent to $\partial(\star\, v)=0$, hence to $\partial v=0$ by (\ref{eqn:prim-form-star-formula}). The equivalence for $\partial^{\star}v=0$ is inferred similarly.  \hfill $\Box$

\begin{Cor}\label{Cor:prim-closed-harm} Under the assumptions of Lemma \ref{Lem:prim-closed-harm}, we have\!\!: \\

\noindent $(i)$\, if $v\in C^{\infty}_{n-1, \, 1}(X, \, \C)_{prim}$ and $\bar\partial v=0$, then

$$dv=0 \Longleftrightarrow \Delta''v=0.$$

\noindent $(ii)$\, if $v\in C^{\infty}_{n-1, \, 1}(X, \, \C)_{prim}$ and $\partial v=0$, then

$$dv=0 \Longleftrightarrow \Delta'v=0.$$

\noindent $(iii)$\, if $v\in C^{\infty}_{n-1, \, 1}(X, \, \C)_{prim}$ and $dv=0$, then

$$\Delta'v=0, \hspace{1ex} \Delta''v=0 \hspace{1ex}\mbox{and}\hspace{1ex} \Delta v=0.$$

\end{Cor}

\noindent {\it Proof.} Since $X$ is compact, we have $\ker\Delta''=\ker\bar\partial\cap\ker\bar\partial^{\star}$ and $\ker\Delta'=\ker\partial\cap\ker\partial^{\star}$. Since for any {\it pure-type} form $v$, the equivalence 

$$dv=0 \Longleftrightarrow \partial v=0 \hspace{1ex} \mbox{and} \hspace{1ex} \bar\partial v=0$$

\noindent holds, $(i)$ and $(ii)$ follow immediately from the two equivalences in (\ref{eqn:prim-dstar-d}). Now $(i)$ and $(ii)$ obviously give $\Delta'v=0$ and $\Delta''v=0$ under the assumptions of $(iii)$. To infer that $\Delta v=0$, it suffices to notice that for any {\it pure-type} form $v$ on a compact {\it Hermitian} manifold $(X, \, \omega)$, we have

\begin{eqnarray}\label{eqn:Delta-sum-eq}\langle\langle\Delta v, \, v\rangle\rangle = \langle\langle\Delta' v, \, v\rangle\rangle + \langle\langle\Delta'' v, \, v\rangle\rangle\end{eqnarray}

\noindent since $\langle\langle\Delta v, \, v\rangle\rangle = ||d\,v||^2 + ||d^{\star}v||^2$, $\langle\langle\Delta' v, \, v\rangle\rangle = ||\partial v||^2 + ||\partial^{\star}v||^2$ and $\langle\langle\Delta'' v, \, v\rangle\rangle$ $= ||\bar\partial v||^2 + ||\bar\partial^{\star}v||^2$, while $||dv||^2 = ||\partial v||^2 + ||\bar\partial v||^2$ (because $\partial v$ and $\bar\partial v$ are pure-type forms of different types, hence orthogonal) and similarly $||d^{\star}v||^2 = ||\partial^{\star} v||^2 + ||\bar\partial^{\star} v||^2$ (because $\partial^{\star} v$ and $\bar\partial^{\star} v$ are orthogonal for the same reason). Since $\Delta'v=0$ and $\Delta''v=0$, from (\ref{eqn:Delta-sum-eq}) we get $\langle\langle\Delta v, \, v\rangle\rangle = 0$ which amounts to $dv=0$ and $d^{\star}v=0$, hence to $\Delta v=0$.  \hfill $\Box$

\vspace{2ex}

 The conclusion $(iii)$ of the above Corollary \ref{Cor:prim-closed-harm} is that if an $(n-1, \, 1)$-form is both {\it primitive} and $d$-{\it closed}, it must be harmonic for each of the Laplacians $\Delta'$, $\Delta''$ and $\Delta$. Thus, if a representative that is both {\it primitive} and $d$-{\it closed} of a primitive $(n-1, \, 1)$-class exists, it can only be the $\Delta''$-harmonic representative. Fortunately we have

\begin{Lem}\label{Lem:harm-prim-harmall} Let $(X, \, \omega)$ be a compact Hermitian manifold ($n:=\mbox{dim}_{\C}X$). Suppose $v$ is a primitive $(n-1, \, 1)$-form such that $\Delta''v=0$. Then $\Delta' v=0$ and $\Delta v=0$. In particular, $dv=0$.

\end{Lem}

\noindent {\it Proof.} The assumption $\Delta''v=0$ means that $\bar\partial v=0$ and $\bar\partial^{\star}v=0$. Then $(i)$ of Corollary \ref{Cor:prim-closed-harm} implies that $dv=0$, i.e. $\partial v=0$. Then $(ii)$ of Corollary \ref{Cor:prim-closed-harm} ensures that $\Delta'v=0$. Then (\ref{eqn:Delta-sum-eq}) ensures that $\Delta v=0$. \hfill $\Box$

\vspace{2ex}

 Thus Question \ref{Question:prim-dclosed} reduces to whether on a balanced Calabi-Yau $\partial\bar\partial$-manifold $(X,\omega)$, the $\Delta''$-harmonic representative of any {\it primitive} $(n-1, \, 1)$-class (in the sense of the ad hoc Definition \ref{Def:primitive-classes_n-11}) is a {\it primitive} form. It will then also be $d$-closed by Lemma \ref{Lem:harm-prim-harmall}. Fix therefore a {\it primitive} $(n-1, \, 1)$-class $[\theta\lrcorner u]$ on $X$, where $[\theta]\in H^{0,\, 1}(X,\, T^{1,\, 0}X)$. By (\ref{eqn:explicit-prim-classes}), this means that

\begin{eqnarray}\label{eqn:prim-theta-u-recall}\theta\lrcorner\omega^{n-1}\in\mbox{Im}\,\bar\partial\end{eqnarray}

\noindent Suppose furthermore that $\Delta''(\theta\lrcorner u)=0$. The question is whether $\theta\lrcorner u$ is primitive, or equivalently (cf. (\ref{eqn:prim-forms-omega})) whether $\theta\lrcorner\omega^{n-1}=0$. Since $\ker\Delta''$ and $\mbox{Im}\,\bar\partial$ are orthogonal subspaces of $C_{n-2,\,n}^{\infty}(X, \,\C)$, (\ref{eqn:prim-theta-u-recall}) reduces the question to determining whether

\begin{eqnarray}\label{eqn:theta-omegan-1-harm}\Delta''(\theta\lrcorner\omega^{n-1})=0 \hspace{2ex} \mbox{or equivalently whether} \hspace{2ex} \bar\partial^{\star}(\theta\lrcorner\omega^{n-1})=0,\end{eqnarray}

\noindent since $\bar\partial(\theta\lrcorner\omega^{n-1})=0$ (trivially since $\theta\lrcorner\omega^{n-1}$ is of type $(n-2,\, n)$). 

 The next lemma transforms identity (\ref{eqn:theta-omegan-1-harm}) whose validity we are trying to determine.

\begin{Lem}\label{Lem:theta-omega-n-1-del} Let $X$ be a compact complex manifold ($\mbox{dim}_{\C}X=n$) equipped with an arbitrary Hermitian metric $\omega$. Fix any $\theta\in C^{\infty}_{0,\, 1}(X,\,T^{1,\, 0}X)$.

\noindent The following equivalence holds\!\!:

$$\bar\partial^{\star}(\theta\lrcorner\omega^{n-1})=0 \Longleftrightarrow \partial(\theta\lrcorner\omega)=0.$$ 

\end{Lem}

\noindent {\it Proof.} Formula (\ref{eqn:prim-form-star-formula-gen}) applied to the (primitive) $(0,\, 2)$-form $v:=\theta\lrcorner\omega$ reads\!\!:

\begin{eqnarray}\label{eqn:star-theta-omega}\star(\theta\lrcorner\omega) = \frac{\omega^{n-2}}{(n-2)!}\wedge(\theta\lrcorner\omega) = \theta\lrcorner\frac{\omega^{n-1}}{(n-1)!},  \hspace{2ex} \mbox{i.e.} \hspace{2ex} \star\bigg(\theta\lrcorner\frac{\omega^{n-1}}{(n-1)!}\bigg) = \theta\lrcorner\omega,\end{eqnarray}

\noindent having also used the property $\star^2=\mbox{Id}$ on $2$-forms. Now, $\bar\partial^{\star}=-\star\partial\star$, hence the condition $\bar\partial^{\star}(\theta\lrcorner\omega^{n-1})=0$ is equivalent to $\partial(\star(\theta\lrcorner\omega^{n-1}))=0$ which in turn is equivalent to $\partial(\theta\lrcorner\omega)=0$ by (\ref{eqn:star-theta-omega}). This proves the contention.  \hfill $\Box$

\vspace{2ex}

 However, we can see no reason why the desired condition $\partial(\theta\lrcorner\omega)=0$ should hold even if we exploit the assumption $\Delta''(\theta\lrcorner u)=0$. Note that if $\mbox{Ric}(\omega)=0$, by (\ref{eqn:Tu-dbarstar}) this assumption means that $\Delta''\theta=0$, i.e. $\bar\partial^{\star}\theta=0$ since we always have $\bar\partial\theta=0$. The most we can make of the property $\bar\partial^{\star}\theta=0$ is expressed in part $(ii)$ of the following lemma. Parts $(i)$ and $(iii)$ show that more can be said about {\it scalar-valued} $(0,\,1)$-forms $v$, although even if that information applied to the $T^{1,\,0}X$-valued $(0,\,1)$-form $\theta$, it would not suffice to deduce that $\partial(\theta\lrcorner\omega)=0$.

\begin{Lem}\label{Lem:01forms-delprim} Let $X$ be a compact complex manifold ($\mbox{dim}_{\C}X=n$) supposed to carry a {\bf balanced} metric $\omega$.

\noindent $(i)$\, For every $v\in C^{\infty}_{0,\, 1}(X,\,\C)$, the following equivalence holds\!\!:

$$\bar\partial^{\star}v=0 \Longleftrightarrow \partial v \hspace{1ex}\mbox{is primitive}.$$

\noindent $(ii)$\, For every $\theta\in C^{\infty}_{0,\, 1}(X,\,T^{1,\, 0}X)$, the following equivalence holds\!\!:

$$\bar\partial^{\star}\theta=0 \Longleftrightarrow (D'\theta)\wedge\omega^{n-1}=0\in C^{\infty}_{n,\, n}(X,\,T^{1,\, 0}X).$$

\noindent $(iii)$\, Suppose, furthermore, that $X$ is a $\partial\bar\partial$-manifold. Then, for every $v\in C^{\infty}_{0,\, 1}(X,\,\C)\cap\ker\bar\partial$, the following equivalence holds\!\!:

$$\Delta''v=0 \Longleftrightarrow \partial v=0 \hspace{2ex} (\Longleftrightarrow \Delta'v=0).$$

\end{Lem}

\noindent {\it Proof.} Since any $(0,\, 1)$-form is primitive, for $\star : C^{\infty}_{0,\, 1}(X,\,\C)\rightarrow C^{\infty}_{n-1,\, n}(X,\,\C)$ formula (\ref{eqn:prim-form-star-formula-gen}) reads

\begin{eqnarray}\label{eqn:prim-form-star-formula-01}\star\,v=i\,v\wedge\frac{\omega^{n-1}}{(n-1)!}, \hspace{3ex} v\in C^{\infty}_{0,\, 1}(X,\,\C).\end{eqnarray} 

\noindent Since $\bar\partial^{\star}=-\star\partial\star$, we see that the condition $\bar\partial^{\star}v=0$ is equivalent to $\partial(v\wedge\omega^{n-1})=0$. Since $\partial\omega^{n-1}=0$ (by the balanced assumption), the last identity is equivalent to $(\partial v)\wedge\omega^{n-1}=0$, which is precisely the condition that the $(1,\, 1)$-form $\partial v$ be primitive. This proves $(i)$. 

 The proof of $(ii)$ runs along the same lines as that of $(i)$ using the formula $\bar\partial^{\star}=-\star D'\star$ when $\bar\partial^{\star}$ acts on $T^{1,\, 0}X$-valued forms and $D'$ is the $(1,\, 0)$-component of the Chern connection $D$ of $(T^{1,\, 0}X,\,\omega)$. Indeed, formula (\ref{eqn:prim-form-star-formula-01}) still holds for $T^{1,\, 0}X$-valued $(0,\, 1)$-forms $\theta$ in place of $v$ and 

$$D'(\theta\wedge\omega^{n-1})=(D'\theta)\wedge\omega^{n-1} - \theta\wedge\partial\omega^{n-1} = (D'\theta)\wedge\omega^{n-1},$$ 

\noindent where the last identity follows from $\omega$ being balanced.

 To prove $(iii)$, fix an arbitrary form $v\in C^{\infty}_{0,\, 1}(X,\,\C)\cap\ker\bar\partial$. Since $\ker\Delta''=\ker\bar\partial\cap\ker\bar\partial^{\star}$, the condition $\Delta''v=0$ is equivalent for this $v$ to $\bar\partial^{\star}v=0$, which is equivalent to $\partial v$ being primitive by $(i)$. We are thus reduced to proving for this $v$ the equivalence\!\!: $\partial v$ is primitive $\Longleftrightarrow\,\,\partial v=0$.

 Notice that $\bar\partial(\partial v)=0$ thanks to the assumption $\bar\partial v=0$. Hence the pure-type form $\partial v$ is $d$-closed and $\partial$-exact, so by the $\partial\bar\partial$-lemma it must be $\partial\bar\partial$-exact\!\!:

$$\partial v = i\partial\bar\partial\varphi  \hspace{2ex}\mbox{for some}\,\, C^{\infty}\,\, \mbox{function}\hspace{1ex} \varphi : X\rightarrow \C.$$ 

\noindent Then we have the equivalences\!\!:

$$\partial v \hspace{1ex} \mbox{is primitive} \Longleftrightarrow \Lambda_{\omega}(i\partial\bar\partial\varphi)=0 \Longleftrightarrow \Delta_{\omega}\varphi=0 \Longleftrightarrow \varphi \hspace{1ex}\mbox{is constant },$$

\noindent where the last equivalence follows by the maximum principle from $X$ being compact. Meanwhile, $\varphi$ being constant is equivalent to the vanishing of $i\partial\bar\partial\varphi$, hence to the vanishing of $\partial v$. \hfill $\Box$

\vspace{3ex}

 The conclusion of these considerations is that Question \ref{Question:prim-dclosed} may have a negative answer in general in the balanced case. Let us now notice that even the answer to the following weaker question may be negative in the balanced case.

\begin{Question}\label{Question:prim-rep-prim-class} Is it true that on a balanced Calabi-Yau $\partial\bar\partial$-manifold, every {\bf primitive} $(n-1, \, 1)$-class (in the sense of the ad hoc Definition \ref{Def:primitive-classes_n-11}) can be represented by a {\bf primitive} form?

\end{Question}

 Let $[\theta\lrcorner u]\in H^{n-1,\,1}_{prim}(X,\,\C)$ be a {\it primitive} class in the ad hoc sense. This means that $\theta\lrcorner\omega^{n-1}$ is $\bar\partial$-exact (for any representative $\theta$ of the class $[\theta]\in H^{0,\,1}(X,\,T^{1,\,0}X)_{[\omega^{n-1}]}$). Pick any representative $\theta$ and any $\bar\partial$-potential $w\in C^{\infty}_{n-2,\,n-1}(X,\,\C)$ of $\theta\lrcorner\omega^{n-1}$, i.e. $\bar\partial w=\theta\lrcorner\omega^{n-1}$. Since

$$L_{\omega}^{n-3}:C^{\infty}_{1,\,2}(X,\,\C)\rightarrow C^{\infty}_{n-2,\,n-1}(X,\,\C), \hspace{3ex} \alpha\mapsto\omega^{n-3}\wedge\alpha,$$

\noindent is an isomorphism (see e.g. [Voi02, lemma 6.20, p. 146]), since there is a Lefschetz decomposition (cf. [Voi02, proposition 6.22, p. 147])

$$\Lambda^{1,\,2}=\Lambda^{1,\,2}_{prim}\oplus\bigg(\omega\wedge\Lambda^{0,\,1}\bigg)$$

\noindent and since every $C^{\infty}$ $(0,\,1)$-form can be written as $(n-1)\,\xi\lrcorner\omega$ for a unique vector field $\xi\in C^{\infty}(X,\,T^{1,\,0}X)$ (because $\omega$ is non-degenerate), we see that there is a unique {\it primitive} $C^{\infty}$ form $\alpha_0$ of type $(1,\,2)$ and a unique $C^{\infty}$ vector field $\xi$ of type $(1,\,0)$ such that

\begin{eqnarray}\label{eqn:w-decomposition}w=\omega^{n-3}\wedge\bigg(\alpha_0 + (n-1)\,\omega\wedge(\xi\lrcorner\omega)\bigg) = \omega^{n-3}\wedge\alpha_0 + \xi\lrcorner\omega^{n-1}.\end{eqnarray}

\noindent Consequently, $\theta\lrcorner\omega^{n-1}= \bar\partial w = \bar\partial(\omega^{n-3}\wedge\alpha_0) + (\bar\partial\xi)\lrcorner\omega^{n-1}$ since $ \bar\partial(\xi\lrcorner\omega^{n-1})=(\bar\partial\xi)\lrcorner\omega^{n-1} - \xi\lrcorner(\bar\partial\omega^{n-1})$ (cf. $(i)$ of Lemma \ref{Lem:dbar-contr-omega}) and here $\bar\partial\omega^{n-1}=0$ by the balanced assumption on $\omega$. Thus we get

$$(\theta-\bar\partial\xi)\lrcorner\omega^{n-1}=\bar\partial(\omega^{n-3}\wedge\alpha_0).$$

\noindent We see that $\theta-\bar\partial\xi$ represents the class $[\theta]\in H^{0,\,1}(X,\,T^{1,\,0}X)_{[\omega^{n-1}]}$, so $(\theta-\bar\partial\xi)\lrcorner u$ represents the class $[\theta\lrcorner u]\in H^{n-1,\,1}_{prim}(X,\,\C)$. We know from Lemma \ref{Lem:prim-forms-omega} that the primitivity condition on the form $(\theta-\bar\partial\xi)\lrcorner u$ is equivalent to $(\theta-\bar\partial\xi)\lrcorner\omega^{n-1}=0$, i.e. to $\bar\partial(\omega^{n-3}\wedge\alpha_0)=0$ in this case. However, we can see no reason why this vanishing should occur, part of the obstruction being the primitive $(1,\,2)$-form $\alpha_0$. 

 Thus in the balanced, non-K\"ahler case, the answer to Question \ref{Question:prim-rep-prim-class} may be negative in general.

\section{Period map and Weil-Petersson metrics}\label{section:period-WP}

 We now fix an arbitrary balanced Calabi-Yau $\partial\bar\partial$-manifold $X$, $\mbox{dim}_{\C}X=n$. All the fibres $(X_t)_{t\in\Delta}$ in the Kuranishi family of $X=X_0$ are again balanced Calabi-Yau $\partial\bar\partial$-manifolds if $t$ is sufficiently close to $0\in\Delta$. This follows from Wu's theorem in [Wu06] and from the deformation openness of the triviality of the canonical bundle $K_{X_t}$ when the dimension of $H^{n,\,0}(X_t,\,\C)$ is locally independent of $t$ (as the $\partial\bar\partial$ assumption ensures this to be the case here). Thus $H^{n,\, 0}(X_t,\,\C)$ is a complex line varying holomorphically with $t$ inside the fixed complex vector space $H^n(X,\,\C)$. The canonical injection $H^{n,\, 0}(X_t,\,\C)\subset H^n(X,\,\C)$ is induced by the $\partial\bar\partial$-lemma property of $X_t$ (cf. Lemma \ref{Lem:d-closed-rep} and comments thereafter). The associated period map $\Delta\ni t\mapsto H^{n,\, 0}(X_t,\,\C)$ takes values in the complex projective space $\Proj H^n(X,\,\C)$ after identifying each complex line $H^{n,\, 0}(X_t,\,\C)$ with the point it defines therein.

\subsection{Period domain and the local Torelli theorem}\label{subsection:period-domain}

 Most of the material in this subsection before Theorem \ref{The:localTorelli} is essentially known, but we take this oportunity to stress that only minimal assumptions are needed and to fix the notation for the rest of the paper.

 Let $\omega$ be a Hermitian metric on $X$. All the formal adjoint operators and Laplacians will be calculated w.r.t. $\omega$. The Hodge $\star$-operator defined by $\omega$ on $n$-forms

$$\star \,\,:\,\, C^{\infty}_n(X,\,\C)\longrightarrow C^{\infty}_n(X,\,\C)$$

\noindent satisfies $\star^2=(-1)^n$, so it induces a decomposition

\begin{eqnarray}\label{eqn:Lambdan-decomposition}C_n^{\infty}(X,\,\C) =  \Lambda^n_{+}\oplus \Lambda^n_{-},\end{eqnarray}

\noindent  where $\Lambda^n_{\pm}$ stand for the eigenspaces of $\star$ corresponding to the eigenvalues $\pm 1$ (if $n$ is even), $\pm i$ (if $n$ is odd). This decomposition is easily seen to be orthogonal for the $L^2$ scalar product induced by $\omega$\!: for any $u\in\Lambda^n_{+}$ and any $v\in\Lambda^n_{-}$, one easily checks that $\langle\langle u,\, v\rangle\rangle = - \langle\langle u,\, v\rangle\rangle$ by writing $u=\star u$ (if $n$ is even) and $u=-i\,(\star u)$ (if $n$ is odd) and using the easy-to-check identity $\langle\langle\star u,\, v\rangle\rangle = (-1)^n\, \langle\langle u,\, \star v\rangle\rangle$ for any $n$-forms $u, v$.

When $\star$ is restricted to $\Delta$-harmonic forms, it assumes $\Delta$-harmonic values\!:

$$\star \,\,:\,\, {\cal H}^n_{\Delta}(X,\,\C)\longrightarrow {\cal H}^n_{\Delta}(X,\,\C)$$

\noindent since $\Delta:=dd^{\star} + d^{\star}d$ commutes with $\star$ as is well known to follow from the standard formula $d^{\star}=-\star\,d\,\star$. Thus the Hodge isomorphism $H^n(X,\,\C)\simeq {\cal H}^n_{\Delta}(X,\,\C)$ mapping any De Rham class to its $\Delta$-harmonic representative extends the definition of $\star$ to the De Rham cohomology of degree $n$\!\!:

\begin{eqnarray}\label{eqn:star-cohom-def}\star \,\,:\,\, H^n(X,\,\C)\longrightarrow H^n(X,\,\C)\end{eqnarray}

\noindent and we get a decomposition in cohomology analogous to (\ref{eqn:Lambdan-decomposition})\!:

\begin{eqnarray}\label{eqn:Hn-decomposition}H^n(X,\,\C) =  H^n_{+}(X,\,\C)\oplus H^n_{-}(X,\,\C),\end{eqnarray}

\noindent where $H^n_{\pm}(X,\,\C)$ are the eigenspaces of $\star$ corresponding to the eigenvalues $\pm 1$ (if $n$ is even), $\pm i$ (if $n$ is odd). Thus $H^n_{+}(X,\,\C)$ (resp. $H^n_{-}(X,\,\C)$) consists of the De Rham classes $\{\alpha\}$ of degree $n$ whose $\Delta$-harmonic representative $\alpha$ lies in $\Lambda^n_{+}$ (resp. $\Lambda^n_{-}$). Note that no assumption whatsoever (either K\"ahler or balanced) is needed on the Hermitian metric $\omega$.

\vspace{2ex}

 On the other hand, the Hodge-Riemann bilinear form can always be defined on the De Rham cohomology of degree $n$\!\!:

\begin{eqnarray}\label{eqn:Q-def}\nonumber Q & : & H^n(X,\,\C)\times H^n(X,\,\C)\longrightarrow\C,\\
         & &  (\{\alpha\}, \{\beta\})\longmapsto (-1)^{\frac{n(n-1)}{2}}\,\int\limits_X \alpha\wedge\beta :=Q(\{\alpha\},\,\{\beta\}).\end{eqnarray}

\noindent It is clear that $Q(\cdot,\,\cdot)$ is independent of the choice of representatives $\alpha$ and $\beta$ of the respective De Rham classes of degree $n$ since no power of $\omega$ is involved in the definition of $Q$, so no K\"ahler or balanced or any other assumption is needed on $\omega$ unlike the case of the De Rham cohomology in degree $k<n$. Thus $Q$ is independent of $\omega$ and of the complex structure of $X$, depending only on the differential structure of $X$. It is also clear that $Q$ is non-degenerate since for any $\Delta$-harmonic $n$-form $\alpha$, $\star\bar\alpha$ is again $\Delta$-harmonic and

\begin{eqnarray}\nonumber Q(\{\alpha\},\,\{\star\bar\alpha\}) & = & (-1)^{\frac{n(n-1)}{2}}\,\int\limits_X \alpha\wedge\star\bar\alpha=(-1)^{\frac{n(n-1)}{2}}\,\int\limits_X\langle\alpha,\,\alpha\rangle_{\omega}\,dV_{\omega}\\
\nonumber  & = & (-1)^{\frac{n(n-1)}{2}}\,||\alpha||^2_{\omega}\neq 0 \hspace{3ex} \mbox{if}\hspace{1ex}\alpha\neq 0.\end{eqnarray}

\noindent Hence the associated sesquilinear form

\begin{eqnarray}\label{eqn:H-def}\nonumber H & : & H^n(X,\,\C)\times H^n(X,\,\C)\longrightarrow\C,\\
         & &  (\{\alpha\}, \{\beta\})\longmapsto (-1)^{\frac{n(n+1)}{2}}\,i^n\,\int\limits_X \alpha\wedge\bar\beta = (-i)^n\,Q(\{\alpha\},\,\{\bar\beta\})\end{eqnarray}

\noindent is non-degenerate.

\begin{Lem}\label{Lem:H-pos} $(a)$\, $H(\{\alpha\},\,\{\alpha\})>0$ for every class $\{\alpha\}\in H^n_{+}(X,\,\C)\setminus\{0\}$. Hence $H$ defines a positive definite sesquilinear form (i.e. a Hermitian metric) on $H^n_{+}(X,\,\C)$.

\noindent $(b)$\, $H(\{\alpha\},\,\{\alpha\})<0$ for every class $\{\alpha\}\in H^n_{-}(X,\,\C)\setminus\{0\}$.

\noindent $(c)$\, $H(\{\alpha\},\,\{\beta\})=0$ for every class $\{\alpha\}\in H^n_{+}(X,\,\C)$ and every class $\{\beta\}\in H^n_{-}(X,\,\C)$. Hence the decomposition (\ref{eqn:Hn-decomposition}) is orthogonal for $H$. 

\end{Lem}

\noindent {\it Proof.} $(a)$\, Let $\alpha$ be a $\Delta$-harmonic $n$-form such that the class $\{\alpha\}\in H^n_{+}(X,\,\C)$.    

\noindent If $n$ is even, $\star\,\alpha=\alpha$, hence taking conjugates we get $\star\,\bar\alpha=\bar\alpha$. Thus

$$H(\{\alpha\},\,\{\alpha\})=(-1)^{\frac{n(n+1)}{2}}\,i^n\,\int\limits_X \alpha\wedge\star\bar\alpha = \int\limits_X |\alpha|^2_{\omega}\,dV_{\omega}=||\alpha||^2_{\omega}>0 $$

\noindent if $\alpha\neq 0$, since $(-1)^{\frac{n(n+1)}{2}}\,i^n=i^{n^2 + 2n} = 1$ when $n$ is even. (Indeed, $n^2 + 2n\in 4\Z$ when $n$ is even.)

\noindent If $n$ is odd, $\star\,\alpha=i\,\alpha$, hence taking conjugates we get $\star\,\bar\alpha=-i\,\bar\alpha$. Equivalently, $\bar\alpha=i\star\bar\alpha$. On the other hand, $(-1)^{\frac{n(n+1)}{2}}\,i^n=i^{n^2 + 2n} = -i$ when $n$ is odd since $n^2+2n\in 4\Z+3$ in this case. We then get as above that again $H(\{\alpha\},\,\{\alpha\})=||\alpha||^2_{\omega}>0$ if $\alpha\neq 0$. This proves $(a)$. The proof of $(b)$ is very similar and is left to the reader.

$(c)$\, Let $\alpha$ and $\beta$ be $\Delta$-harmonic $n$-forms such that $\{\alpha\}\in H^n_{+}(X,\,\C)$ and $\{\beta\}\in H^n_{-}(X,\,\C)$. If $n$ is even, this means that $\star\,\alpha=\alpha$ and $\star\,\beta=-\beta$. Using the property $\star\,\beta=-\beta$, we get

\begin{eqnarray}\label{H-perp2}H(\{\alpha\},\,\{\beta\})=-(-1)^{\frac{n(n+1)}{2}}\,i^n\,\int\limits_X \alpha\wedge\star\bar\beta = -(-1)^{\frac{n(n+1)}{2}}\,i^n\,\langle\langle\alpha,\,\beta\rangle\rangle_{\omega},\end{eqnarray}

\noindent while using the property $\star\,\alpha=\alpha$, we get

\begin{eqnarray}\label{H-perp2'}\nonumber H(\{\alpha\},\,\{\beta\}) & = & (-1)^{\frac{n(n+1)}{2}}\,i^n\,\int\limits_X \star\,\alpha\wedge\bar\beta = (-1)^{n^2}(-1)^{\frac{n(n+1)}{2}}\,i^n\,\int\limits_X \bar\beta\wedge\star\,\alpha \\
  & = & (-1)^{\frac{n(n+1)}{2}}\,i^n\,\int\limits_X \overline{\langle\beta,\,\alpha\rangle_{\omega}}\, dV_{\omega} = (-1)^{\frac{n(n+1)}{2}}\,i^n\,\langle\langle\alpha,\,\beta\rangle\rangle_{\omega},\end{eqnarray}

\noindent having used the fact $(-1)^{n^2}=1$ since $n$ is even and the identity $\overline{\langle\beta,\,\alpha\rangle_{\omega}}= \langle\alpha,\,\beta\rangle_{\omega}$. The expressions (\ref{H-perp2}) and (\ref{H-perp2'}) for $H(\{\alpha\},\,\{\beta\})$ are now seen to differ only by a sign, hence $H(\{\alpha\},\,\{\beta\})=0$.

\noindent When $n$ is odd, we have $\star\,\alpha=i\,\alpha$ (hence $\alpha=-i\star\alpha$) and $\star\,\beta=-i\,\beta$ (hence $\bar\beta=-i\star\bar\beta$). Using the former and then the latter of these two pieces of information, we get as above two expressions for $H(\{\alpha\},\,\{\beta\})$ that differ only by a sign. Hence $H(\{\alpha\},\,\{\beta\})=0$.  \hfill $\Box$

\vspace{2ex}

 We now bring in the complex structure of $X$ (that is supposed to have the $\partial\bar\partial$ property which induces the inclusion $H^{n,\,0}(X,\,\C)\subset H^n(X,\,\C)$).

\begin{Lem}\label{Lem:Hn0-pos} Let $X$ be a compact complex $\partial\bar\partial$-manifold ($\mbox{dim}_{\C}X=n$). Then the following inclusions hold\!\!:

\vspace{2ex}

\noindent $H^{n,\,0}(X,\,\C)\subset H^n_{+}(X,\,\C) \hspace{1ex}\mbox{if}\,\,n\,\,\mbox{is {\bf even},} \hspace{1ex} H^{n,\,0}(X,\,\C)\subset H^n_{-}(X,\,\C) \hspace{1ex}\mbox{if}\,\,n\,\,\mbox{is {\bf odd}}.$

\vspace{2ex}

\noindent In particular, the restriction $H : H^{n,\,0}(X,\,\C)\times H^{n,\,0}(X,\,\C)\rightarrow \C$ of $H$ to $H^{n,\,0}(X,\,\C)$ is {\bf positive definite if} $n$ {\bf is even} and is {\bf negative definite if} $n$ {\bf is odd} thanks to Lemma \ref{Lem:H-pos} (hence we get a Hermitian metric on $H^{n,\,0}(X,\,\C)$ defined by the scalar product induced by $H$ when $n$ is even and by $-H$ when $n$ is odd).

\end{Lem}

Before proving this statement, we make a trivial but useful observation.

\begin{Lem}\label{Lem:n0delta''-'} Let $(X,\,\omega)$ be any compact complex Hermitian manifold ($\mbox{dim}_{\C}X=n$). For every $(n,\,0)$-form $\alpha$, the following equivalence and implication hold\!\!:

$$\Delta''\alpha = 0 \Longleftrightarrow \Delta'\alpha=0 \Longrightarrow \Delta\alpha=0.$$

\end{Lem}

\noindent {\it Proof.} Since $X$ is compact, $\ker\,\Delta''=\ker\bar\partial\cap\ker\bar\partial^{\star}$ and $\ker\,\Delta'=\ker\partial\cap\ker\partial^{\star}$. However, $\partial\alpha=0$ and $\bar\partial^{\star}\alpha=0$ for any $(n,\,0)$-form $\alpha$ for trivial bidegree reasons. Hence, for any $\alpha\in C^{\infty}_{n,\,0}(X,\,\C)$,  the following equivalences hold\!\!:

 $$\Delta'\alpha=0 \Leftrightarrow \partial^{\star}\alpha=0  \hspace{3ex} \mbox{and} \hspace{3ex} \Delta''\alpha=0 \Leftrightarrow \bar\partial\alpha=0.$$

\noindent Consequently, from the identity $\partial^{\star} = -\star\bar\partial\star$ (cf. e.g. [Dem97, VI, $\S. 5.1$]) and from the fact that $\star$ is an isomorphism, we get the equivalence\!\!: $\Delta'\alpha=0 \,\,\Leftrightarrow\,\,\bar\partial(\star\alpha)=0$. Since $\alpha$ is of type $(n,\,0)$, it is primitive (w.r.t. any metric, hence also w.r.t. $\omega$), so formula (\ref{eqn:prim-form-star-formula-gen}) applied to $\alpha$ reads\!\!: $\star\,\alpha= (-1)^{n(n+1)/2}\,i^n\,\alpha$. Thus the previous equivalence implies the following equivalence\!\!:

 $$\Delta'\alpha=0 \,\,\Leftrightarrow\,\,\bar\partial\alpha=0,$$  

\noindent while the equivalence $\bar\partial\alpha=0\,\,\Leftrightarrow\,\,\Delta''\alpha=0$ has already been observed. We have thus proved the equivalence claimed in the statement. The implication claimed in the statement now follows from identity (\ref{eqn:Delta-sum-eq}) applied to the pure-type form $\alpha$ and the fact that $\langle\langle\Delta\alpha,\,\alpha\rangle\rangle\geq 0$ with equality if and only if $\Delta\alpha=0$.  \hfill $\Box$

\vspace{3ex}

\noindent {\it Proof of Lemma \ref{Lem:Hn0-pos}.} Let $[\alpha]\in H^{n,\,0}(X,\,\C)$ be an arbitrary Dolbeault cohomology class of type $(n,\,0)$. Since the only $\bar\partial$-exact form of type $(n,\,0)$ is the zero form, the class $[\alpha]$ contains a unique representative $\alpha$. Clearly, $\alpha$ is of type $(n,\,0)$ and $\Delta''$-harmonic, so from Lemma \ref{Lem:n0delta''-'} we get $\Delta\alpha=0$. On the other hand, formula (\ref{eqn:prim-form-star-formula-gen}) applied to $\alpha$ (which is primitive since it is of type $(n,\, 0)$) reads\!: $\star\,\alpha = (-1)^{n(n+1)/2}\,i^n\alpha=i^{n(n+2)}\,\alpha$. Hence, if $n$ is even, $\alpha\in\Lambda^n_{+}$ since $i^{n(n+2)}=1$, while if $n$ is odd, $\alpha\in\Lambda^n_{-}$ since $i^{n(n+2)}=-i$. Therefore the De Rham cohomology class $\{\alpha\}\in H^n(X,\,\C)$ represented by the $\Delta$-harmonic form $\alpha$ must belong to $H^n_{+}(X,\C)$ when $n$ is even, resp. to $H^n_{-}(X,\C)$ when $n$ is odd.  \hfill $\Box$

\vspace{3ex}

 Let us now consider a holomorphic family $(J_t)_{t\in\Delta}$ of Calabi-Yau $\partial\bar\partial$ complex structures on a compact differential manifold $X$. We set $X_t:=(X,\,J_t)$ and let $n:=\mbox{dim}_{\C}X_t$ for all $t\in\Delta$. Notice that $Q$ and $H$ (cf. (\ref{eqn:Q-def}) and (\ref{eqn:H-def})) depend only on the differential structure of $X$. Thus,

\vspace{1ex}

\hspace{3ex} $C_{+}:=\bigg\{\{\alpha\}\in H^n(X,\,\C)\,\,/\,\, H(\{\alpha\},\,\{\alpha\})>0\bigg\}\subset H^n(X,\,\C),$

\vspace{1ex}

\noindent and

\vspace{1ex}

\hspace{3ex} $C_{-}:=\bigg\{\{\alpha\}\in H^n(X,\,\C)\,\,/\,\, H(\{\alpha\},\,\{\alpha\})<0\bigg\}\subset H^n(X,\,\C)$

\vspace{1ex}

\noindent are {\it open} subsets of $H^n(X,\,\C)$ and depend only on the differential structure of $X$. Furthermore, if we equip the fibres $X_t$ with a $C^{\infty}$ family of arbitrary Hermitian metrics $(\omega_t)_{t\in\Delta}$, the corresponding Hodge $\star$ operator $\star = \star_t$ has eigenspaces $H^n_{+}(X_t,\,\C)$ and $H^n_{-}(X_t,\,\C)$ (cf. (\ref{eqn:Hn-decomposition})) depending on the complex structure $J_t$ via the metric $\omega_t$ (which is in particular a $J_t$-type $(1,\,1)$-form). Lemma \ref{Lem:H-pos} ensures that

\vspace{1ex}

 $H^n_{+}(X_t,\,\C)\setminus\{0\} \subset C_{+} \hspace{3ex} \mbox{and} \hspace{3ex} H^n_{-}(X_t,\,\C)\setminus\{0\} \subset C_{-}  \hspace{3ex} \mbox{for all} \hspace{1ex} t\in\Delta.$

\vspace{1ex}

\noindent Moreover, Lemmas \ref{Lem:H-pos} and \ref{Lem:Hn0-pos} imply the following inclusions\!\!:

\begin{eqnarray}\nonumber H^{n,\,0}(X_t,\,\C)\setminus\{0\}\subset H^n_{+}(X_t,\,\C)\setminus\{0\} \subset C_{+}\subset H^n(X,\,\C) \hspace{2ex} \mbox{if}\,n\,\mbox{is even}, & & \\
\label{eqn:Hn0+-} H^{n,\,0}(X_t,\,\C)\setminus\{0\}\subset H^n_{-}(X_t,\,\C)\setminus\{0\} \subset C_{-}\subset H^n(X,\,\C) \hspace{2ex} \mbox{if}\,n\,\mbox{is odd}. & & \end{eqnarray}

\vspace{2ex}

\noindent It is clear that for any class $\varphi_t=[\alpha_t]\in H^{n,\,0}(X_t,\,\C)$, $Q(\varphi_t,\,\varphi_t)=0$ since $\alpha_t\wedge\alpha_t=0$ for any form of $J_t$-type $(n,\,0)$. Thus the {\bf period domain}, containing the complex lines $H^{n,\,0}(X_t,\,\C)$ varying inside $H^n(X,\,\C)$ when $J_t$ varies, can be defined as in the standard (i.e. K\"ahler) case as

\vspace{2ex}

\noindent $D=\{\C\mbox{-line}\,\, l\subset H^n(X,\,\C)\,\,;\,\,\forall\varphi\in l\setminus\{0\},\,\,Q(\varphi,\,\varphi)=0\,\,\mbox{and}\,\, H(\varphi,\,\varphi)>0\}$

\vspace{1ex}

\noindent if $n$ is even (so, in particular, $l\subset C_{+}$ whenever $l\in D$), and as

\vspace{2ex}

\noindent $D=\{\C\mbox{-line}\,\, l\subset H^n(X,\,\C)\,\,;\,\,\forall\varphi\in l\setminus\{0\},\,\,Q(\varphi,\,\varphi)=0\,\,\mbox{and}\,\, H(\varphi,\,\varphi)<0\}$

\vspace{1ex}

\noindent if $n$ is odd (so, in particular, $l\subset C_{-}$ whenever $l\in D$). Given the natural holomorphic embedding $D\subset\Proj H^n(X,\,\C)$, the complex manifold $D$ is projective and is contained in the quadric defined by $Q$ in $\Proj H^n(X,\,\C)$.

\vspace{3ex}

 We can now show that the {\bf local Torelli theorem} holds in this context.

\begin{The}\label{The:localTorelli} Let $X$ be a compact Calabi-Yau $\partial\bar\partial$-manifold, $\mbox{dim}_{\C}X=n$, and let $\pi:{\cal X}\longrightarrow\Delta$ be its Kuranishi family. Then the associated period map

$${\cal P}\,\,:\,\,\Delta\longrightarrow D\subset\Proj H^n(X,\,\C), \hspace{3ex}\Delta\ni t\mapsto H^{n,\,0}(X_t,\,\C),$$

\noindent is a local holomorphic immersion.

\end{The}

\noindent {\it Proof.} As usual, we denote by $(X_t)_{t\in\Delta}$ the fibres of the Kuranishi family of $X=X_0$. They are all $C^{\infty}$-diffeomorphic to $X$ and the holomorphic family $(X_t)_{t\in\Delta}$ can be seen as a fixed $C^{\infty}$ manifold $X$ equipped with a holomorphic family of complex structures $(J_t)_{t\in\Delta}$. Let $(u_t)_{t\in\Delta}$ be a holomorphic family of nowhere vanishing $n$-forms on $X$ such that for every $t\in\Delta$, $u_t$ is of type $(n,\,0)$ for the complex structure $J_t$ and $\bar\partial_tu_t=0$. The form $u_t$ identifies with the class $[u_t]$ it defines in $H^{n,\,0}(X_t,\,\C)$, hence with the whole space $H^{n,\,0}(X_t,\,\C)=\C\,u_t$. Thus the period map identifies with the map

$$\Delta\ni t\mapsto u_t.$$

 It suffices to prove that ${\cal P}$ is a local immersion at $t=0$. Recall that in the present situation the Kodaira-Spencer map $\rho\,\,:\,\,T_0\Delta\rightarrow H^{0,\,1}(X,\,T^{1,\,0}X)$ is an isomorphism (thanks to Theorem \ref{The:BTT}) and that for any tangent vector $\partial/\partial t\in T_0\Delta$, the choice of a representative $\theta$ in the class $\rho(\partial/\partial t)=[\theta]\in H^{0,\,1}(X,\,T^{1,\,0}X)$ determines a $C^{\infty}$ trivialisation $\Phi\,\,:\,\,{\cal X}\longrightarrow\Delta\times X_0$ (after possibly shrinking $\Delta$ about $0$), which in turn determines about any pre-given point $x\in X$ a choice of local $J_t$-holomorphic coordinates $z_1(t),\dots , z_n(t)$ for every $t\in\Delta$.

 Denote $u=u_0$. Fix an arbitrary tangent vector $\partial/\partial t\in T_0\Delta\setminus\{0\}$ and choose a representative $\theta$ of the class $\rho(\partial/\partial t)\in H^{0,\,1}(X,\,T^{1,\,0}X)$ such that the representative $\theta\lrcorner u$ of the class $[\theta\lrcorner u]\in H^{n-1,\,1}(X,\,\C)$ is $d$-closed. This is possible by the $\partial\bar\partial$ assumption on $X$ and by Lemmas \ref{Lem:d-closed-rep} and \ref{Lem:3space-decomp-contr}. The associated local $C^{\infty}$ trivialisation $\Phi:{\cal X}\rightarrow\Delta\times X_0$ induces $C^{\infty}$ diffeomorphisms $\Phi_t^{-1}: X_0\rightarrow X_t$, $t\in\Delta$, so the differential of the period map at $t=0$ in the $\partial/\partial t$-direction identifies with

\begin{eqnarray}\label{eqn:diff-periodmap}\frac{\partial(\Phi_t^{-1})^{\star}u_t}{\partial t}_{|t=0}=\theta\lrcorner u + v    \hspace{3ex} \mbox{on}\,\,X,\end{eqnarray}

\noindent where $v$ is some $(n,\,0)$-form on $X=X_0$. The identity in (\ref{eqn:diff-periodmap}) can be proved in the usual way (see e.g. [Tia87, proof of Lemma 7.2])\!\!: having fixed an arbitrary point $x\in X$, one writes

\begin{eqnarray}\label{eqn:ut-local}u_t = f_t\,dz_1(t)\wedge\dots\wedge dz_n(t)\end{eqnarray}

\noindent where $f_t$ is a holomorphic function in a neighbourhood of $x$ in $X_t$ and $z_1(t),\dots , z_n(t)$ are the local $J_t$-holomorphic coordinates about $x$ determined by the choice of $\theta$ in the class $\rho(\partial/\partial t)$. Taking $\partial/\partial t$ at $t=0$ in (\ref{eqn:ut-local}), one finds on the right-hand side the sum of the form $v=(\partial f_t/\partial t)_{|t=0} \,dz_1(0)\wedge\dots\wedge dz_n(0)$ of $J_0$-type $(n,\,0)$ with the form $\theta\lrcorner u$ of $J_0$-type $(n-1,\, 1)$. The latter form is easily seen to be the sum of the terms obtained by deriving one of the $dz_j(t)$ in (\ref{eqn:ut-local}) since, with the above choices of $\theta$ and $z_1(t),\dots , z_n(t)$, we have

$$\frac{\partial}{\partial t}(dz_j(t))_{|t=0} = \theta\lrcorner\, dz_j(0),  \hspace{3ex} j=1, \dots , n.$$

\noindent Now, $du_t=0$ for all $t$, hence the left-hand term in (\ref{eqn:diff-periodmap}) is a $d$-closed $n$-form on $X$. Thus $d(\theta\lrcorner u + v)=0$. By our choice of $\theta$ (based on a key application of the $\partial\bar\partial$ lemma), $d(\theta\lrcorner u)=0$, hence $dv=0$. In particular, $v$ is a $\bar\partial_0$-closed form of $J_0$-type $(n,\,0)$, so $v=c\,u$ for some constant $c\in\C$.

 It is now clear that if $(d{\cal P})_0(\partial/\partial t)=0$, then $\theta\lrcorner u=0$ and $v=c\,u=0$, so $\theta=0$ (since $T_u(\theta)=\theta\lrcorner u$ and $T_u$ is an isomorphism -- see (\ref{eqn:u-isom})), hence $\partial/\partial t=0$ (since the Kodaira-Spencer map is an isomorphism here). This last vanishing contradicts the choice of $\partial/\partial t\neq 0$. We have thus shown that ${\cal P}$ is a local immersion at $t=0$.   \hfill $\Box$

\subsection{Weil-Petersson metrics on $\Delta$}\label{subsection:WP-metrics}

 We start with a refinement of Lemma \ref{Lem:d-closed-rep} singling out a particular $d$-closed representative of a given Dolbeault cohomology class on a $\partial\bar\partial$-manifold.

\begin{Def}\label{def:min-d-closed} Let $X$ be a compact $\partial\bar\partial$-manifold equipped with an arbitrary Hermitian metric $\omega$. Given any Dolbeault cohomology class $[\alpha]\in H^{p,\,q}(X,\,\C)$, let $\alpha$ be its $\Delta''_{\omega}$-harmonic representative and let $v_{min}\in\mbox{Im}(\partial\bar\partial)^{\star}\subset C^{\infty}_{p,\,q-1}(X,\,\C)$ be the solution of minimal $L^2$ norm (w.r.t. $\omega$) of equation (\ref{eqn:del-prooflemma}). 

The $d$-closed $(p,\,q)$-form $\alpha_{min}:=\alpha + \bar\partial v_{min}$ will be called the {\bf $\omega$-minimal  $d$-closed representative} of the class $[\alpha]$. (It coincides with the $\Delta''_{\omega}$-harmonic representative if $\omega$ is K\"ahler.)

\end{Def}

 A word of explanation is in order. Recall that the Aeppli cohomology group of type $(p,\,q)$ is standardly defined as

$$H^{p, \, q}_A(X, \, \C)=\frac{\ker(\partial\bar\partial:C^{\infty}_{p, \, q}(X)\rightarrow C^{\infty}_{p+1, \, q+1}(X))}{\mbox{Im}(\partial:C^{\infty}_{p-1, \, q}(X)\rightarrow C^{\infty}_{p, \, q}(X)) + \mbox{Im}(\bar\partial:C^{\infty}_{p, \, q-1}(X)\rightarrow C^{\infty}_{p, \, q}(X))}$$

\noindent and that the fourth-order Aeppli Laplacian $\tilde\Delta^{p, \, q}_A:C^{\infty}_{p, \, q}(X, \, \C)\rightarrow C^{\infty}_{p, \, q}(X, \, \C)$ (cf. [KS60], also [Sch07, 2.c., p. 9-10]) defined by

$$\tilde\Delta^{p, \, q}_A:=\partial\partial^{\star} + \bar\partial\bar\partial^{\star} + (\partial\bar\partial)^{\star}(\partial\bar\partial) + (\partial\bar\partial)(\partial\bar\partial)^{\star} + \partial\bar\partial^{\star}(\partial\bar\partial^{\star})^{\star} + (\partial\bar\partial^{\star})^{\star}\partial\bar\partial^{\star}$$

\noindent is elliptic and thus induces a three-space decomposition

$$C^{\infty}_{p, \, q}(X, \C)=\ker\tilde\Delta^{p, \, q}_A \oplus (\mbox{Im}\partial + \mbox{Im}\bar\partial) \oplus \mbox{Im}(\partial\bar\partial)^{\star}$$

\noindent that is orthogonal w.r.t. the $L^2$ scalar product defined by $\omega$ and in which 

\begin{equation}\label{eqn:A-kernel}\ker(\partial\bar\partial)=\ker\tilde\Delta^{p, \, q}_A \oplus (\mbox{Im}\partial + \mbox{Im}\bar\partial),\end{equation}

\noindent yielding the Hodge isomorphism $H^{p, \, q}_A(X, \, \C)\simeq \ker\tilde\Delta^{p, \, q}_A$. Since the solution $v$ of equation (\ref{eqn:del-prooflemma}) is unique only modulo $\ker(\partial\bar\partial)$, the solution of minimal $L^2$ norm is the unique solution lying in $\ker(\partial\bar\partial)^{\perp}=\mbox{Im}(\partial\bar\partial)^{\star}$. Note that if the $\Delta''$-harmonic representative $\alpha$ of the class $[\alpha]$ happens to be $d$-closed (for example, this is the case if the metric $\omega$ is K\"ahler), then $\partial\alpha=0$ and $v_{min}=0$, so $\alpha_{min}=\alpha$. Thus $\alpha_{min}$ can be seen as the minimal $d$-closed correction in a given Dolbeault class of the $\Delta''$-harmonic representative of that class.

\vspace{2ex}

 Recall that if we fix a compact balanced Calabi-Yau $\partial\bar\partial$-manifold $(X,\,\omega)$ ($\mbox{dim}_{\C}X=n$), the base space $\Delta_{[\omega^{n-1}]}$ of the local universal family $(X_t)_{t\in\Delta_{[\omega^{n-1}]}}$ of deformations of $X$ that are co-polarised by the balanced class $[\omega^{n-1}]\in H^{n-1,\,n-1}(X,\,\C)$ identifies to an open subset of $H^{0,\,1}(X,\,T^{1,\,0}X)_{[\omega^{n-1}]}$ and

\vspace{2ex}

$T_t\Delta_{[\omega^{n-1}]}\simeq H^{0,\,1}(X_t,\,T^{1,\,0}X_t)_{[\omega^{n-1}]}\simeq H^{n-1,\,1}_{prim}(X_t,\,\C), \hspace{3ex} t\in\Delta_{[\omega^{n-1}]}.$

\vspace{2ex}

\noindent We shall now define two Weil-Petersson metrics on $\Delta_{[\omega^{n-1}]}$ induced by pre-given balanced metrics on the fibres $X_t$ whose $(n-1)^{st}$ powers lie in the co-polarising balanced class.

\begin{Def}\label{Def:WP1} Fix any holomorphic family of nonvanishing holomorphic $n$-forms $(u_t)_{t\in\Delta}$ on the fibres $(X_t)_{t\in\Delta}$. Let $(\omega_t)_{t\in\Delta_{[\omega^{n-1}]}}$ be a $C^{\infty}$ family of balanced metrics on the fibres $(X_t)_{t\in\Delta_{[\omega^{n-1}]}}$ such that $\omega_t^{n-1}\in\{\omega^{n-1}\}$ for all $t$ and $\omega_0=\omega$. The associated {\bf Weil-Petersson metrics} $G^{(1)}_{WP}$ and $G^{(2)}_{WP}$ on $\Delta_{[\omega^{n-1}]}$ are defined as follows. For any $t\in\Delta_{[\omega^{n-1}]}$ and any $[\theta_t], [\eta_t]\in H^{0,\,1}(X_t,\,T^{1,\,0}X_t)_{[\omega^{n-1}]}$, let

\begin{eqnarray}\label{eqn:WP1}G^{(1)}_{WP}([\theta_t],\,[\eta_t]) & := & \frac{\langle\langle\theta_t,\, \eta_t\rangle\rangle}{\int\limits_{X_t}dV_{\omega_t}},\hspace{3ex} (\mbox{where}\,\,\,dV_{\omega_t}\!:=\frac{\omega_t^n}{n!})\\
\label{eqn:WP2}G^{(2)}_{WP}([\theta_t],\,[\eta_t]) & := & \frac{\langle\langle\theta_t\lrcorner u_t,\, \eta_t\lrcorner u_t\rangle\rangle}{i^{n^2}\,\int_{X_t}u_t\wedge\bar{u}_t},\end{eqnarray}

\noindent where $\theta_t$ (resp. $\eta_t$) is chosen in its class $[\theta_t]$ (resp. $[\eta_t]$) such that $\theta_t\lrcorner u_t$ (resp. $\eta_t\lrcorner u_t$) is the $\omega_t$-minimal $d$-closed representative of the class $[\theta_t\lrcorner u_t]\in H^{n-1,\,1}(X_t,\,\C)$ (resp. $[\eta_t\lrcorner u_t]\in H^{n-1,\,1}(X_t,\,\C)$), while $\langle\langle\,\,,\,\,\rangle\rangle$ stands for the $L^2$ scalar product induced by $\omega_t$ on the spaces involved. 

 The $C^{\infty}$ positive definite $(1,\, 1)$-forms on $\Delta_{[\omega^{n-1}]}$ associated with $G^{(1)}_{WP}$ and $G^{(2)}_{WP}$ are denoted by

\vspace{1ex}

\hspace{10ex} $\omega^{(1)}_{WP}>0 \hspace{2ex} \mbox{and} \hspace{2ex} \omega^{(2)}_{WP}>0 \hspace{2ex} \mbox{on} \hspace{1ex}\Delta_{[\omega^{n-1}]}.$

\end{Def}  

 Since every $u_t$ is unique up to a constant factor, the definition of $G^{(2)}_{WP}$ is independent of the choice of the family $(u_t)_{t\in\Delta}$. From Lemma \ref{Lem:u-parallel} we infer

\begin{Obs}\label{Obs:WP-metrics-coincide} If the balanced metrics can be chosen such that $\mbox{Ric}(\omega_t)=0$ for all $t\in\Delta_{[\omega^{n-1}]}$, then

\vspace{1ex}

\hspace{20ex} $\omega^{(1)}_{WP} = \omega^{(2)}_{WP} \hspace{2ex}\mbox{on}\hspace{1ex}\Delta_{[\omega^{n-1}]}.$

\end{Obs}

\subsection{Metric on $\Delta$ induced by the period map}\label{subsection:PM-metric}

 Let $L={\cal O}_{\Proj H^n(X,\,\C)}(-1)$ be the tautological line bundle on $\Proj H^n(X,\,\C)$. We will endow the restrictions of $L$ to two open subsets of $\Proj H^n(X,\,\C)$ with Hermitian fibres metrics induced by $H$. We set\!\!:

\vspace{1ex}

\noindent  $U^n_{+}:=\bigg\{[l]\in\Proj H^n(X,\,\C)\,\,/\,\,l \hspace{1ex} \mbox{is a}\,\C\mbox{-line such that}\hspace{1ex} l\subset C_{+}\bigg\}\subset\Proj H^n(X,\,\C),$

\vspace{1ex}

\noindent and

\vspace{1ex}

 \noindent $U^n_{-}:=\bigg\{[l]\in\Proj H^n(X,\,\C)\,\,/\,\,l \hspace{1ex} \mbox{is a}\,\C\mbox{-line such that}\hspace{1ex} l\subset C_{-}\bigg\}\subset\Proj H^n(X,\,\C),$

\vspace{1ex}

\noindent where $[l]$ denotes the point in $\Proj H^n(X,\,\C)$ defined by the line $l\subset H^n(X,\,\C)$. It follows from the discussion of $C_{+}$ and $C_{-}$ in $\S.$\ref{subsection:period-domain} that $U^n_{+}$ and $U^n_{-}$ are {\it open} subsets of $\Proj H^n(X,\,\C)$ and depend only on the differential structure of $X$.
 
 Moreover, for every $[l]\in U^n_{+}$, the fibre $L_{[l]} = l\subset C_{+}$ is endowed with the scalar product defined by the restriction of $H$. Thus $L_{|U^n_{+}}$ has a Hermitian fibre metric $h^{+}_L$ induced by $H$. The (negative) curvature form $i\,\Theta_{h^{+}_L}(L_{|U^n_{+}})$ defines the associated Fubini-Study metric on $U^n_{+}$ by

$$\omega^{+}_{FS}=-i\,\Theta_{h^{+}_L}(L_{|U^n_{+}})>0 \hspace{3ex} \mbox{on}\hspace{1ex} U^n_{+}\subset\Proj H^n(X,\,\C).$$

\noindent Likewise, for every $[l]\in U^n_{-}$, the fibre $L_{[l]} = l\subset C_{-}$ is endowed with the scalar product defined by the restriction of $-H$. Thus $L_{|U^n_{-}}$ has a Hermitian fibre metric $h^{-}_L$ induced by $-H$. The (negative) curvature form $i\,\Theta_{h^{-}_L}(L_{|U^n_{-}})$ defines the associated Fubini-Study metric on $U^n_{-}$ by

$$\omega^{-}_{FS}=-i\,\Theta_{h^{-}_L}(L_{|U^n_{-}})>0 \hspace{3ex} \mbox{on}\hspace{1ex} U^n_{-}\subset\Proj H^n(X,\,\C).$$

\noindent It follows from the above discussion that $\omega^{+}_{FS}$ and $\omega^{-}_{FS}$ depend only on the differential structure of $X$. Composing the period map with the holomorphic embedding $D\stackrel{\iota}{\hookrightarrow}\Proj H^n(X,\,\C)$, we obtain a local holomorphic immersion $\iota\circ{\cal P}:\Delta\rightarrow\Proj H^n(X,\,\C)$ (cf. Theorem \ref{The:localTorelli}). From (\ref{eqn:Hn0+-}), we get\!\!:

\vspace{2ex}

\hspace{6ex} $\mbox{Im}\,(\iota\circ{\cal P})\subset U^n_{+} \hspace{2ex}\mbox{if}\,\,n\,\,\mbox{is even}, \hspace{2ex} \mbox{Im}\,(\iota\circ{\cal P})\subset U^n_{-} \hspace{2ex}\mbox{if}\,\,n\,\,\mbox{is odd}.$

\vspace{2ex}

\noindent Taking the inverse image of $\omega^{+}_{FS}$ when $n$ is even, resp. of $\omega^{-}_{FS}$ when $n$ is odd, we get a Hermitian metric (i.e. a positive definite $C^{\infty}$ $(1,\,1)$-form) $\gamma$ on $\Delta$ which is actually K\"ahler\!\!:

\vspace{2ex}

$\gamma:=(\iota\circ{\cal P})^{\star}(\omega^{+}_{FS})>0  \hspace{2ex}\mbox{if}\,\,n\,\,\mbox{is even}, \hspace{2ex} \gamma:=(\iota\circ{\cal P})^{\star}(\omega^{-}_{FS})>0 \hspace{2ex}\mbox{if}\,\,n\,\,\mbox{is odd}.$ 

\vspace{2ex} 

\noindent {\bf Computation of $\gamma$.} We can compute $\gamma$ at any point $t\in\Delta$ (e.g. at $t=0$) in the same way as in [Tia87, $\S.7$]. We spell out the details for the reader's convenience.  Let $(u_t)_{t\in\Delta}$ be a holomorphic family of nonvanishing holomorphic $n$-forms on the fibres $(X_t)_{t\in\Delta}$. Recall that a tangent vector $(\partial/\partial t)_{|t=0}$ to $\Delta$ at $0$ identifies via the Kodaira-Spencer map with a class $[\theta]\in H^{0,\,1}(X,\,T^{1,\, 0}X)$. Fix any such class $[\theta]$. We will compute $\gamma_0([\theta],\,[\theta])$.

 We have\!: $L_{u_t}=\C\cdot u_t=H^{n,\,0}(X_t,\,\C)$. Thus\!:

\vspace{1ex}

\noindent $(i)$\,\,if $n$ is {\it even}, then $L_{u_t}\subset H^n_{+}(X_t,\,\C)$ and $(-i)^n\,Q(u_t,\,\bar{u}_t)=H(u_t,\,u_t)=|u_t|_{h^{+}_L}^2=e^{-\rho(t)}$, where $\rho$ denotes the local weight function of the fibre metric $h^{+}_L$ of $L_{|U^n_{+}}$. We get $\rho(t)=-\log((-i)^n\,Q(u_t,\,\bar{u}_t));$

\vspace{1ex}

\noindent $(ii)$\,\,if $n$ is {\it odd}, then $L_{u_t}\subset H^n_{-}(X_t,\,\C)$ and $-(-i)^n\,Q(u_t,\,\bar{u}_t) = -H(u_t,\,u_t)=|u_t|_{h^{-}_L}^2=e^{-\rho(t)}$, where $\rho$ denotes the local weight function of the fibre metric $h^{-}_L$ of $L_{|U^n_{-}}$. We get $\rho(t)=-\log(-(-i)^n\,Q(u_t,\,\bar{u}_t))$.

\vspace{1ex}

 Now suppose that $n$ is {\it even}. The curvature form of $(L,\,h^{+}_L)$ on a $\C$-line $\C\!\cdot\! t$ in a small neighbourhood of $0$ equals $i\partial_t\bar\partial_t\rho(t)$, which in turn equals\!\!:

$$-i\partial_t\bar\partial_t\log ((-i)^n\,Q(u_t,\,\bar{u}_t)) = -i\,\frac{\partial^2\log ((-i)^n\,Q(u_t,\,\bar{u}_t))}{\partial t\,\partial\bar{t}}\,dt\wedge d\bar{t},$$

\noindent This means that for $[\theta]=\rho(\partial/\partial t_{|t=0})$, using the fact that $\frac{\partial u_t}{\partial\bar{t}}=0$ (since $u_t$ varies holomorphically with $t$), we get\!\!:

\begin{eqnarray}\nonumber \gamma_0([\theta],\,[\theta]) & = & -\frac{\partial^2\log ((-i)^n\,Q(u_t,\,\bar{u}_t))}{\partial t\,\partial\bar{t}}_{|t=0} = -\frac{\partial}{\partial t}\bigg((-1)^n\,\frac{Q(u_t,\,\frac{\partial\bar{u}_t}{\partial\bar{t}})}{Q(u_t,\,\bar{u}_t)}\bigg)_{|t=0} \\
\nonumber & = & (-1)^{n+1}\,\bigg[\frac{Q(\frac{\partial u_t}{\partial t}_{|t=0},\,\frac{\partial\bar{u}_t}{\partial\bar{t}}_{|t=0})}{Q(u_0,\,\bar{u}_0)} - \frac{Q(\frac{\partial u_t}{\partial t}_{|t=0},\,\bar{u}_0)\cdot Q(u_0,\, \frac{\partial\bar{u}_t}{\partial\bar{t}}_{|t=0})}{Q(u_0,\,\bar{u}_0)^2}\bigg].\end{eqnarray}

\noindent Now recall that in the proof of Theorem \ref{The:localTorelli} a key application of the $\partial\bar\partial$ lemma enabled us to choose the representative $\theta$ of the class $[\theta]$ such that $d(\theta\lrcorner u)=0$. With this choice, if $u\!\!:=u_0$, in formula (\ref{eqn:diff-periodmap}) we had $v=c\,u$ and

\vspace{1ex}

\hspace{20ex} $\displaystyle \frac{\partial u_t}{\partial t}_{|t=0}=\theta\lrcorner u + c\,u,$

\vspace{1ex} 

\noindent where $c\in\C$ is a constant, if we identify $u_t$ with $(\Phi_t^{-1})^{\star}u_t$ when $\Phi_t:X_t\rightarrow X_0$ ($t\in\Delta$) denote the $C^{\infty}$ isomorphisms induced by the choice of $\theta$ in $[\theta]$. Using this, the above formula for $\gamma_0([\theta],\,[\theta])$ translates to

\begin{eqnarray}\label{eqn:gamma0-2}\nonumber\gamma_0([\theta],\,[\theta]) & = & (-1)^{n+1}\,\frac{Q(u,\,\bar{u})\cdot Q(\theta\lrcorner u,\, \overline{\theta\lrcorner u}) + |c|^2\,Q(u,\,\bar{u})^2 - |c|^2\,Q(u,\,\bar{u})^2}{Q(u,\,\bar{u})^2}\\
\nonumber  & = & (-1)^{n+1}\,\frac{Q(\theta\lrcorner u,\, \overline{\theta\lrcorner u})}{Q(u,\,\bar{u})} = \frac{-H(\{\theta\lrcorner u\},\,\{\theta\lrcorner u\})}{i^{n^2}\,\int\limits_Xu\wedge\bar{u}}.\end{eqnarray}

 In the case when $n$ is odd, the formula for $\gamma_0([\theta],\,[\theta])$ gets an extra $(-1)$ factor. The conclusion of these calculations is summed up in the following

\begin{Lem}\label{Lem:concl-calc-gamma} The K\"ahler metric $\gamma$ defined on $\Delta$ by $\gamma:=(\iota\circ{\cal P})^{\star}(\omega^{+}_{FS})>0$ when $n$ is even and by $\gamma:=(\iota\circ{\cal P})^{\star}(\omega^{-}_{FS})>0$ when $n$ is odd, is independent of the choice of any metrics on $(X_t)_{t\in\Delta}$ and is explicitly given by the formula\!\!:

$$\gamma_t([\theta_t],\,[\theta_t])=\frac{-\int\limits_X(\theta_t\lrcorner u_t)\wedge\overline{(\theta_t\lrcorner u_t)}}{i^{n^2}\,\int\limits_Xu_t\wedge\bar{u_t}} = \frac{-H(\{\theta_t\lrcorner u_t\},\,\{\theta_t\lrcorner u_t\})}{i^{n^2}\,\int\limits_Xu_t\wedge\bar{u_t}}, \hspace{2ex} \mbox{if}\,\,n\,\,\mbox{is even},$$

$$\gamma_t([\theta_t],\,[\theta_t])=\frac{-i\,\int\limits_X(\theta_t\lrcorner u_t)\wedge\overline{(\theta_t\lrcorner u_t)}}{i^{n^2}\,\int\limits_Xu_t\wedge\bar{u_t}} = \frac{H(\{\theta_t\lrcorner u_t\},\,\{\theta_t\lrcorner u_t\})}{i^{n^2}\,\int\limits_Xu_t\wedge\bar{u_t}}, \hspace{2ex} \mbox{if}\,\,n\,\,\mbox{is odd},$$

\noindent for every $t\in\Delta$ and every $[\theta_t]\in H^{0,\,1}(X_t,\,T^{1,\,0}X_t)$.

\noindent In particular, we see that $\gamma_t([\theta_t],\,[\theta_t])$ is independent of the choice of representative $\theta_t$ in the class $[\theta_t]\in H^{0,\,1}(X_t,\,T^{1,\,0}X_t)$ such that $\theta_t\lrcorner u_t$ is $d$-closed. Since for every $t\in\Delta$, $u_t$ is unique up to a constant factor, $\gamma$ is independent of the choice of holomorphic family $(u_t)_{t\in\Delta}$ of $J_t$-holomorphic $n$-forms.

\end{Lem}

 Notice that $i^{n^2}\,u_t\wedge\bar{u_t}>0$ at every point of $X_t$ for any non-vanishing $(n,\,0)$-form $u_t$. On the other hand, it follows from  Lemma \ref{Lem:n1-1decompositions} below that $H(\{\theta_t\lrcorner u_t\},\,\{\theta_t\lrcorner u_t\})<0$ when $n$ is {\it even} and that  $H(\{\theta_t\lrcorner u_t\},\,\{\theta_t\lrcorner u_t\})>0$ when $n$ is {\it odd} if a $d$-closed representative $\theta_t\lrcorner u_t$ of the class $[\theta_t\lrcorner u_t]\in H^{n-1,\,1}(X_t,\,\C)\subset H^n(X,\,\C)$ can be chosen to be {\it primitive}. This reproves that $\gamma_t([\theta_t],\,[\theta_t])>0$ in this case (which does occur if primitivess is taken w.r.t. a K\"ahler metric).

 \subsection{Comparison of metrics on $\Delta$}\label{subsection:comparison-metrics} 

We shall now compare the Weil-Petersson metric $\omega^{(2)}_{WP}$ with the period-map metric $\gamma$ on $\Delta_{[\omega^{n-1}]}$. We need a general fact first.

\vspace{2ex}

 Let $X$ be a compact complex manifold ($\mbox{dim}_{\C}X=n$) equipped with a Hermitian metric $\omega$ and let $\star:\Lambda^{n-1,\, 1}\rightarrow\Lambda^{n-1,\, 1}$ be the Hodge $\star$ operator defined by $\omega$ on $(n-1,\, 1)$-forms. (Here $\Lambda^{n-1,\, 1}$ stands for the space $C^{\infty}_{n-1,\,1}(X,\,\C)$ of global smooth forms of bidegree $(n-1,\,1)$ on $X$ although $\star$ acts even pointwise on forms.) Since $\star^2=(-1)^n$, $\star$ induces a decomposition that is orthogonal for the $L^2$ scalar product defined by $\omega$ on $X$ (cf. $\S.$\ref{subsection:period-domain})\!\!:

\begin{eqnarray}\label{eqn:duality-decomposition} \Lambda^{n-1,\, 1}=\Lambda^{n-1,\, 1}_{-}\oplus\Lambda^{n-1,\, 1}_{+}  \hspace{3ex}\mbox{(the duality decomposition)},\end{eqnarray}

\noindent where $\Lambda^{n-1,\, 1}_{\pm}$ stand for the eigenspaces of $\star$ corresponding to the eigenvalues $\pm 1 $ (if $n$ is even), $\pm i $ (if $n$ is odd). On the other hand, the Hermitian metric $\omega$ induces the Lefschetz decomposition (cf. [Voi02, proposition 6.22, p. 147]) 

\begin{eqnarray}\label{eqn:Lefschetz-decomposition} \Lambda^{n-1,\, 1}=\Lambda^{n-1,\, 1}_{prim}\oplus\bigg(\omega\wedge\Lambda^{n-2,\, 0}\bigg),\end{eqnarray}

\noindent which is again orthogonal for the $L^2$ scalar product defined by $\omega$ on $X$, where $\Lambda^{n-1,\, 1}_{prim}$ denotes the space of {\it primitive} $(n-1,\, 1)$-forms $u$ (i.e. those $u\in\Lambda^{n-1,\, 1}$ for which $\omega\wedge u=0$ or, equivalently, $\Lambda u=0$), while $\omega\wedge\Lambda^{n-2,\, 0}$ denotes the space of forms $\omega\wedge v$ with $v$ an arbitrary form of bidegree $(n-2,\, 0)$.

\begin{Lem}\label{Lem:n1-1decompositions} The decompositions (\ref{eqn:duality-decomposition}) and (\ref{eqn:Lefschetz-decomposition}) coincide up to order, i.e. \\

\hspace{3ex} $\Lambda^{n-1,\, 1}_{-}=\Lambda^{n-1,\, 1}_{prim} \hspace{3ex} \mbox{and} \hspace{3ex} \Lambda^{n-1,\, 1}_{+}=\omega\wedge\Lambda^{n-2,\, 0} \hspace{3ex} \mbox{if}\,\,n\,\,\mbox{is even},$

\vspace{1ex}

\hspace{3ex} $\Lambda^{n-1,\, 1}_{+}=\Lambda^{n-1,\, 1}_{prim} \hspace{3ex} \mbox{and} \hspace{3ex} \Lambda^{n-1,\, 1}_{-}=\omega\wedge\Lambda^{n-2,\, 0} \hspace{3ex} \mbox{if}\,\,n\,\,\mbox{is odd}.$

\end{Lem}

\noindent {\it Proof.} It suffices to prove the inclusions\!\!: 

\vspace{1ex}

\hspace{3ex} $(A)\hspace{2ex}\Lambda^{n-1,\, 1}_{prim}\subset\Lambda^{n-1,\, 1}_{-}  \hspace{3ex} \mbox{and} \hspace{3ex} (B)\hspace{2ex}\omega\wedge\Lambda^{n-2,\, 0}\subset\Lambda^{n-1,\, 1}_{+}  \hspace{3ex} \mbox{if}\,\,n\,\,\mbox{is even},$

\vspace{1ex}

\hspace{3ex} $(A)\hspace{2ex}\Lambda^{n-1,\, 1}_{prim}\subset\Lambda^{n-1,\, 1}_{+}  \hspace{3ex} \mbox{and} \hspace{3ex} (B)\hspace{2ex}\omega\wedge\Lambda^{n-2,\, 0}\subset\Lambda^{n-1,\, 1}_{-}  \hspace{3ex} \mbox{if}\,\,n\,\,\mbox{is odd}.$ \\

 Let $u\in\Lambda^{n-1,\, 1}_{prim}$. Formula (\ref{eqn:prim-form-star-formula-gen}) gives $\star u=(-1)^{n(n+1)/2}\,i^{n-2}\,u=i^{n^2+2n-2}\,u$. If $n$ is even, $n^2+2n-2\in 4\Z-2$, hence $i^{n^2+2n-2}=i^{-2}=-1$, so $u\in\Lambda^{n-1,\, 1}_{-}$. If $n$ is odd, $n^2+2n-2\in 4\Z+1$, hence $i^{n^2+2n-2}=i$, so $u\in\Lambda^{n-1,\, 1}_{+}$. This proves inclusions $(A)$.

 To prove inclusions $(B)$, we first prove the following formula

\begin{eqnarray}\label{eqn:star-omega-v}\star(\omega\wedge v)=i^{n(n-2)}\,\omega\wedge v  \hspace{3ex} \mbox{for all}\,\,v\in\Lambda^{n-2,\, 0}.\end{eqnarray}

\noindent Pick any $v\in\Lambda^{n-2,\, 0}$. Then $\omega\wedge v\in\Lambda^{n-1,\, 1}$. For every $u\in\Lambda^{n-1,\,1}$, we have

\begin{eqnarray}\label{eqn:ustar-omega-v1}\int\limits_Xu\wedge\star\overline{(\omega\wedge v)}=\int\limits_X\langle u,\,\omega\wedge v\rangle\,dV_{\omega} = \langle\langle u,\,\omega\wedge v\rangle\rangle = \langle\langle \Lambda u,\, v\rangle\rangle.\end{eqnarray} 

\noindent On the other hand, the following formula holds

\begin{eqnarray}\label{eqn:omega-u-lambda}\omega\wedge u = \frac{\omega^2}{2!}\wedge\Lambda u   \hspace{3ex} \mbox{for all}\,\, u\in\Lambda^{n-1,\,1}.\end{eqnarray}

\noindent Indeed, $\omega^2\wedge\Lambda u = [L^2,\,\Lambda]\,u = 2(n-n+2-1)\,Lu=2\,\omega\wedge u,$ where for the first identity we have used the fact that $L^2u=0$ since $L^2u$ is of type $(n+1,\, 3)$, while for the second identity we have used the standard formula (\ref{eqn:LrOmega}) with $r=2$ and $k=n$. 

\noindent Applying (\ref{eqn:omega-u-lambda}) on the top line below, for every $u\in\Lambda^{n-1,\,1}$ we get

\begin{eqnarray}\label{eqn:ustar-omega-v2}\nonumber\int\limits_Xu\wedge\overline{(\omega\wedge v)} & = &\int\limits_X(\omega\wedge u)\wedge\bar{v}=\int\limits_X\bigg(\frac{\omega^2}{2!}\wedge\Lambda u\bigg)\wedge\bar{v} \\
 \nonumber & = & \int\limits_X(\Lambda u)\wedge\overline{\bigg(\frac{\omega^2}{2!}\wedge v\bigg)} = i^{n(n-2)}\,\int\limits_X(\Lambda u)\wedge\star\bar{v}\\ 
  & = & i^{n(n-2)}\,\int\limits_X\langle\Lambda u,\,v\rangle\,dV_{\omega} = i^{n(n-2)}\,\langle\langle\Lambda u,\,v\rangle\rangle,\end{eqnarray}

\noindent where the last identity on the second line above has followed from the formula

\vspace{1ex}

$\star v = i^{n(n-2)}\,\frac{\omega^2}{2!}\wedge v,  \hspace{3ex} v\in\Lambda^{n-2,\, 0} \hspace{3ex} (\mbox{cf.}\,\,(\ref{eqn:prim-form-star-formula-gen})\,\,\mbox{with}\,\,(p,\,q)=(n-2,\,0)).$

\vspace{1ex}

\noindent It is clear that the combination of (\ref{eqn:ustar-omega-v1}) and (\ref{eqn:ustar-omega-v2}) proves formula (\ref{eqn:star-omega-v}).

 With (\ref{eqn:star-omega-v}) in place, inclusions $(B)$ follow immediately. Indeed, if $n$ is even, $n(n-2)\in 4\Z$, hence $i^{n(n-2)}=1$, so $\omega\wedge v\in\Lambda^{n-1,\, 1}_{+}$ for all $v\in\Lambda^{n-2,\, 0}$. If $n$ is odd, $n(n-2)\in 4\Z -1$, hence $i^{n(n-2)}=-i$, so $\omega\wedge v\in\Lambda^{n-1,\, 1}_{-}$ for all $v\in\Lambda^{n-2,\, 0}$.   \hfill $\Box$

\vspace{3ex}

  For any $\theta\in C^{\infty}_{0,\, 1}(X,\,T^{1,\, 0}X)$, we denote by

\begin{eqnarray}\label{eqn:theta-L-decomp}\theta\lrcorner u = \theta'\lrcorner u + \omega\wedge\zeta\end{eqnarray}

\noindent the decomposition of $\theta\lrcorner u\in\Lambda^{n-1,\, 1}$ induced by the Lefschetz decomposition (\ref{eqn:Lefschetz-decomposition}). Thus $\theta'\lrcorner u\in\Lambda^{n-1,\, 1}_{prim}$ and $\zeta\in\Lambda^{n-2,\, 0}$. By orthogonality we have $||\theta\lrcorner u||^2 = ||\theta'\lrcorner u||^2 + ||\omega\wedge\zeta||^2$. Now

$$||\omega\wedge\zeta||^2=\langle\langle\Lambda(\omega\wedge\zeta),\,\zeta\rangle\rangle = \langle\langle[\Lambda,\, L]\,\zeta,\,\zeta\rangle\rangle = 2||\zeta||^2,$$

\noindent since $\Lambda\zeta=0$ for bidegree reasons (hence $[\Lambda,\, L]\,\zeta = \Lambda(\omega\wedge\zeta) - \omega\wedge\Lambda\zeta = \Lambda(\omega\wedge\zeta)$) and $[\Lambda,\, L]\,\zeta = 2\zeta$ (by formula (\ref{eqn:LrOmega}) with $r=1$ and $k=n-2$).

\begin{The}\label{The:final-metric-formulae} Let $X$ be a compact balanced Calabi-Yau $\partial\bar\partial$-manifold of complex dimension $n$. Then the metrics $G^{(2)}_{WP}$ and $\gamma$ on the base space $\Delta_{[\omega^{n-1}]}$ of the local universal family of deformations of $X$ that are co-polarised by a given balanced class $[\omega^{n-1}]\in H^{n-1,\, n-1}(X,\,\C)\subset H^{2n-2}(X,\,\C)$ are given at every point $t\in\Delta_{[\omega^{n-1}]}$ by the formulae (see notation (\ref{eqn:theta-L-decomp}))\!\!:

\noindent \begin{eqnarray}\label{eqn:WP2-final}G^{(2)}_{WP,\, t}([\theta_t],\,[\theta_t]) = \frac{||\theta'_t\lrcorner u_t||^2 + 2||\zeta_t||^2}{i^{n^2}\,\int_Xu_t\wedge\bar{u_t}}, \hspace{2ex}  [\theta_t]\in H^{0,\,1}(X_t,\,T^{1,\,0}X_t)_{[\omega^{n-1}]}, & & \\ 
\label{eqn:gamma-final}\gamma_t([\theta_t],\,[\theta_t]) = \frac{||\theta'_t\lrcorner u_t||^2 - 2||\zeta_t||^2}{i^{n^2}\,\int_Xu_t\wedge\bar{u_t}}, \hspace{2ex} [\theta_t]\in H^{0,\,1}(X_t,\,T^{1,\,0}X_t)_{[\omega^{n-1}]}. & & \end{eqnarray}

\noindent Here $\theta_t$ is chosen in its class $[\theta_t]$ such that $\theta_t\lrcorner u_t$ is the $\omega_t$-minimal $d$-closed representative of the class $[\theta_t\lrcorner u_t]\in H^{n-1,\,1}(X_t,\,\C)$ (where the $\omega_t\in\{\omega^{n-1}\}$ are balanced metrics in the co-polarising balanced class given beforehand).

\end{The}

\noindent {\it Proof.} We may assume that $t=0$. Formula (\ref{eqn:WP2-final}) follows immediately from (\ref{eqn:WP2}) and from the above considerations. To get (\ref{eqn:gamma-final}), notice that Lemma \ref{Lem:n1-1decompositions} shows that if $n$ is even, then $\theta\lrcorner u=\star(-\theta'\lrcorner u + \omega\wedge\zeta)$, from which we get

$$\int\limits_X(\theta\lrcorner u)\wedge\overline{(\theta\lrcorner u)} = \int\limits_X\bigg(\theta'\lrcorner u + \omega\wedge\zeta\bigg)\wedge\bigg(-\star\overline{(\theta'\lrcorner u)} + \star\overline{(\omega\wedge\zeta)}\bigg) = -||\theta'\lrcorner u||^2 + 2\,||\zeta||^2,$$

\noindent while if $n$ is odd, then $\theta\lrcorner u=\star(-i\,\theta'\lrcorner u +i\, \omega\wedge\zeta)$, from which we get

$$\int\limits_X(\theta\lrcorner u)\wedge\overline{(\theta\lrcorner u)} = \int\limits_X\bigg(\theta'\lrcorner u + \omega\wedge\zeta\bigg)\wedge\bigg(i\star\overline{(\theta'\lrcorner u)} - i\star\overline{(\omega\wedge\zeta)}\bigg) = i\,||\theta'\lrcorner u||^2 - 2i\,||\zeta||^2.$$ 

\noindent Now (\ref{eqn:gamma-final}) follows from these expressions and from Lemma \ref{Lem:concl-calc-gamma}.  \hfill $\Box$

\begin{Cor}\label{Cor:metric-comparison} For all $[\theta_t]\in H^{0,\,1}(X_t,\,T^{1,\,0}X_t)_{[\omega^{n-1}]}\setminus\{0\}$, we have

$$(G^{(2)}_{WP}-\gamma)_t([\theta_t],\,[\theta_t])=\frac{4\,||\zeta_t||^2}{i^{n^2}\,\int_{X_t}u_t\wedge\bar{u}_t}\geq 0,  \hspace{3ex} t\in\Delta_{[\omega^{n-1}]},$$

\noindent hence the Hermitian metric $\omega^{(2)}_{WP}$ on $\Delta_{[\omega^{n-1}]}$ defined by $G^{(2)}_{WP}$ is bounded below by the K\"ahler metric $\gamma$.

\end{Cor}

 It is now clear that the obstruction to the metrics $\omega^{(2)}_{WP}$ and $\gamma$ coinciding on $\Delta_{[\omega^{n-1}]}$ is the possible negative answer to Question \ref{Question:prim-dclosed} in the case of balanced, non-K\"ahler fibres. Indeed, if every class in $H^{n-1,\,1}_{prim}(X_t,\,\C)$ could be represented by a form $\eta_t\lrcorner u_t$ that is both {\it primitive} and $d$-{\it closed}, we would have, thanks to Lemma \ref{Lem:n1-1decompositions}, that $\star\overline{(\eta_t\lrcorner u_t)} = c\,\overline{(\eta_t\lrcorner u_t)}$ with $c=-1$ (if $n$ is even), $c=-i$ (if $n$ is odd). Hence, from Lemma \ref{Lem:concl-calc-gamma}, we would get $\omega^{(2)}_{WP}=\gamma$ as in the case of K\"ahler polarised deformations of [Tia87] since formula (\ref{eqn:WP2}) can be re-written in the following obvious way\!\!:

\vspace{2ex}

\hspace{10ex} $\displaystyle G^{(2)}_{WP}([\theta_t],\,[\eta_t]) = \frac{\int\limits_{X_t}(\theta_t\lrcorner u_t)\wedge\star\overline{(\eta_t\lrcorner u_t)}}{i^{n^2}\,\int\limits_{X_t}u_t\wedge\bar{u}_t}$

\section{Balanced holomorphic symplectic $\partial\bar\partial$-manifolds}\label{section:hol-symp}

\subsection{Primitive $(1,\, 1)$-classes on balanced manifolds}\label{subsection:primitive11}

 Let $(X, \, \omega)$ be a compact, balanced manifold ($\mbox{dim}_{\C}X=n$). The balanced class $[\omega^{n-1}]\in H^{n-1,\,n-1}(X,\,\C)$ enables one to define the notion of {\it primitive} $2$-classes on $X$ in the same way as in the standard K\"ahler case. Indeed, at the level of Dolbeault cohomology, the linear operator

\begin{eqnarray}\label{eqn:Ln-1}L_{\omega}^{n-1}: H^{1,\,1}(X,\,\C)\rightarrow H^{n,\, n}(X,\,\C)\simeq\C, \hspace{3ex} [\alpha]\mapsto [\omega^{n-1}\wedge\alpha],\end{eqnarray}

\noindent is well defined because, thanks to the {\it balanced} assumption on $\omega$, $\bar\partial(\omega^{n-1}\wedge\alpha)=0$ whenever $\bar\partial\alpha=0$ and $\omega^{n-1}\wedge\alpha=\bar\partial(\omega^{n-1}\wedge\beta)$ whenever $\alpha=\bar\partial\beta$ is $\bar\partial$-exact. We can then call {\it primitive} those classes that are in the kernel of $L_{\omega}^{n-1}$, i.e.

\begin{eqnarray}\label{eqn:primclasses11}H^{1,\,1}_{prim}(X,\,\C):=\{[\alpha]\in H^{1,\,1}(X,\,\C)\,\,;\,\,\omega^{n-1}\wedge\alpha \hspace{1ex} \mbox{is}\hspace{1ex}\bar\partial-\mbox{exact}\}.\end{eqnarray}

\noindent Analogous definitions can be made for De Rham $2$-classes and Dolbeault $(2,\, 0)$ and $(0,\, 2)$-classes, but all $(2,\, 0)$ and $(0,\, 2)$-classes are primitive for trivial bidegree reasons. Thus, if the $\partial\bar\partial$-lemma is supposed to hold on $X$, the Hodge decomposition $H^2(X,\, \C) = H^{2,\,0}(X,\, \C)\oplus H^{1,\,1}(X,\, \C)\oplus H^{0,\,2}(X,\, \C)$ shows that only the $H^{1,\,1}(X,\, \C)$ component supports a nontrivial notion of primitivity. Notice that for $k>2$, there is no corresponding notion of primitive $k$-classes if $\omega$ is only balanced since $\omega^{n-k+1}$ is not closed unless $\omega$ is K\"ahler. It had to be replaced in bidegree $(n-1,\,1)$ by the ad-hoc definition \ref{Def:primitive-classes_n-11} using the Calabi-Yau isomorphism when $K_X$ was assumed to be trivial.

\begin{Lem}\label{Lem:prim11class-form} Let $(X, \, \omega)$ be a compact, balanced manifold ($\mbox{dim}_{\C}X=n$). Then a class $[\alpha]\in H^{1,\,1}(X,\,\C)$ is primitive if and only if it can be represented by a primitive form.

\end{Lem}

\noindent {\it Proof.} By the standard definition (applicable to any Hermitian metric $\omega$), a $(1,\,1)$-form $\alpha$ is primitive if $\omega^{n-1}\wedge\alpha=0$. It is thus obvious that any class representable by a primitive form is primitive. To see the converse, pick any class $[\alpha]\in H^{1,\,1}_{prim}(X,\,\C)$ and any representative $\alpha$. We have to prove the existence of a $(1,\, 0)$-form $u$ such that the representative $\alpha + \bar\partial u$ of $[\alpha]$ is primitive. This amounts to $\omega^{n-1}\wedge(\alpha+\bar\partial u)=0$, which is equivalent to $\bar\partial(\omega^{n-1}\wedge u)=-\omega^{n-1}\wedge\alpha$ thanks to the balanced assumption $\bar\partial\omega^{n-1}=0$. Now, $\omega^{n-1}\wedge\alpha$ is $\bar\partial$-exact by the primitivity assumption on the class $[\alpha]$. Pick any $w\in C^{\infty}_{n,\,n-1}(X,\,\C)$ such that $\bar\partial w=-\omega^{n-1}\wedge\alpha$. It thus suffices to prove the existence of a $(1,\, 0)$-form $u$ such that $\omega^{n-1}\wedge u=w$. The linear operator

\begin{eqnarray}\label{eqn:HL10}L_{\omega}^{n-1} : C^{\infty}_{1,\,0}(X,\,\C)\rightarrow C^{\infty}_{n,\,n-1}(X,\,\C), \hspace{3ex} u\mapsto\omega^{n-1}\wedge u,\end{eqnarray}

\noindent is an isomorphism (for any Hermitian metric $\omega$), so there is a unique $(1,\, 0)$-form $u$ such that $\omega^{n-1}\wedge u=w$.  \hfill $\Box$

\vspace{2ex}

 The primitive representative of a primitive class $[\alpha]\in H^{1,\,1}(X,\,\C)$ need not be unique, but we can single out a particular one that is uniquely determined by the metric $\omega$ in the given primitive class in the following way.

\vspace{2ex}

\noindent {\bf Choice of a primitive representative\!\!:} {\it given a primitive $(1,\,1)$-class, let $\alpha$ be its $\Delta''_{\omega}$-harmonic representative. Then choose $w\in C^{\infty}_{n\,n-1}(X,\,\C)$ to be the solution of minimal $L^2$-norm (w.r.t. $\omega$) of the equation $\bar\partial w=-\omega^{n-1}\wedge\alpha$. Since the map (\ref{eqn:HL10}) is an isomorphism, the $(1,\, 0)$-form $u$ such that $\omega^{n-1}\wedge u=w$ is uniquely determined by $w$. Since the above choices of $\alpha$ and $w$ make them unique, the primitive representative $\alpha + \bar\partial u$ of the primitive class $[\alpha]$ is uniquely determined in this way by $\omega$ and $[\alpha]\in H^{1,\,1}_{prim}(X,\,\C)$.}  \hfill $(\star)$

\vspace{2ex}

 When $\omega$ is K\"ahler, the $\Delta''_{\omega}$-harmonic representative $\alpha$ of a primitive class is a {\it primitive} form, a standard fact that follows from $\Delta''_{\omega}$ and $L_{\omega}$ commuting (as can be easily seen from the K\"ahler identities). Thus $\omega^{n-1}\wedge\alpha=0$, hence $w=0$ is the minimal $L^2$-norm solution of equation $\bar\partial w=-\omega^{n-1}\wedge\alpha$. Consequently, $u=0$ and $\alpha + \bar\partial u = \alpha$, showing that our choice $(\star)$ of primitive representative coincides with the standard $\Delta''_{\omega}$-harmonic choice when $\omega$ is K\"ahler. However, when $\omega$ is only balanced, it is not clear whether the $\Delta''_{\omega}$-harmonic representative of a primitive class is a {\it primitive} form. This accounts for the need of introducing the choice $(\star)$.

\subsection{Co-polarised deformations of holomorphic symplectic manifolds}\label{subsection:copol-def}

 Let $(X, \, \omega)$ be a compact, balanced $\partial\bar\partial$-manifold ($\mbox{dim}_{\C}X=n$). Suppose there exists a $C^{\infty}$ $\bar\partial$-closed $(2,\, 0)$-form $\sigma$ that is non-degenerate at every point of $X$ and that such a $\sigma$ is unique up to a nonzero constant factor. Thus $H^{2,\, 0}(X,\, \C)\simeq\C$ and $\sigma$ defines a holomorphic symplectic structure on $X$. The form $\sigma$ naturally identifies with the class $[\sigma]\in H^{2,\,0}(X,\,\C)$. 

 It follows from the $\partial\bar\partial$-assumption on $X$ that $\sigma$ is actually $d$-closed by the following observation which is standard when $X$ is K\"ahler (and probably also under the weaker $\partial\bar\partial$-assumption). The standard K\"ahler-case proof, using the Laplacian equality $\Delta'=\Delta''$, no longer holds in the $\partial\bar\partial$-case for which we spell out the argument below for the sake of completeness.

\begin{Lem}\label{Lem:holpforms--d-closed} Every holomorphic $p$-form is $d$-closed on any compact complex $\partial\bar\partial$-manifold $X$ for any $0\leq p\leq n=\mbox{dim}_{\C}X$.

\end{Lem}

\noindent {\it Proof.} Fix any $p$ and let $\alpha\in C^{\infty}_{p,\, 0}(X,\, \C)$ be $\bar\partial$-closed. To show that $d\alpha=0$, it suffices to show that $\partial\alpha=0$. Now, $\partial\alpha$ is $\bar\partial$-closed since $\alpha$ is, while $\partial$ and $\bar\partial$ anti-commute. Thus $\partial\alpha$ is a $d$-closed, $\partial$-exact form of pure type $(p+1,\, 0)$. By the $\partial\bar\partial$-lemma, $\partial\alpha$ must be $\partial\bar\partial$-exact, i.e. $\partial\alpha = \partial\bar\partial\beta$ for some $(p,\, -1)$-form $\beta$. Since $\beta$ must vanish for type reasons, $\partial\alpha$ vanishes.  \hfill $\Box$

\vspace{2ex}

 We are now ready to connect the primitive $(1,\,1)$-cohomology to the parameter space of co-polarised deformations defined by a balanced class via the natural isomorphism associated with the holomorphic symplectic structure.

\begin{Lem}\label{Lem:hol-symp-iso} Let $X$ be a compact complex manifold ($\mbox{dim}_{\C}X=n$) admitting a holomorphic symplectic structure $\sigma$ that is unique up to a constant factor.

\noindent $(i)$\, The linear map defined by $\sigma$ as

\begin{eqnarray}\label{eqn:hol-symp-iso-forms}T_{\sigma}\,\,:\,\,C^{\infty}_{0,\,1}(X,\,T^{1,\,0}X)\stackrel{\cdot\lrcorner\sigma}{\longrightarrow}C^{\infty}_{1,\,1}(X,\,\C), \hspace{3ex} \theta\mapsto T_{\sigma}(\theta)\!:=\theta\lrcorner\sigma,\end{eqnarray}

\noindent is an isomorphism satisfying the following properties\!\!:

\begin{eqnarray}\label{eqn:Tsigma-dbar}T_{\sigma}(\ker\bar\partial)=\ker\bar\partial \hspace{3ex} \mbox{and} \hspace{3ex} T_{\sigma}(\mbox{Im}\,\bar\partial)=\mbox{Im}\,\bar\partial.\end{eqnarray}

\noindent Consequently, $T_{\sigma}$ induces an isomorphism in cohomology

\begin{eqnarray}\label{eqn:hol-symp-iso-classes}T_{[\sigma]}\,\,;\,\,H^{0,\,1}(X,\,T^{1,\,0}X)\stackrel{\cdot\lrcorner[\sigma]}{\longrightarrow}H^{1,\,1}(X,\,\C)\end{eqnarray}

\noindent defined by $T_{[\sigma]}([\theta])=[\theta\lrcorner\sigma]$ for all $[\theta]\in H^{0,\,1}(X,\,T^{1,\,0}X)$.

\vspace{2ex}

\noindent $(ii)$\, If $\omega$ is a {\bf balanced} metric on $X$, then the image under $T_{[\sigma]}$ of the subspace $H^{0,\,1}(X,\,T^{1,\,0}X)_{[\omega^{n-1}]}\subset H^{0,\,1}(X,\,T^{1,\,0}X)$ defined in (\ref{eqn:co-pol-space}) is the subspace $H^{1,\,1}_{prim}(X,\,\C)\subset H^{1,\,1}(X,\,\C)$ of primitive $(1,\,1)$-classes defined in (\ref{eqn:primclasses11}), i.e.

\begin{eqnarray}\label{eqn:hol-symp-iso-primclasses}T_{[\sigma]}\,\,:\,\,H^{0,\,1}(X,\,T^{1,\,0}X)_{[\omega^{n-1}]}\stackrel{\simeq}{\longrightarrow}H^{1,\,1}_{prim}(X,\,\C).\end{eqnarray}

\end{Lem}

\noindent {\it Proof.} It is clear that $T_{\sigma}$ is an isomorphism. As in the proof of Lemma \ref{Lem:3space-decomp-contr}, the rest of $(i)$ follows from the easy-to-check formulae

\begin{eqnarray}\label{eqn:Leibniz-dbar-sigma}\bar\partial(\theta\lrcorner\sigma)=(\bar\partial\theta)\lrcorner\sigma + \theta\lrcorner(\bar\partial\sigma)=(\bar\partial\theta)\lrcorner\sigma, \hspace{1ex} \bar\partial(\xi\lrcorner\sigma)=(\bar\partial\xi)\lrcorner\sigma - \xi\lrcorner(\bar\partial\sigma)=(\bar\partial\xi)\lrcorner\sigma\end{eqnarray}

\noindent for all $\theta\in C^{\infty}_{0, \, 1}(X, \, T^{1, \, 0}X)$ and all $\xi\in C^{\infty}(X, \, T^{1, \, 0}X)$ which readily imply the inclusions $T_{\sigma}(\ker\bar\partial)\subset\ker\bar\partial$ and $T_{\sigma}(\mbox{Im}\,\bar\partial)\subset\mbox{Im}\,\bar\partial$.

 Let us prove, for example, the identity $\mbox{Im}\,\bar\partial = T_{\sigma}(\mbox{Im}\,\bar\partial)$. This amounts to proving that $\theta$ is $\bar\partial$-exact if and only if $\theta\lrcorner\sigma$ is $\bar\partial$-exact. Having fixed local holomorphic coordinates $z_1, \dots , z_n$ on some open subset $U\subset X$, let

$$\theta=\sum\limits_{\alpha,\,\beta}\theta^{\alpha}_{\beta}\,\frac{\partial}{\partial z^{\alpha}}\,d\bar{z}^{\beta} \hspace{2ex}  \mbox{and}  \hspace{2ex} \sigma=\sum\limits_{\alpha,\,\delta}\sigma_{\alpha,\,\delta}\,dz^{\alpha}\wedge dz^{\delta},$$ 

\noindent where the coefficients $\sigma_{\alpha,\,\delta}$ are holomorphic functions (since $\sigma$ is holomorphic) and the matrix $(\sigma_{\alpha,\,\delta})_{\alpha,\,\delta}$ is invertible at every point since $\sigma$ is non-degenerate at every point. Then $\theta\lrcorner\sigma=\sum\limits_{\alpha,\,\beta,\,\delta}\theta^{\alpha}_{\beta}\,(\sigma_{\alpha,\,\delta}-\sigma_{\delta,\,\alpha})\, d\bar{z}^{\beta}\wedge dz^{\delta}$. Thus $\theta\lrcorner\sigma$ is $\bar\partial$-exact if and only if there exists a $(1,\,0)$-form $v=\sum\limits_{\delta}v_{\delta}\,dz^{\delta}$ such that $\theta\lrcorner\sigma=\bar\partial v$, which amounts to

$$\sum\limits_{\alpha,\,\beta,\,\delta}\theta^{\alpha}_{\beta}\,(\sigma_{\alpha,\,\delta}-\sigma_{\delta,\,\alpha})\, d\bar{z}^{\beta}\wedge dz^{\delta}=\sum\limits_{\delta,\,\beta}\frac{\partial v_{\delta}}{\partial\bar{z}^{\beta}}\,d\bar{z}^{\beta}\wedge dz^{\delta} \hspace{1ex} \Leftrightarrow \hspace{1ex} \sum\limits_{\alpha}\theta^{\alpha}_{\beta}\,(\sigma_{\alpha,\,\delta}-\sigma_{\delta,\,\alpha})= \frac{\partial v_{\delta}}{\partial\bar{z}^{\beta}}$$

\noindent for all $\beta, \,\delta$. The last identity is equivalent to $\theta^{\alpha}_{\beta}=\sum\limits_{\delta}\frac{\partial v_{\delta}}{\partial\bar{z}^{\beta}}\,(\sigma^{\delta,\,\alpha}-\sigma^{\alpha,\,\delta})=\frac{\partial}{\partial\bar{z}^{\beta}}(\sum\limits_{\delta}(\sigma^{\delta,\,\alpha}-\sigma^{\alpha,\,\delta})\,v_{\delta})$ for all $\alpha,\,\beta$, where the matrix $(\sigma^{\delta,\,\alpha})_{\alpha,\,\delta}$ is the inverse of $(\sigma_{\alpha,\,\delta})_{\alpha,\,\delta}$. (We have used the fact that the $\sigma^{\delta,\,\alpha}$'s are holomorphic functions since the $\sigma_{\alpha,\,\delta}$'s are.) This, in turn, is equivalent to $\theta=\bar\partial(\sum\limits_{\alpha}(\sum\limits_{\delta}(\sigma^{\delta,\,\alpha}-\sigma^{\alpha,\,\delta})\,v_{\delta})\,\frac{\partial}{\partial z^{\alpha}})$, i.e. to $\theta$ being $\bar\partial$-exact. We have thus proved that $\theta\lrcorner\sigma$ is $\bar\partial$-exact if and only if $\theta$ is $\bar\partial$-exact, i.e. the latter identity in (\ref{eqn:Tsigma-dbar}). 

 The remaining inclusion in the former identity of (\ref{eqn:Tsigma-dbar}) is proved in a similar way. 

 The proof of $(ii)$ will run in two steps. First we prove the inclusion

\begin{eqnarray}\label{eqn:direct-incl-prim-hol}T_{[\sigma]}\bigg(H^{0,\,1}(X,\,T^{1,\,0}X)_{[\omega^{n-1}]}\bigg)\subset H^{1,\,1}_{prim}(X,\,\C),\end{eqnarray}  

\noindent which amounts to proving that for every class $[\theta]\in H^{0,\,1}(X,\,T^{1,\,0}X)$ for which $\theta\lrcorner\omega^{n-1}$ is $\bar\partial$-exact, $\omega^{n-1}\wedge(\theta\lrcorner\sigma)$ is also $\bar\partial$-exact. Now, we always have

$$0=\theta\lrcorner(\omega^{n-1}\wedge\sigma) = (\theta\lrcorner\omega^{n-1})\wedge\sigma + \omega^{n-1}\wedge(\theta\lrcorner\sigma),$$

\noindent where the first identity follows from the fact that $\omega^{n-1}\wedge\sigma$ is of type $(n+1,\, n-1)$, hence vanishes. Thus

$$(\theta\lrcorner\omega^{n-1})\wedge\sigma = - \omega^{n-1}\wedge(\theta\lrcorner\sigma) \hspace{2ex} \mbox{for all}\hspace{1ex} \theta\in C^{\infty}_{0,\,1}(X,\,T^{1,\,0}X).$$ 

\noindent Now, if $\theta\lrcorner\omega^{n-1}$ is supposed to be $\bar\partial$-exact, then $(\theta\lrcorner\omega^{n-1})\wedge\sigma$ is $\bar\partial$-exact, too, since $\sigma$ is $\bar\partial$-closed. Hence $\omega^{n-1}\wedge(\theta\lrcorner\sigma)$ is $\bar\partial$-exact whenever $\theta\lrcorner\omega^{n-1}$is, proving the inclusion (\ref{eqn:direct-incl-prim-hol}). 

 Since $T_{[\sigma]}$ is injective by $(i)$, it suffices to prove the dimension equality

\begin{eqnarray}\label{eqn:dim-equality}\mbox{dim}\,H^{0,\,1}(X,\,T^{1,\,0}X)_{[\omega^{n-1}]} = \mbox{dim}\, H^{1,\,1}_{prim}(X,\,\C)\end{eqnarray}

\noindent to be able to conclude that the inclusion (\ref{eqn:direct-incl-prim-hol}) is actually an identity. 

 By definition (\ref{eqn:primclasses11}), we have

$$H^{1,\,1}_{prim}(X,\,\C) = \ker\bigg(L_{\omega}^{n-1}\,\,:\,\,H^{1,\,1}(X,\,\C)\rightarrow H^{n,\,n}(X,\,\C)\simeq\C \bigg).$$

\noindent The linear map (\ref{eqn:Ln-1}) cannot vanish identically, so it is surjective. Hence

\begin{eqnarray}\label{eqn:dim-h11prim}\mbox{dim}\, H^{1,\,1}_{prim}(X,\,\C) = h^{1,\,1}-1,\end{eqnarray}

\noindent where $h^{1,\,1}\!:=\mbox{dim}\, H^{1,\,1}(X,\,\C)$. Meanwhile, definition (\ref{eqn:co-pol-space}) translates to

$$H^{0,\,1}(X,\,T^{1,\,0}X)_{[\omega^{n-1}]} = \ker\bigg(H^{0,\,1}(X,\,T^{1,\,0}X)\ni [\theta]\stackrel{T_{[\omega^{n-1}]}}{\longmapsto}[\theta\lrcorner\omega^{n-1}]\in H^{n-2,\,n}(X,\,\C)\bigg),$$

\noindent while $H^{n-2,\,n}(X,\,\C)\simeq H^{2,\,0}(X,\,\C)\simeq\C$ by Serre duality and the uniqueness (up to a constant factor) assumption on the holomorphic symplectic structure $[\sigma]\in H^{2,\,0}(X,\,\C)$. It is clear that the linear map $T_{[\omega^{n-1}]}$ does not vanish identically, so it must be surjective. Thus we get

\begin{eqnarray}\label{eqn:dim-h01co-pol}\mbox{dim}\, H^{0,\,1}(X,\,T^{1,\,0}X)_{[\omega^{n-1}]} = \mbox{dim}\, H^{0,\,1}(X,\,T^{1,\,0}X) - 1 =  h^{1,\,1}-1,\end{eqnarray}

\noindent where the last identity follows from the isomorphism (\ref{eqn:hol-symp-iso-classes}) dealt with under $(i)$. It is now clear that the dimension equality (\ref{eqn:dim-equality}) is a consequence of the combined identities (\ref{eqn:dim-h11prim}) and (\ref{eqn:dim-h01co-pol}). The proof is complete.  \hfill $\Box$

\vspace{3ex}

 The use of the isomorphism $T_{[\sigma]}$ in (\ref{eqn:hol-symp-iso-classes}) in the holomorphic symplectic case may be an alternative to the use of the isomorphism $T_{[u]}$ in (\ref{eqn:u-isom-classes}) of the more general Calabi-Yau case while running the construction of the Weil-Petersson metrics of section \ref{section:period-WP}.

\section{Appendix}\label{section:appendix}

We start by briefly recalling Wu's argument in [Wu06] proving the deformation openness of the simultaneous occurence of the $\partial\bar\partial$ and balanced properties (and even more). 


\begin{The}(C.-C. Wu [Wu06])\label{The:Wu_bal-ddbar-open_appendix} Let $(X_t)_{t\in\Delta}$ be a holomorphic family of compact complex manifolds. 

If the fibre $X_0$ is a balanced $\partial\bar\partial$-manifold, the fibre $X_t$ is again a balanced $\partial\bar\partial$-manifold for every $t\in\Delta$ sufficiently close to $0$. 

Moreover, every balanced metric $\omega_0$ on $X_0$ deforms to a family of balanced metrics $\omega_t$ on $X_t$ varying in a $C^\infty$ way with $t$ for $t$ in a small enough neighbourhood of $0$.

\end{The}

\noindent {\it Proof.} We reproduce Wu's arguments in a slightly different notation. Let $(\gamma_t)_{t\in\Delta}$ be an arbitrary $C^\infty$ family of Hermitian metrics on the fibres $(X_t)_{t\in\Delta}$. If $\Delta_{BC}(t)$ denotes the Bott-Chern Laplacian (cf. [KS60]) induced by the metric $\gamma_t$, the following $3$-space orthogonal decomposition is well-known (see e.g. [Sch07] or [Pop15] for some background) in every bidegree $(p,\,q)$: \begin{equation}\label{eqn:3-space_decomp_appendix} C^\infty_{p,\,q}(X_t,\,\C) = \ker\Delta_{BC}(t)\oplus\mbox{Im}\,(\partial_t\bar\partial_t)\oplus (\mbox{Im}\,\partial_t^\star + \mbox{Im}\,\bar\partial_t^\star),  \hspace{3ex} t\in\Delta,  \end{equation}

\noindent where $\ker\partial_t\cap\ker\bar\partial_t = \ker\Delta_{BC}(t)\oplus\mbox{Im}\,(\partial_t\bar\partial_t)$. Letting $F_t$ stand for the orthogonal projection w.r.t. the $L^2_{\gamma_t}$ inner product onto $\ker\Delta_{BC}(t)$ and letting $\Delta_{BC}^{-1}(t)$ stand for the Green operator of the elliptic operator $\Delta_{BC}(t)$, every form $\alpha_t\in C^\infty_{p,\,q}(X_t,\,\C)$ splits uniquely as $\alpha_t = F_t\alpha_t + \Delta_{BC}(t)\Delta_{BC}^{-1}(t)\alpha_t$. Moreover, if $\alpha_t\in\ker\partial_t\cap\ker\bar\partial_t$, this splitting reduces to

$$\alpha_t = F_t\alpha_t + \partial_t\bar\partial_t(\partial_t\bar\partial_t)^\star\Delta_{BC}^{-1}(t)\,\alpha_t.$$

\noindent (See Wu's original argument or the later Theorem 4.1 in [Pop15].)

Let $\omega_0$ be a balanced metric on $X_0$ and $n$ the complex dimension of $X_t$. Then $\omega^{n-1}_0\in\ker\partial_0\cap\ker\bar\partial_0$, so $\omega^{n-1}_0 = F_0\,\omega^{n-1}_0 + \partial_0\bar\partial_0(\partial_0\bar\partial_0)^\star\Delta_{BC}^{-1}(0)\,\omega^{n-1}_0$. Extend $\omega_0$ in an arbitrary way to Hermitian metrics $\widetilde\omega_t$ varying in a $C^\infty$ way with $t$ on the nearby fibres $X_t$ such that $\widetilde\omega_0=\omega_0$. Put

$$\Omega_t:=\mbox{Re}\, (F_t\,\widetilde\omega_t^{n-1} + \partial_t\bar\partial_t(\partial_t\bar\partial_t)^\star\Delta_{BC}^{-1}(t)\,\widetilde\omega_t^{n-1}),  \hspace{3ex} t\in\Delta.$$

\noindent By construction, every $\Omega_t$ is a $C^\infty$, real, $J_t$-type $(n-1,\,n-1)$-form on $X_t$ such that $d\Omega_t=0$ for every $t$. Moreover, $\Omega_0=\omega^{n-1}_0$.  

Now, since $X_0$ is a $\partial\bar\partial$-manifold, the fibres $X_t$ are again $\partial\bar\partial$-manifolds for every $t$ close to $0$ by Wu's first main result in [Wu06] and the dimensions $h^{p,\,q}(t)$ of the Bott-Chern cohomology spaces $H^{p,\,q}_{BC}(t)$ are independent of $t$ close to $0$ by Wu's main technical preliminary result. Thanks to the Hodge isomorphism $H^{p,\,q}_{BC}(t)\simeq\ker\Delta_{BC}(t)$ and to the classical Kodaira-Spencer theory for smooth families of elliptic operators, this implies that the operators $F_t$ and $\Delta_{BC}(t)^{-1}$ vary in a $C^\infty$ way with $t$. Therefore, the real differential forms $\Omega_t$ vary in a $C^\infty$ way with $t$. Since $\Omega_0=\omega^{n-1}_0>0$, we get by continuity that $\Omega_t>0$ for every $t$ sufficiently close to $0$. 

Taking the (unique) $(n-1)^{st}$ root $\omega_t>0$ of $\Omega_t>0$ for $t$ close to $0$, we get a $C^\infty$ family of {\it balanced} metrics $\omega_t$ on the fibres $X_t$ whose element corresponding to $t=0$ coincides with the original $\omega_0$. \hfill $\Box$

\vspace{3ex}

We can now prove the following observation that was used in the paper. While independent of the above approach of Wu, the proof uses similar techniques and, in particular, reproves Theorem \ref{The:Wu_bal-ddbar-open_appendix}.

\begin{Obs}\label{Obs:balanced-def-class_appendix} Let $(X_t)_{t\in\Delta}$ be a holomorphic family of $n$-dimensional compact complex manifolds such that the fibre $X_0$ is a balanced $\partial\bar\partial$-manifold. We denote by $X$ the differentiable manifold underlying the fibres $X_t$ (after possibly shrinking $\Delta$ about $0$.)

Let $\omega_0$ be a balanced metric on $X_0$ and suppose that the De Rham class $\{\omega_0^{n-1}\}_{DR}\in H^{2n-2}_{DR}(X,\,\C)$ is of type $(n-1,\,n-1)$ for the complex structure $J_t$ of $X_t$ for all $t$ close to zero and lying on a path through $0$ in $\Delta$.

Then, the De Rham class $\{\omega_0^{n-1}\}_{DR}$ contains a $J_t$-balanced metric for every $t$ as above sufficiently close to $0$.

\end{Obs}

\noindent {\it Proof.} Since $X_t$ is a $\partial\bar\partial$-manifold for every $t$ close to $0$, there are canonical isomorphisms $H_{BC}^{p,\,q}(X_t,\,\C)\simeq H_A^{p,\,q}(X_t,\,\C)$ (for every $(p,\,q)$) and \begin{eqnarray}\label{eqn:BC-A_splittings_appendix}\nonumber H^{2n-2}_{DR}(X,\,\C) & \simeq & H_{BC}^{n,\,n-2}(X_t,\,\C)\oplus H_{BC}^{n-1,\,n-1}(X_t,\,\C)\oplus H_{BC}^{n-2,\,n}(X_t,\,\C) \\
\nonumber & \simeq & H_A^{n,\,n-2}(X_t,\,\C)\oplus H_A^{n-1,\,n-1}(X_t,\,\C)\oplus H_A^{n-2,\,n}(X_t,\,\C).\end{eqnarray}

 Now, let $\omega_0^{n-1} = \Omega_t^{n,\,n-2} + \Omega_t^{n-1,\,n-1} + \Omega_t^{n-2,\,n}$ be the splitting of $\omega_0^{n-1}$ into components of pure $J_t$-types. In particular, $\Omega_t^{n-1,\,n-1}$ is a {\it real} $J_t$-type $(n-1,\,n-1)$-form that varies in a $C^\infty$ way with $t$ and is positive definite for every $t$ sufficiently close to $0$ since $\Omega_0^{n-1,\,n-1} = \omega_0^{n-1}>0$.

Meanwhile, since $d\omega_0^{n-1}=0$, it is easy to see (cf. e.g. [Pop15]) that $\partial_t\bar\partial_t\Omega_t^{n-1,\,n-1}=0$ and that the Aeppli cohomology class $[\Omega_t^{n-1,\,n-1}]_A$ is the image of $\{\omega_0^{n-1}\}_{DR}$ under the projection $H^{2n-2}_{DR}(X,\,\C)\longrightarrow H_A^{n-1,\,n-1}(X_t,\,\C)$ defined by the latter cohomology splitting above. 

To construct the image of $[\Omega_t^{n-1,\,n-1}]_A\in H_A^{n-1,\,n-1}(X_t,\,\C)$ in $H_{BC}^{n-1,\,n-1}(X_t,\,\C)$ under the canonical isomorphism $H_A^{n-1,\,n-1}(X_t,\,\C)\simeq H_{BC}^{n-1,\,n-1}(X_t,\,\C)$, we can proceed as in [Pop15] and look for the ``most economic choice'' of a $J_t$-$(n-2,\,n-1)$-form $u_t$ and a $J_t$-$(n-1,\,n-2)$-form $v_t$ such that the following $J_t$-$(n-1,\,n-1)$-form

$$\widetilde{\Omega}_t^{n-1,\,n-1}:=\Omega_t^{n-1,\,n-1} + \partial_t u_t + \bar\partial_t v_t$$

\noindent is $d$-closed. This amounts to $\partial_t\bar\partial_t u_t = \bar\partial_t\Omega_t^{n-1,\,n-1}$ and $\partial_t\bar\partial_t v_t = -\partial_t\Omega_t^{n-1,\,n-1}$. If we choose $v_t:=\bar{u}_t$, the latter equation becomes redundant, while the minimal $L^2_{\gamma_t}$-norm solution of the former equation (which is solvable since $X_t$ is a $\partial\bar\partial$-manifold) is given by the following Neumann-type formula (see [Wu06] or [Pop15])\!:

$$u_t=(\partial_t\bar\partial_t)^\star\Delta_{BC}^{-1}(t)\bar\partial_t\Omega_t^{n-1,\,n-1},  \hspace{3ex} t\in\Delta,$$ 

\noindent after possibly shrinking $\Delta$ about $0$ to ensure that $X_t$ is a $\partial\bar\partial$-manifold. (As usual, we have fixed an arbitrary $C^\infty$ family $(\gamma_t)_{t\in\Delta}$ of Hermitian metrics on the fibres $(X_t)_{t\in\Delta}$.) 

Then, for all $t$ close to $0$, we get $$\widetilde{\Omega}_t^{n-1,\,n-1}:=\Omega_t^{n-1,\,n-1} + \partial_t(\partial_t\bar\partial_t)^\star\Delta_{BC}^{-1}(t)\bar\partial_t\Omega_t^{n-1,\,n-1} + \bar\partial_t(\bar\partial_t\partial_t)^\star\overline{\Delta_{BC}^{-1}(t)}\partial_t\Omega_t^{n-1,\,n-1}.$$

\noindent When $t=0$, $\partial_0\bar\partial_0 u_0 = \bar\partial_0\Omega_0^{n-1,\,n-1} = \bar\partial_0\omega_0^{n-1} = 0$ (the last identity being a consequence of $\omega_0$ being balanced), so the minimal $L^2$-norm solution of this equation is $u_0=0$. Note that $u_t$, hence also $\widetilde{\Omega}_t^{n-1,\,n-1}$, depends in a $C^\infty$ way on $t$ for the same reason as in Wu's proof of Theorem \ref{The:Wu_bal-ddbar-open_appendix}: the $\partial\bar\partial$-assumption implies the invariance w.r.t. $t$ of the Bott-Chern numbers $h^{p,\,q}_{BC}(t)$, which implies the smooth dependence on $t$ of $\Delta_{BC}^{-1}(t)$.

We have thus constructed a $C^\infty$ family of real $d$-closed $J_t$-$(n-1,\,n-1)$-forms $\widetilde{\Omega}_t^{n-1,\,n-1}$ such that $\widetilde{\Omega}_0^{n-1,\,n-1} = \omega_0^{n-1}>0$. By continuity, we must have $\widetilde{\Omega}_t^{n-1,\,n-1}>0$, hence $\widetilde{\Omega}_t^{n-1,\,n-1}$ defines a balanced metric on $X_t$, for all $t$ close to $0$. (In particular, this gives another proof of Wu's Theorem \ref{The:Wu_bal-ddbar-open_appendix}.) 

Moreover, $[\widetilde{\Omega}_t^{n-1,\,n-1}]_{BC}$ is the image in $H^{n-1,\,n-1}_{BC}(X_t,\,\C)$ of $[\Omega_t^{n-1,\,n-1}]_A$ under the canonical isomorphism $H^{n-1,\,n-1}_A(X_t,\,\C)\rightarrow H^{n-1,\,n-1}_{BC}(X_t,\,\C)$. Since $[\Omega_t^{n-1,\,n-1}]_A$ is the image of $\{\omega_0^{n-1}\}_{DR}$ under the canonical projection of $H^{2n-2}_{DR}(X,\,\C)$ onto $H_A^{n-1,\,n-1}(X_t,\,\C)$, we infer that $[\widetilde{\Omega}_t^{n-1,\,n-1}]_{BC}$ is the image in $H^{n-1,\,n-1}_{BC}(X_t,\,\C)$ of $\{\omega_0^{n-1}\}_{DR}$ under the canonical projection of $H^{2n-2}_{DR}(X,\,\C)$ onto $H_{BC}^{n-1,\,n-1}(X_t,\,\C)$. Meanwhile, if the class $\{\omega_0^{n-1}\}_{DR}\in H^{2n-2}_{DR}(X,\,\C)$ is supposed to be of $J_t$-type $(n-1,\,n-1)$, it coincides with its projection $[\widetilde{\Omega}_t^{n-1,\,n-1}]_{BC}$ (after the obvious canonical identification of $H^{n-1,\,n-1}_{BC}(X_t,\,\C)$ with its image in $H^{2n-2}_{DR}(X,\,\C)$). This means that $\{\omega_0^{n-1}\}_{DR} = \{\widetilde{\Omega}_t^{n-1,\,n-1}\}_{DR}$ for all $t$ sufficiently close to $0$ and lying on the path through $0$ in $\Delta$ along which $\{\omega_0^{n-1}\}_{DR}$ is assumed to be of $J_t$-type $(n-1,\,n-1)$. Thus, the class $\{\omega_0^{n-1}\}_{DR}$ contains the $J_t$-balanced metric $ \widetilde{\Omega}_t^{n-1,\,n-1}$ for all these $t$'s.  \hfill $\Box$

\vspace{3ex}

The other issue dealt with in this appendix is the following computation.

\vspace{2ex}

\noindent {\it Proof of Lemma \ref{Lem:dbar-contr-omega}.} Fix an arbitrary point $x_0\in X$ and let $z_1, \dots , z_n$ be local holomorphic coordinates about $x_0$. If we denote

$$\omega^{n-1}=i^{n-1}\,\sum\limits_{\alpha, \,\beta}\gamma_{\alpha\beta}\,\widehat{dz_{\alpha}\wedge d\bar{z}_{\beta}}  \hspace{3ex}\mbox{and}\hspace{3ex} \xi=\sum\limits_j\xi_j\,\frac{\partial}{\partial z_j},$$

\noindent as $\widehat{dz_{\alpha}\wedge d\bar{z}_{\beta}}:=dz_1\wedge\dots\wedge\widehat{dz_{\alpha}}\wedge\dots\wedge dz_n\wedge d\bar{z}_1\wedge\dots\wedge\widehat{d\bar{z}_{\beta}}\wedge\dots\wedge d\bar{z}_n$, we get \\

\noindent $\xi\lrcorner\omega^{n-1}  = $ 
\begin{eqnarray}\label{eqn:xi-omega_n-1}\nonumber i^{n-1}\,\sum\limits_{\stackrel{\beta}{j<\alpha}}(-1)^{j-1}\,\xi_j\gamma_{\alpha\beta}\,\widehat{(dz_j\wedge dz_{\alpha}\wedge d\bar{z}_{\beta})} + i^{n-1}\,\sum\limits_{\stackrel{\beta}{j>\alpha}}(-1)^j\,\xi_j\gamma_{\alpha\beta}\,\widehat{(dz_{\alpha}\wedge dz_j\wedge d\bar{z}_{\beta})} & & \\
\nonumber = i^{n-1}\,\sum\limits_{\stackrel{\beta}{j<\alpha}}\bigg((-1)^{j-1}\,\xi_j\gamma_{\alpha\beta} + (-1)^{\alpha}\,\xi_{\alpha}\gamma_{j\beta}\bigg)\,\widehat{(dz_j\wedge dz_{\alpha}\wedge d\bar{z}_{\beta})},\end{eqnarray}

\noindent where we have used the notation

\noindent $\widehat{(dz_j\wedge dz_{\alpha}\wedge d\bar{z}_{\beta})}:=dz_1\wedge\dots\wedge\widehat{dz_j}\wedge\dots\wedge\widehat{dz_{\alpha}}\wedge\dots\wedge dz_n\wedge d\bar{z}_1\wedge\dots\wedge\widehat{d\bar{z}_{\beta}}\wedge\dots\wedge d\bar{z}_n.$

\noindent Hence, by applying $\bar\partial$, we get

\begin{eqnarray}\label{eqn:dbar-xi-omega_n-1}\nonumber\bar\partial(\xi\lrcorner\omega^{n-1}) =  i^{n-1}\,\sum\limits_{\stackrel{\beta}{j<\alpha}}(-1)^{n+\beta-1}\bigg[(-1)^{j-1}\,\xi_j\frac{\partial\gamma_{\alpha\beta}}{\partial\bar{z}_{\beta}} & + & (-1)^{j-1}\,\frac{\partial\xi_j}{\partial\bar{z}_{\beta}}\gamma_{\alpha\beta} + (-1)^{\alpha}\xi_{\alpha}\frac{\partial\gamma_{j\beta}}{\partial\bar{z}_{\beta}} \\
 & + & (-1)^{\alpha}\frac{\partial\xi_{\alpha}}{\partial\bar{z}_{\beta}}\gamma_{j\beta}\bigg]\,\widehat{dz_j\wedge dz_{\alpha}}.\end{eqnarray}

\noindent Similar calculations yield

\begin{eqnarray}\label{eqn:dbar-xi_with_omega_n-1}(\bar\partial\xi)\lrcorner\omega^{n-1} = (-1)^ni^{n-1}\,\sum\limits_{\stackrel{\beta}{j<\alpha}}(-1)^{\beta}\bigg((-1)^j\frac{\partial\xi_j}{\partial\bar{z}_{\beta}}\gamma_{\alpha\beta} - (-1)^{\alpha}\frac{\partial\xi_{\alpha}}{\partial\bar{z}_{\beta}}\gamma_{j\beta}\bigg)\widehat{dz_j\wedge dz_{\alpha}},\end{eqnarray}

\noindent showing that $(\bar\partial\xi)\lrcorner\omega^{n-1}$ equals the sum of the second and fourth groups of terms in the expression (\ref{eqn:dbar-xi-omega_n-1}) for $(\bar\partial\xi)\lrcorner\omega^{n-1}$. On the other hand, we get

\begin{eqnarray}\label{eqn:(dbar-xi)-omega_n-1}\nonumber\bar\partial\omega^{n-1}=(-1)^ni^{n-1}\,\sum\limits_{\alpha, \beta}(-1)^{\beta}\frac{\partial\gamma_{\alpha\beta}}{\partial\bar{z}_{\beta}}\,\widehat{dz_{\alpha}},\end{eqnarray}

\noindent leading to

\begin{eqnarray}\label{eqn:xi-dbar-omega_n-1}\nonumber\xi\lrcorner\bar\partial\omega^{n-1} & = & (-1)^ni^{n-1}\,\sum\limits_{\stackrel{\beta}{j<\alpha}}(-1)^{j+\beta-1}\xi_j\frac{\partial\gamma_{\alpha\beta}}{\partial\bar{z}_{\beta}}\,\widehat{dz_j\wedge dz_{\alpha}}\\
\nonumber & + & (-1)^ni^{n-1}\,\sum\limits_{\stackrel{\beta}{j>\alpha}}(-1)^{j+\beta}\xi_j\frac{\partial\gamma_{\alpha\beta}}{\partial\bar{z}_{\beta}}\,\widehat{dz_{\alpha}\wedge dz_j}\\
\nonumber  & = & (-1)^ni^{n-1}\,\sum\limits_{\stackrel{\beta}{j<\alpha}}(-1)^{\beta}\bigg((-1)^{j-1}\xi_j\frac{\partial\gamma_{\alpha\beta}}{\partial\bar{z}_{\beta}} + (-1)^{\alpha}\xi_{\alpha}\frac{\partial\gamma_{j\beta}}{\partial\bar{z}_{\beta}}\bigg)\,\widehat{dz_j\wedge dz_{\alpha}}.\end{eqnarray}

\noindent Thus $\xi\lrcorner\bar\partial\omega^{n-1}$ equals the sum multiplied by $(-1)$ of the first and third groups of terms in the expression (\ref{eqn:dbar-xi-omega_n-1}) for $\bar\partial(\xi\lrcorner\omega^{n-1})$. Combining with (\ref{eqn:dbar-xi-omega_n-1}) and (\ref{eqn:dbar-xi_with_omega_n-1}), we get the identity claimed in $(i)$. Similar calculations prove $(ii)$.  \hfill $\Box$

\vspace{2ex}

\noindent {\bf References.} 

\vspace{1ex}

\noindent [AB90]\, L. Alessandrini, G. Bassanelli --- {\it Small Deformations of a Class of Compact Non-K\"ahler Manifolds} --- Proc. Amer. Math. Soc. {\bf 109} (1990), no. 4, 1059–1062. 

\vspace{1ex}

\noindent [AB93]\, L. Alessandrini, G. Bassanelli --- {\it Metric Properties of Manifolds Bimeromorphic to Compact K\"ahler Spaces} --- J. Diff. Geom. {\bf 37} (1993), 95-121.

\vspace{1ex}

\noindent [AB95]\, L. Alessandrini, G. Bassanelli --- {\it Modifications of Compact Balanced Manifolds} --- C.R. Acad. Sci. Paris, t. 320, S\'erie I Math. {\bf 320} (1995), no. 12, 1517–1522.

\vspace{1ex}

\noindent [AB96]\, L. Alessandrini, G. Bassanelli --- {\it The class of Compact Balanced Manifolds Is Invariant under Modifications} --- Complex Analysis and Geometry (Trento, 1993), 1–17, Lecture Notes in Pure and Appl. Math., {\bf 173}, Dekker, New York, 1996.

\vspace{1ex}

\noindent [Cam91a]\, F. Campana --- {\it The Class ${\cal C}$ Is Not Stable by Small Deformations} --- Math. Ann. {\bf 290} (1991) 19-30.

\vspace{1ex}

\noindent [Cam91b]\, F. Campana --- {\it On Twistor Spaces of the Class ${\cal C}$} --- J. Diff. Geom. {\bf 33} (1991) 541-549.

\vspace{1ex}

\noindent [Cam95]\, F. Campana --- {\it Remarques sur les groupes de K\"ahler nilpotents} --- Ann. Sci. \'Ecole Norm. Sup. (4) {\bf 28} (1995), no. 3, 307-316.

\vspace{1ex}

\noindent [Cam04]\, F. Campana --- {\it Orbifolds, Special Varieties and Classification Theory} --- Ann. Inst. Fourier (Grenoble) {\bf 54} (2004), no. 3, 499-630.

\vspace{1ex}

\noindent [Chi14]\, I. Chiose --- {\it Obstructions to the Existence of K\"ahler Structures on Compact Complex Manifolds} --- Proc AMS {\bf 142} (2014), no. 10, 3561–3568).

\vspace{1ex}

\noindent [DGMS75]\, P. Deligne, Ph. Griffiths, J. Morgan, D. Sullivan --- {\it Real Homotopy Theory of K\"ahler Manifolds} --- Invent. Math. {\bf 29} (1975), 245-274.

\vspace{1ex}

\noindent [Dem 97]\, J.-P. Demailly --- {\it Complex Analytic and Algebraic Geometry}---http://www-fourier.ujf-grenoble.fr/~demailly/books.html

\vspace{1ex}

\noindent [FOU14]\, A. Fino, A. Otal, L. Ugarte --- {\it Six Dimensional Solvmanifolds with Holomorphically Trivial Canonical Bundle} ---  Int. Math. Res. Not. IMRN 2015, no. {\bf 24}, 13757-13799.

\vspace{1ex}

\noindent [Fri17]\, R. Friedman --- {\it The $\partial\bar\partial$-Lemma for General Clemens Manifolds} --- arXiv e-print AG 1708.00828v1. 

\vspace{1ex}

\noindent [FLY12]\, J.Fu, J.Li, S.-T. Yau --- {\it Balanced Metrics on Non-K\"ahler Calabi-Yau Threefolds} --- J. Diff. Geom. {\bf 90} (2012) 81-129.

\vspace{1ex}

\noindent [Gau77a]\, P. Gauduchon --- {\it Le th\'eor\`eme de l'excentricit\'e nulle} --- C.R. Acad. Sc. Paris, S\'erie A, t. {\bf 285} (1977), 387-390.

\vspace{1ex}

\noindent [Gau77b]\, P. Gauduchon --- {\it Fibr\'es hermitiens \`a endomorphisme de Ricci non n\'egatif} -- Bull. SMF {\bf 105} (1977) 113-140.

\vspace{1ex}

\noindent [Gau91]\, P. Gauduchon --- {\it Structures de Weyl et th\'eor\`emes d'annulation sur une vari\'et\'e conforme autoduale} --- Ann. Scuola Norm. Sup. Pisa Cl. Sci. (4) {\bf 18} (1991), no. 4, 563–629.

\vspace{1ex}

\noindent [Kas13]\, H. Kasuya --- {\it Techniques of Computations of Dolbeault Cohomology of Solvmanifolds.  } ---  Math. Z. {\bf 273} (2013), no. 1-2, 437-447.

\vspace{1ex}

\noindent [KS60]\, K. Kodaira, D.C. Spencer --- {\it On Deformations of Complex Analytic Structures, III. Stability Theorems for Complex Structures} --- Ann. of Math. {\bf 71}, no.1 (1960), 43-76.

\vspace{1ex}

\noindent [Kur62]\, M. Kuranishi --- {\it On the Locally Complete Families of Complex Analytic Structures} --- Ann. of Math. {\bf 75}, no. 3 (1962), 536-577.

\vspace{1ex}

\noindent [KW70]\, S. Kobayashi, H.-H. Wu --- {\it On Holomorphic Sections of Certain Hermitian Vector Bundles} --- Math. Ann. {\bf 189} 1-4 (1970).

\vspace{1ex}

\noindent [LP92]\, C. Lebrun, Y.-S. Poon --- {\it Twistors, K\"ahler Manifolds, and Bimeromorphic Geometry. II} --- J. Amer. Math. Soc. {\bf 5}, No. 2 (1992) 317-325.

\vspace{1ex}

\noindent [Mic82]\, M. L. Michelsohn --- {\it On the Existence of Special Metrics in Complex Geometry} --- Acta Math. {\bf 149} (1982), no. 3-4, 261-295.

\vspace{1ex}

\noindent [Pop11]\, \, D. Popovici ---{\it Deformation Openness and Closedness of Various Classes of Compact Complex Manifolds; Examples} --- Ann. Sc. Norm. Super. Pisa Cl. Sci. (5), Vol. XIII (2014), 255-305.

\vspace{1ex}

\noindent [Pop15]\, D. Popovici --- {\it Aeppli Cohomology Classes Associated with Gauduchon Metrics on Compact Complex Manifolds} --- Bull. Soc. Math. France {\bf 143} (3), (2015), p. 1-37.

\vspace{1ex}

\noindent [Sch07]\, M. Schweitzer --- {\it Autour de la cohomologie de Bott-Chern} --- arXiv e-print math.AG/0709.3528v1.

\vspace{1ex}

\noindent [STW17]\, G. Sz\'ekelyhidi, V. Tosatti, B. Weinkove --- {\it Gauduchon Metrics with Prescribed Volume Form} --- Acta Math. {\bf 219} (2017), no.1, 181-211.

\vspace{1ex}

\noindent [Tia87]\, G. Tian --- {\it Smoothness of the Universal Deformation Space of Compact Calabi-Yau Manifolds and Its Petersson-Weil Metric} --- Mathematical Aspects of String Theory (San Diego, 1986), Adv. Ser. Math. Phys. 1, World Sci. Publishing, Singapore (1987), 629--646.

\vspace{1ex}

\noindent [Tod89]\, A. N. Todorov --- {\it The Weil-Petersson Geometry of the Moduli Space of $SU(n\geq 3)$ (Calabi-Yau) Manifolds I} --- Comm. Math. Phys. {\bf 126} (1989), 325-346.

\vspace{1ex}

\noindent [TW18]\, V. Tosatti, B. Weinkove --- {\it Hermitian Metrics, $(n-1, n-1)$-forms and Monge-Amp\`ere Equations} ---  arXiv e-print math.DG 1310.6326v2, to appear in J. Reine Angew. Math. 2018. 

\vspace{1ex}

\noindent [Voi02]\, C. Voisin --- {\it Hodge Theory and Complex Algebraic Geometry. I.} --- Cambridge Studies in Advanced Mathematics, 76, Cambridge University Press, Cambridge, 2002.         

\vspace{1ex}

\noindent [Yos01]\, K. Yoshioka --- {\it Moduli Spaces of Stable Sheaves on Abelian Surfaces} --- Math.Ann {\bf 321} (2001), 817-884.

\vspace{1ex}

\noindent [Wu06]\, C.-C. Wu --- {\it On the Geometry of Superstrings with Torsion} --- thesis, Department of Mathematics, Harvard University, Cambridge MA 02138, (April 2006).

\vspace{2ex}

\noindent Institut de Math\'ematiques de Toulouse, Universit\'e Paul Sabatier,

\noindent 118 route de Narbonne,  31 062 Toulouse Cedex 9, France

\noindent Email\!\!: popovici@math.univ-toulouse.fr

\end{document}